\documentclass[11pt,twoside, a4paper]{article}

\usepackage{latexsym}
\usepackage[english]{babel}
\usepackage{t1enc}
\usepackage{mathrsfs}
\usepackage{ifthen}
\usepackage{url}
\usepackage{epsfig}
\usepackage{amsbsy}
\usepackage{extarrows}
\usepackage{amsfonts}
\usepackage{amsmath}
\usepackage{amssymb}
\usepackage{amsthm}
\usepackage{amsxtra}
\usepackage{bbm}
\usepackage{color}
\usepackage{relsize} 
\usepackage{enumitem}
\usepackage{graphicx}
\usepackage{caption}
\usepackage{subcaption}
\usepackage{placeins}
\usepackage{wrapfig}

\usepackage[margin=2.5cm]{geometry}

\usepackage{verbatim}

\usepackage{hyperref} 

\numberwithin{equation}{section}
\numberwithin{figure}{section}

\newtheoremstyle{thm-style-oskari}
{7pt}      
{7pt}      
{\itshape} 
{}         
{\scshape} 
{.}        
{.5em}     
{}         

\theoremstyle{thm-style-oskari}
    \newtheorem{theorem}{Theorem}[section]
    \newtheorem{proposition}[theorem]{Proposition}
    \newtheorem{corollary}[theorem]{Corollary}
    \newtheorem{lemma}[theorem]{Lemma}
    \newtheorem{definition}[theorem]{Definition}
    
    \newtheorem{convention}[theorem]{Convention}
    \newtheorem{remark}[theorem]{Remark}

\newenvironment{Proof}[1][Proof]{\begin{proof}[\sc{#1}]}{\end{proof}}

\newcommand{\NTheorem}[2] {
        \begin{theorem}[#1] \label{thr:#1}
                #2
        \end{theorem}
        }

\newcommand{\NLemma}[2] {
        \begin{lemma}[#1] \label{lmm:#1}
                #2
        \end{lemma}
        }

\newcommand{\NCorollary}[2] {
        \begin{corollary}[#1] \label{crl:#1}
                #2
        \end{corollary}
        }

\newcommand{\NProposition}[2] {
        \begin{proposition}[#1] \label{prp:#1}
                #2
        \end{proposition}
        }

\newcommand{\NDefinition}[2] {
        \begin{definition}[#1] \label{def:#1}
                #2
        \end{definition}
        }

\newcommand{\NRemark}[2] {
        \begin{remark}[#1] \label{remark:#1}
                #2
        \end{remark}
        }

\newcommand{\bels}[2] {
        \begin{equation} \label{#1} \begin{split} 
                #2 
        \end{split} \end{equation}
        }
\newcommand{\bea}[1]{
	\begin{align*}
		#1
	\end{align*}
	}

\definecolor{olivegreen}{rgb}{0,0.6,0.1}
\newcommand{\nc}{\normalcolor}

\newcommand{\cob}{\color{blue}}

\newcommand{\Lp}[1]{\mathrm{L}^#1}					

\newcommand{\Ind}{\mathbbm{1}}

\newcommand{\vect}[1]{\mathbf{#1}}					

\newcommand{\1} {\mspace{1 mu}}
\newcommand{\2} {\mspace{2 mu}}
\newcommand{\msp}[1] {\mspace{#1 mu}}

\newcommand{\Err} {\Lambda}
\newcommand{\Bound}{\mcl{E}}

\newcommand{\la} {\langle}
\newcommand{\ra} {\rangle}
\newcommand{\avg}[1] {\la #1 \ra}
\newcommand{\avgb}[1] {\bigl\la #1 \bigr\ra}

\newcommand{\eps}{\varepsilon}

\newcommand{\dist} {\mrm{dist}}                   
\DeclareMathOperator*{\Spec}{Spec}						

\newcommand{\Id}{{\bf 1 }}          

\newcommand{\mrm}[1] {\mathrm{#1}}
\newcommand{\mcl}[1] {\mathcal{#1}}

\newcommand{\brm}[1] {\boldsymbol{\mathrm{#1}}}

\newcommand{\wti}[1] {\widetilde{#1}}
\newcommand{\wht}[1] {\widehat{#1}}
\newcommand{\ul}[1] {\underline{#1}}

\newcommand{\EE} {\mathbbm{E}}
\newcommand{\PP}  {\mathbbm{P}}

\newcommand{\diag} {\mrm{diag}}
\newcommand{\Tr} {\mrm{Tr}}

\DeclareMathOperator{\supp} {supp}

\newcommand{\ins} {\msp{1}\in\msp{1}}

\newcommand{\abs}[1]{\lvert #1 \rvert}
\newcommand{\absb}[1]{\big\lvert #1 \big\rvert}
\newcommand{\absB}[1]{\Bigl\lvert #1 \Bigr\rvert}
\newcommand{\absbb}[1]{\biggl\lvert #1 \biggr\rvert}

\newcommand{\norm}[1]{\lVert #1 \rVert}

\newcommand{\ceil}  [1] {      \lceil  {#1}         \rceil}

\newcommand{\R} {\mathbb{R}}
\newcommand{\C} {{\mathbb{C}}}

\newcommand{\N} {\mathbb{N}}
\newcommand{\Z} {\mathbb{Z}}

\newcommand{\Cp} {\mathbb{H}}
\newcommand{\G} {\mathbb{G}}

\newcommand{\D} {\mathbb{D}}

\newcommand{\sett}[1] { \{ {#1} \} }

\newcommand{\setb}[1] { \bigl\{ {#1} \bigl\} }
\newcommand{\setB}[1] { \Bigl\{ {#1} \Bigr\} }
\newcommand{\setbb}[1] { \biggl\{\, {#1} \,\biggr\} }

\newcommand{\titem}[1] {\item[\emph{(#1)}]} 

\newcommand{\dif} {\mathrm{d}}
\newcommand{\cI} {\mathrm{i}}
\newcommand{\nE} {\mathrm{e}}
\newcommand{\Ord} {\mathcal{O}}

\renewcommand{\Im}{\mathrm{Im}}
\renewcommand{\Re}{\mathrm{Re}}

\begin{document}
\renewcommand{\thefootnote}{\fnsymbol{footnote}}

\title{{\bf Universality for general Wigner-type matrices}
}
\author{
\begin{minipage}{0.3\textwidth}
 \begin{center}
Oskari H. Ajanki\footnotemark[1]\\
\footnotesize 
{IST Austria}\\
{\url{oskari.ajanki@iki.fi}}
\end{center}
\end{minipage}
\begin{minipage}{0.3\textwidth}
\begin{center}
L\'aszl\'o Erd{\H o}s\footnotemark[2]  \\
\footnotesize {IST Austria}\\
{\url{lerdos@ist.ac.at}}
\end{center}
\end{minipage}
\begin{minipage}{0.3\textwidth}
 \begin{center}
Torben Kr\"uger\footnotemark[3]\\
\footnotesize 
{IST Austria}\\
{\url{torben.krueger@ist.ac.at}}
\end{center}
\end{minipage}
}

\date{} 

\maketitle
\thispagestyle{empty} 

\footnotetext[1]{Partially supported by ERC Advanced Grant RANMAT No.\ 338804, and SFB-TR 12 Grant of the German Research Council.}
\footnotetext[2]{Partially supported by ERC Advanced Grant RANMAT No.\ 338804.}
\footnotetext[3]{Partially supported by ERC Advanced Grant RANMAT No.\ 338804, and SFB-TR 12 Grant of the German Research Council}
	
\renewcommand{\thefootnote}{\fnsymbol{footnote}}

\vspace{-0.8cm}
\begin{abstract}
We consider the local eigenvalue distribution of large self-adjoint $N\times N$ random matrices $\mathbf{H}=\mathbf{H}^*$ with centered independent entries. In contrast to previous works the matrix of variances $s_{i j}=\mathbbm{E}\mspace{2 mu} |h_{i j}|^2$ is not assumed to be stochastic. Hence the density of states is not the Wigner semicircle law. Its possible shapes are described in the companion paper \cite{AEK1}. 
We show that as $N$ grows, the resolvent, $\mathbf{G}(z)=(\mathbf{H}-z)^{-1}$, converges to a diagonal matrix, $\mathrm{diag}(\mathbf{m}(z))$, where $\mathbf{m}(z)=(m_1(z),\dots,m_N(z))$ solves the vector equation $-1/m_i(z) = z+\sum_j s_{i j}  \mspace{2 mu} m_j(z) $ that has been analyzed in \cite{AEK1}. We  prove a local law down to the smallest spectral resolution scale, and bulk universality for both real symmetric and complex hermitian symmetry classes. 
%
\end{abstract}
%
{\small {\bf Keywords:} Eigenvector delocalization, Rigidity, Anisotropic local law, Local spectral statistics\\
{\bf AMS Subject Classification (2010):} \texttt{60B20}, \texttt{15B52}}\\

\vspace{-0.5cm}
\tableofcontents

\section{Introduction}

In the seminal paper \cite{Wigner-ORIG} Wigner introduced random self-adjoint matrices, $\vect{H}=\vect{H}^*$,  with centered, identically distributed and independent entries (subject to the symmetry constraint). He proved that the empirical density of the eigenvalues converges to the semicircle distribution. Wigner also conjectured that the distribution of the distance between consecutive eigenvalues (\emph{gap statistics}) is universal, 
hence it is the same as in the Gaussian model. 
His revolutionary observation was that these universality phenomena hold for much larger classes of physical systems and only the basic symmetry type determines local spectral statistics. It is generally believed, but mathematically unproven, that random matrix theory (RMT), among many other examples, also describes the local  statistics of random Schr\"odinger operators in the delocalized regime and quantization of chaotic classical Hamiltonians. 

The first rigorous results on the local eigenvalue statistics in the bulk spectrum were given by Dyson, Mehta and Gaudin in the 60's. These concerned the Gaussian models and identified their local correlation functions. According to  Wigner's   universality hypothesis, these statistics should hold independently of the common law of the matrix elements. This conjecture was resolved recently in a series of works. The strongest result on Wigner matrices in the bulk spectrum is Theorem 7.2 in \cite{EKYY2}, see \cite{EYBull} and \cite{TV}  for a summary of the history and related results. In fact, the {\it three-step} approach developed in \cite{ESYY, EYY, EKYY} also applies for \emph{generalized Wigner matrices} that allow for non-identically distributed matrix elements as long as the variance matrix $s_{ij}:= \EE |h_{ij}|^2$ is stochastic, i.e. $\sum_j s_{ij}=1$ (in particular, independent of the index $i$). The stochasticity of $\vect{S}$ guarantees that the eigenvalue density is given by the semicircle law and the diagonal elements $G_{ii} = G_{ii}(z) $ of the resolvent 
\bels{def of G(z)}{
\brm{G}(z) \2=\2 (\vect{H}-z\1)^{-1}\,,
\qquad\Im \, z>0
\,, 
}
become not only deterministic but also  independent of $i$ as the 
the matrix size $ N $ goes to infinity.
They asymptotically satisfy a system of self-consistent equations
\bels{GSCE}{
-\frac{1}{ G_{ii}\!} \,\approx\, z+ \sum_{j} s_{ij}\2 G_{jj}
\,,
}
that reduces to a particularly simple scalar equation 
\bels{scemsc}{
-\frac{1}{m} =  z + m\,,
}
for the common value $m\approx G_{ii}$ for all $i$ as $N\to\infty$. 
\nc
The  solution $ m = m(z) $ of \eqref{scemsc} is the Stieltjes transform of the  Wigner semicircle law.

In this paper we consider a general variance matrix $\brm{S}$ without stochasticity condition. We show that the approximate self-consistent equation \eqref{GSCE} still holds, but it does not simplify to a scalar equation. 
In fact, $G_{ii}$ remains $i$-dependent even as $N\to\infty$ and it is close to the solution $m_i $ of the \emph{Quadratic Vector Equation} (QVE)
\bels{qev}{
-\frac{1}{ m_i} = z+ \sum_{j} s_{ij} \2m_j
\,,
}
under the additional condition that $\Im \, m_i > 0$.

In the context of random matrices importance of this equation has been realized by Girko \cite{Girko-book}, Shlyakhtenko \cite{ShlyakhtenkoGBM}, Khorunzhy and Pastur \cite{KhorunzhyPastur94}, see also Guionnet \cite{Guionnet-GaussBand}, as well as Anderson and Zeitouni \cite{AZind,AZdep}, but no detailed study has been initiated. 
In the companion paper \cite{AEK1} we analyzed \eqref{qev} in full detail. See also Section~3 of \cite{AEK1cpam} for how the QVE is related to other random matrix models. 
We showed that $\avg{\1m} := N^{-1}\sum_i m_i $ is the Stieltjes transform of a probability density $\rho $ that is supported on a finite number of intervals, inside of which it is a real analytic function.
We also described  the behavior of $ \rho $ near the edges of its support; it features only square root or cubic root (cusp) singularities and an explicit one parameter family of profiles interpolating between them as a gap in the support closes.
   
The main result of the current paper is the  universality of the local eigenvalue statistics in the bulk for Wigner-type matrices with a general variance matrix (cf. Theorem \ref{thr:Universality}). This extends Wigner's vision towards full universality by considering a much larger class of matrix ensembles than previously studied. In particular, we demonstrate  that local statistics, as expected, are fully independent of the global density. This fact has already been established for very general 
 $\beta$-ensembles in \cite{BEY} (see also \cite{BFG} and \cite{MSch}) and for additively deformed Wigner ensembles having a density with a single interval support  \cite{LSSY2014}. Our class admits a  general variance matrix and allows for densities with several intervals (we do not, however, consider non-centered distributions here; an extension to matrices with non-centered entries on the diagonal may be incorporated in our analysis with additional technical effort). 
 
Following the three-step approach, we first   prove \emph{local laws} for  $\brm{G}$ on the scale $\eta=\Im \, z \gg N^{-1}$, i.e. down to the optimal scale just slightly above the eigenvalue spacing. 
This is the main and novel part of our analysis. The previous proofs (see \cite{EKYY} for a pedagogical presentation) heavily relied on properties of the semicircle law, especially on its square root edge singularity. With possible cubic root singularities and small gaps in the support of $\rho$ an additional scale appears which needs to be controlled. 
The second step is to prove universality for Wigner-type matrices with a tiny Gaussian component via {\it Dyson Brownian motion} (DBM). The method of {\it local relaxation flow}, introduced first in \cite{ESY, ESYY}, also heavily relies on the semicircle law since it requires that the global density remain unchanged along DBM. In \cite{EK} a new method was developed to localize the DBM that proves universality of the gap statistics around a fixed energy $\tau$ in the bulk, assuming that the local law holds near $\tau$ (we remark that a similar result was obtained independently in \cite{LY}). Since Wigner-type matrices were one of the main motivations for  \cite{EK},  it was formulated such that it could be directly applied once the local laws  are available. Finally, the third step is a perturbation result to remove the tiny Gaussian component using the {\it Green function comparison} method that first appeared in \cite{EYY} and can be applied to our case basically without any modifications.

In a separate paper \cite{AEK3} we  apply the results of this work and \cite{AEK1} to treat Gaussian random matrices with correlated entries. Assuming translation invariance of the correlation structure in these Gaussian matrix ensembles we prove an optimal local law, bulk universality and non-trivial decay of off-diagonal resolvent entries.

\medskip
\emph{Acknowledgement.}
We thank the anonymous referee for several useful comments and suggestions. We are grateful to Johannes Alt for pointing out several typos.

\subsection{Set-up and main results}
\label{ssec:Set-up and main results}

Let $\vect{H}^{(N)}\in \C^{N \times N}$ be a sequence of  self-adjoint random matrices.  In particular, if the entries of $\vect{H}$ are real then $\vect{H}^{(N)}$ is symmetric. 
The matrix ensemble $\vect{H}=\vect{H}^{(N)}$ is said to be of {\bf Wigner-type} if  its entries $h_{i j}$ are independent for 
$i \leq j$ and centered, i.e.,
\bels{EE h^2 = s}{
\EE \2 h_{ij} \,=\, 0 \quad\text{ for all } \quad i,j \,=\, 1,\dots,N\,.
}
The dependence of $\vect{H}$ and other quantities on the dimension $N$ will be suppressed in our notation.
The matrix of variances, $\vect{S}=(s_{i j})_{i,j=1}^N$, is defined through  
\bels{definition of S}{
s_{ij} := \EE \2 |h_{i j}|^2\,.
}
It is symmetric with non-negative entries.
In \cite{AEK1} it was shown that for every such matrix $\vect{S}$ the 
{\bf quadratic vector equation} (QVE), 
\bels{SCE}{
-\, \frac{1}{m_i(z)} \,=\, z+\sum_{j=1}^N s_{ij}\1 m_j(z)\,, 
\quad \text{for all } \; i=1,\dots,N \text{ and } z \in \Cp\,,
}
for a function $\vect{m}=(m_1,\dots,m_N) : \Cp \to \Cp^N$ on the complex upper half plane, $\Cp=\{z\in\C: \Im \1 z >0  \}$, has a unique solution.
The main result of this paper is to establish the local law for Wigner-type matrices, i.e. that for large $N$ the resolvent, 
$\vect{G}(z) = (\vect{H}-z)^{-1}$, with spectral parameter $z= \tau +\cI \1\eta \in \Cp$, is close to the diagonal matrix, $\diag(\vect{m}(z))$, as long as $\eta \gg N^{-1}$. As a consequence, we obtain rigidity estimates on the eigenvalues and complete delocalization for the eigenvectors. Combining this information with the Dyson-Brownian motion, we obtain universality of the eigenvalue gap statistics in the bulk.

We now list the assumptions on the variance matrices $\vect{S}=\vect{S}^{(N)}$.
The restrictions on $\vect{S}$ are controlled by three {\bf model parameters}, $p, P >0$ and $L \in \N$, which do not depend on $N$. These parameters will remain fixed throughout this paper. 
\begin{enumerate}
\item[(A)] For all $N$ the matrix $\vect{S}$ is {\bf flat}, i.e.,
\bels{flat}{
s_{i j} \;\leq\; \frac{1}{N}\,, \qquad  i,j=1,\dots,N\,.
}
\item[(B)] For all $N$ the matrix $\vect{S}$ is {\bf uniformly primitive}, i.e.,
\bels{uniform primitivity}{
(\1\brm{S}^{L})_{i j} \,\geq\, \frac{p}{N}\,, \qquad  i,j=1,\dots,N\,.
}
\item[(C)] For all $N$ the matrix $\vect{S}$ induces a {\bf bounded solution} of the QVE, i.e., the unique solution $\vect{m}$ of \eqref{SCE} corresponding to $\vect{S}$ is bounded,
\bels{boundedness of m}{
|m_i(z)|\;\leq\; P\, , \qquad  i=1,\dots,N\,, \quad z \in \Cp\,.
\,.
}
\end{enumerate}
%
\NRemark{Boundedness and normalization}{
The assumption on the boundedness of $\vect{m}$ is an implicit condition in the sense that it can be checked only after solving \eqref{SCE}. In Theorem~6.1 of  \cite{AEK1} we list sufficient, explicitly checkable conditions on $\vect{S}$, which ensure \eqref{boundedness of m}. We also remark that the assumption \eqref{flat} can be replaced by $s_{ij} \leq C/N $ for some positive constant $C$. This will lead to a rescaling (cf. Remark~2.2 of \cite{AEK1}) of $\vect{m}$. We pick the normalization $ C=1 $ just for convenience. 
}

\NRemark{Primitivity}{
The primitivity condition \eqref{uniform primitivity} excludes some
important models, e.g. matrices of the form
\[
\brm{H} = \begin{pmatrix} \brm{0} & \brm{X} \cr \brm{X^*} &\brm{0} \end{pmatrix}
\,,
\]
whose eigenvalues  yield  the singular values of the Gram matrix $\brm{X}\brm{X}^*$, where $\brm{X}$ 
has independent centered entries with an arbitrary variance profile.
Condition \emph{(B)} is not a mere technicality; Gram matrices may have singularities in the spectrum near 0 (often referred to as the 'hard-edge') that require separate treatment; but even away from 0 some new ideas are needed. The complete analysis is presented in  \cite{AltEK}, where we prove local laws for Gram matrices.
}

In addition to the assumptions on the variances of $h_{i j}$, we also require uniform 
boundedness of higher moments. This leads to another basic model parameter, 
$\ul{\mu}=(\mu_1, \mu_2, \dots )$, which is a sequence of non-negative real numbers.

\begin{itemize}
\item[(D)] For all $N$ the entries $h_{i j}$ of the random matrix $\vect{H}$ have {\bf bounded moments},
\bels{bounded moments}{
\EE \2 | h_{i j} |^k \,\leq\,  \mu_k \1s_{i j}^{\1k/2}\,,
\qquad k\in \N\,,\;i,j=1, \dots, N\,.
}
\end{itemize}

In order to state our main result, in the next corollary we collect a few facts about the solution of the QVE that are proven in \cite{AEK1}. Although these properties are sufficient for the formulation of our results, for their proofs we will need much more precise information about the solution of the QVE. Theorems \ref{thr:Solution of the QVE} and \ref{thr:Stability around support} summarize everything that is needed from \cite{AEK1} besides the existence and uniqueness of the solution of the QVE. In particular, the statement of Corollary \ref{crl:Solution of QVE} follows easily from Theorem \ref{thr:Solution of the QVE}  below. 

\NCorollary{Solution of QVE}{ 
Suppose  $ \brm{S} $ satisfies \emph{(A)--(C)}.
Let $\vect{m} : \Cp \to \Cp^N$ be the solution the QVE \eqref{SCE} corresponding to $\vect{S}$. Then $\vect{m}$ is analytic and has a continuous extension (denoted again by $\vect{m}$) to the closed upper half plane,
$\vect{m}: \overline \Cp \to \overline \Cp^N$, with $\overline \Cp:=\Cp \cup \R$. The function $\rho: \R \to [0,\infty)$, defined by
\bels{def of rho}{
\rho(\tau)\,:=\, \lim_{\eta \downarrow 0}\frac{1}{\pi N}\sum_{i=1}^N \Im\2 m_i(\tau+\cI\1\eta\1)\,,
}
is a probability density. Its support is contained in $[-2,2]$ and is a union of closed disjoint intervals
\bels{DOS support}{
\supp \rho\,=\, \bigcup_{k=1}^K[\1\alpha_k,\beta_k]\,,\quad \text{ where }\quad
\alpha_k<\beta_k< \alpha_{k+1}\,.
}
There exists a positive constant $\delta_*$, depending only on the model parameters $p$, $P$ and $L$, such that the sizes of these intervals are bounded from below by
\bels{lower bound on interval size}{
\beta_k-\alpha_k\,\geq\, 2\2\delta_*\,.
}
}

\begin{wrapfigure}{r}{0.6\textwidth}
	\centering
	\vspace{-0.5cm}
	\includegraphics[width=0.55\textwidth]{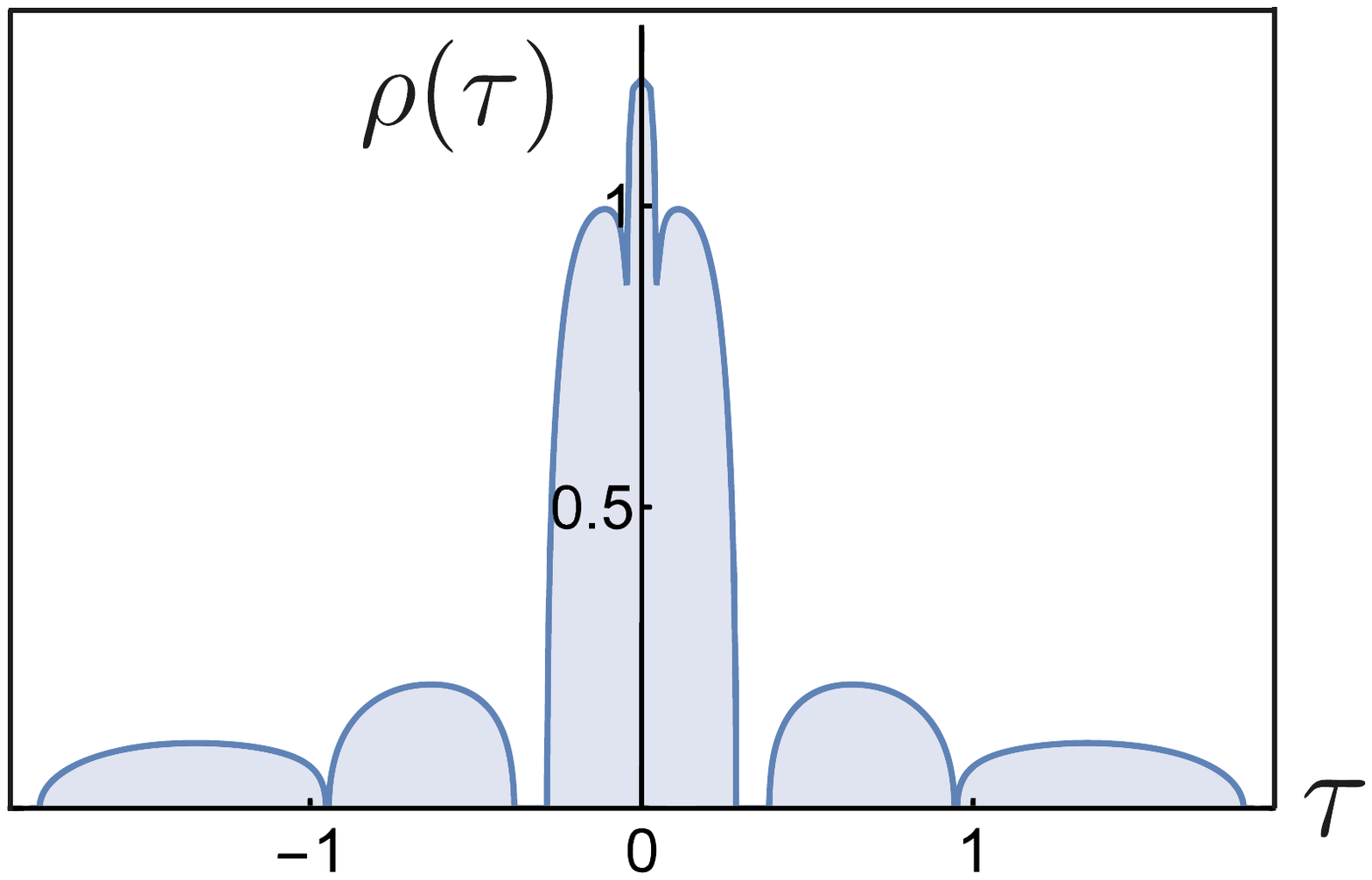}
	\caption{The density of states may have gaps, cusps and local minima. It is always a symmetric function around zero, i.e., $\rho(\tau)=\rho(-\tau)$.}
	\label{DOSFigure}
	\vspace{-0.5cm}	
\end{wrapfigure}
Note that \eqref{lower bound on interval size} provides a lower bound on the length of the intervals that constitute $\supp \rho$, while the length of the gaps, $\alpha_{k+1}-\beta_k$, between neighboring intervals can be arbitrarily small. 
Figure \ref{DOSFigure} shows a shape that the density of states typically might have. In particular, $\rho$ may have gaps in its support and may have additional zeros (cusps) in the interior of $\supp\rho$.
However, the behavior of $ \rho $ on the domain $ \rho \leq \eps $, for some sufficiently small $ \eps > 0 $, can be completely characterized by universal shape functions. For more details  see Theorem~2.6 of \cite{AEK1}.

\NDefinition{Density of states}{
The function $\rho$ defined in \eqref{def of rho} is called the {\bf density of states}. Its harmonic extension to the upper half plane
\bels{harmonic extension of rho}{
\rho(\1\tau+ \cI\1\eta\1)\,:=\, \int_\R \Pi_\eta(\tau-\sigma)\2\rho(\sigma)\2\dif \sigma
\,,\qquad
\Pi_\eta(\tau)\,:=\, \frac{1}{\pi}\2\frac{\eta}{\tau^2+\eta^{\12}\!}\,; \qquad \tau \in \R\,,\;\eta \1>\10\,,
}
is again denoted by $\rho$. With a slight abuse of notation we still write $\supp \rho$, as in \eqref{DOS support}, for the support of the density of states as a function on the real line.
}

The density of states will be shown to be the eigenvalue distribution of $\vect{H}$ in the large $N$ limit on the macroscopic scale. For any fixed values $\tau_1,\tau_2 \in \R$ with $\tau_1<\tau_2$ it satisfies 

\bels{interpretation DOS}{ 
\lim_{N\to \infty} \frac{\big|\Spec(\2\vect{H}^{(N)})\cap [\1\tau_1,\tau_2\1]\,\big|}{N\int_{\tau_1}^{\tau_2} \!\rho^{(N)}(\tau)\,\dif \tau}\,=\, 1
\,,
}
provided the denominator does not vanish in the limit. The identity \eqref{interpretation DOS} motivates the terminology of density of states for the function $\rho$. The harmonic extension of $\rho$ to $\Cp$ is a version of the density of states, that is smoothed out on the scale $\eta$. It satisfies the identity $\rho(z)= (\pi\1N)^{-1}\sum_k \Im \2m_k(z)$ not just for $z \in \R$ (cf. \eqref{def of rho}) but for all $z \in \overline\Cp$ and it will be used in the statement of our main result, Theorem \ref{thr:Local law}.

In fact, Theorem \ref{thr:Local law}, implies a local version of \eqref{interpretation DOS}, where the length of the interval, $[\tau_1,\tau_2]$, may converge to zero as $N$ tends to infinity. Our estimates and thus the proven speed of convergence depend on the distance of the interval to the edges of $\supp \rho$ and even on the length of the closest gap in this case. We introduce a function $\Delta: \R \to [0,\infty)$, which encodes this relation.

\NDefinition{Local gap size}{
Let $\alpha_k$ and $\beta_k$ be the {\bf edges} of the support of the density of states (cf. \eqref{DOS support}) and $\delta_*$ the constant introduced in Corollary \ref{crl:Solution of QVE}. Then for any $\delta \in [0,\delta_*)$ we set
\bels{def of Delta_delta}{
\Delta_\delta(\tau)\,:=\, 
\begin{cases}
\alpha_{k+1}-\beta_k
& 	
\text{if $\beta_k-\delta \leq \tau \leq \alpha_{k+1}+\delta $ for some $k=1,\dots,K-1$,}
\\
\;1 				
& 	
\text{if $\tau \leq \alpha_1+\delta$ or $\tau \geq \beta_{K}-\delta$,} 
\\
\;0
& 	
\text{ otherwise.}
\end{cases}
}
}

Finally, we fix an arbitrarily small {\bf tolerance exponent} $\gamma \in (0,1)$.
This number will stay fixed throughout this paper in the same fashion as the model parameters $P$, $p$, $L$ and $\ul{\mu}$.
Our main result is stated for spectral parameters $z=\tau+\cI\1\eta$ whose imaginary parts satisfy
\bels{}{
\eta\,\geq\, N^{\gamma-1}\,.
} 

For a compact statement of the main theorem we define the notion of stochastic domination, introduced in  \cite{EKY} and \cite{EKYY}. This notion is designed to compare sequences of random variables in the 
large $N$ limit up to small powers of $N$ on high probability sets. 

\NDefinition{Stochastic domination}{
Suppose $N_0: (0,\infty)^2\to \N$ is a given function, depending only on the model parameters $p$, $P$, $L$ and $\ul{\mu}$, as well as on the tolerance exponent $\gamma$. For two sequences, $\varphi=(\varphi^{(N)})_N$ and $\psi=(\psi^{(N)})_N$, of non-negative random variables we say that $\varphi$ is {\bf stochastically dominated} by $\psi$ if for all $\eps>0$ and $D>0$,
\bels{}{
\PP\Bigl( \1\varphi^{(N)} \1>\1 N^\eps \1 \psi^{(N)}\Bigr) 
\;\leq\; N^{-D},\qquad N\2\geq\2 N_0(\eps,D)\,.
}
In this case we write $\varphi \prec \psi$. 
}
Basic properties of the stochastic domination that are used extensively in this paper are listed in Lemma~\ref{lmm:Basic facts about stochastic domination}. 
The threshold $ N_0(\eps,D) = N_0(\eps,D;P,p,L,\ul{\mu},\gamma) $ will always be an explicit function whose value will be increased throughout the paper, though we will not follow its form. This will happen only finitely many times, ensuring that $N_0$ stays finite. 
The threshold is uniform in all other parameters, e.g. in the spectral parameter $z$, as well as in the indices $ i,j, \dots$ of the  matrix entries, that the sequences $\varphi$ and $\psi$ may depend on. Typically, we will not mention the existence of $N_0$ - it is implicit in the notation $\varphi\prec \psi$.  As an example, we see that the bounded moment condition, (D), implies 
\[
\abs{\1h_{ij}} \2\prec\2 N^{-1/2}
.
\]

Actually, the function $N_0$ depends only on finitely many moment parameters $(\mu_1,\dots, \mu_{M})$ instead of the entire sequence $\ul{\mu}$, where now the number of required moments $M$ $=M(\eps,D;P,p,L,\gamma)$, is an explicit function.

Now we are ready to state our main result on the local law. Suppose $\vect{H}=\vect{H}^{(N)}$ is a sequence of self-adjoint random matrices with the corresponding sequence $\vect{S}=\vect{S}^{(N)}$ of variance matrices and $\rho=\rho^{(N)}$ the induced sequence of densities of state. Recall that $\delta_*$ is the positive constant, depending only on $p$, $P$ and $L$, introduced in Corollary \ref{crl:Solution of QVE} and $\Delta_\delta$ is defined as in Definition~\ref{def:Local gap size}.

\begin{theorem}[Local law]
\label{thr:Local law}
Suppose that assumptions \emph{(A)-(D)} are satisfied and fix an arbitrary $\gamma \in (0,1)$.  
There is a deterministic function $\kappa=\kappa^{(N)}: \Cp \to (0,\infty]$ such that 
uniformly for all $z=\tau +\cI\1\eta\in \Cp$ with $\eta \geq N^{\gamma-1}$ the resolvents \eqref{def of G(z)} of the random matrices $\vect{H}=\vect{H}^{(N)}$ satisfy
\bels{pointwise local law}{
\max_{i,j}|\1G_{ij}(z)-m_i(z)\1\delta_{ij}\1| 
\,\prec\,
\sqrt{\frac{\rho(z)}{N\1\eta}}
+\frac{1}{N\1\eta} 
+ 
\min\setbb{\!\frac{1}{\sqrt{N\1\eta}}\,,\2 \frac{\kappa(z)}{N\1\eta}\!}
\,.
}
Furthermore, for any sequence of deterministic vectors $\vect{w} = \vect{w}^{(N)} \in \C^N$ with $ \max_{i}|w_i|\leq 1 $ the averaged resolvent diagonal has an improved convergence rate,
\bels{averaged local law}{
\absbb{\frac{1}{N}\sum_{i=1}^N \overline{w}_i\,\bigl(G_{ii}(z)-m_i(z)\bigr)}
\,\prec\, 
\min\setbb{\!\frac{1}{\sqrt{N\1\eta}}\,,\2 \frac{\kappa(z)}{N\eta}\!}
\,.
}
In particular, for $w_i=1$ this implies
\bels{}{
\absbb{\frac{1}{N}\2\Im\2\Tr \,\vect{G}(z)-\pi\2\rho(z)}
\,\prec\, 
\min\setbb{\!\frac{1}{\sqrt{N\1\eta}}\,,\2 \frac{\kappa(z)}{N\eta}\!}
\,.
}
The function $ \kappa $ may be chosen to be
\bels{kappa bound}{
\kappa(z)\,=\,
\frac{1}{\Delta(\tau)^{1/3}+\rho(z)}
\,,
}
where $\Delta=\Delta_\delta $, with some $\delta \in (0,\delta_*)$ that depends only on the model parameters $p$, $P$ and $L$.
 
In the regime, where $ z $ is not too close to the support of the density of states in the sense that
\bels{condition for improved kappa bound}{
(\2\Delta(\tau)^{1/3}\!+ \rho(z)\1)\,\dist(\1z, \supp \rho\2)
\;\ge\,
\frac{N^\gamma}{(N\1\eta)^2}\,,
}
$\kappa$ maybe improved to
\bels{improved kappa bound}{
\kappa(z)
\;=\: 
\frac{\eta}{\dist(z, \supp \rho\1)\,(\1\Delta(\tau)^{1/3}\!+\rho(z))}
\,+\,
\frac{1}{N\eta\2\dist(z, \supp \rho\1)^{1/2}\2(\1\Delta(\tau)^{1/3}+ \rho(z))^{1/2}\!}
\;.
}
\end{theorem}

The size of $ \rho(z) $ is described in \eqref{qualitative size of rho} below.
Theorem \ref{thr:Local law} can be localized to a spectral interval $ I \subset \R $, i.e., the statements hold for $ \Re\,z \in I $ provided \eqref{boundedness of m} applies for $ z \in I + \cI\1(0,\infty) $.
In particular, in the bulk of the spectrum Theorem~\ref{thr:Local law} simplifies considerably.
\begin{corollary}[Local law in the bulk]
\label{crl:Local law in the bulk}
Assume \emph{(A)}, \emph{(B)} and \emph{(D)} with $ L = 1 $.
Suppose there is a constant $ \rho_\ast > 0 $ and an interval $ I \subset \supp\,\rho $ such that $  \rho(\tau) \ge \rho_\ast $ for all $ \tau \in I$. 
Then uniformly for all $ z = \tau + \cI\1\eta $, with $ \tau \in I $ and $ \eta \ge N^{\gamma-1} $, and non-random $ \brm{w} \in \C^N $ satisfying $ \max_i \abs{w_i}\leq 1 $, the local laws hold
\bels{local law in the bulk}{
\max_{i,\2j=1}^N\,\abs{\1G_{ij}(z)-m_i(z)\1\delta_{ij}\1} 
\,\prec\,
\frac{1}{\!\sqrt{N\1\eta}}\,,
\quad\text{and }\quad
\absbb{\frac{1}{N}\sum_{i=1}^N \overline{w}_i\,\bigl(G_{ii}(z)-m_i(z)\bigr)}
\,&\prec\, 
\frac{1}{N\1\eta}
\,,
}
where $ \rho_\ast $ is considered as an additional model parameter.
\end{corollary}

Here the additional assumption $ L = 1 $ is only used to guarantee (cf. (i) of Theorem~6.1 in \cite{AEK1}) that the solution $ \brm{m}(z) $ of the QVE stays bounded around $ z = 0 $. 
Indeed, if $ z $ is bounded away from zero then (i) of Lemma~5.4 in \cite{AEK1} implies $ \norm{\brm{m}}_\infty := \max_i\abs{\1m_i}$ is bounded by a constant independent of $ N $ in the bulk of the spectrum. 
Therefore, if $ \dist(I,0) \ge \delta $ for some $ \delta > 0 $, or $ \sup\sett{\,\norm{\1\brm{m}(z)}_\infty:\Re\,z \in I\2} \leq P $ is known for some $ P < \infty $, then the assumption $ L=1 $ can be removed, and \eqref{local law in the bulk} holds with $ \delta $ or $ P $, respectively, considered as model parameters.

Theorem~\ref{thr:Local law} generalizes the previous local laws for stochastic variance matrices $\vect{S}$ (see \cite{EKYY} and references therein). It is valid for densities $\rho$ with an edge behavior different from the square root growth that is known from Wigner's semicircular law. In particular, singularities that interpolate between a square root  and a cubic root are possible. In the bulk of the support of the density of states, i.e., where $\rho$ is bounded away from zero, the function $\kappa$ is bounded. The same is true near the edges, unless the nearby gap is small. The bound deteriorates near small gaps in the support of $\rho$. 

In applications, the sequence $\vect{S}=\vect{S}^{(N)}$ satisfying (A)--(C) may be constructed by discretizing a piecewise $ 1/2$-H\"older \nc continuous limit function (cf. Remark~6.2 in \cite{AEK1}).
As a particularly simple example, suppose $f $ is a smooth, non-negative, symmetric, $f(x,y)=f(y,x)$, function on $[0,1]^2$ with a positive diagonal, $f(x,x)>0$.  Then the sequence of variance matrices,
\[
s_{i j}^{(N)}\,:=\, \frac{1}{N} f\Bigl(\frac{i}{N},\frac{j}{N}\Bigr)\,,\qquad i,j = 1,\dots,N\,,
\]
satisfies conditions (A)--(C). The validity of (C) can be verified by using the general criteria (cf. Theorem~2.10 and Theorem~6.1 of \cite{AEK1}) for uniform boundedness.  
In this case the solution, $\vect{m}=(m_1,\dots,m_N)$, of the QVE converges to a limit in the sense that
\[
\sup_{z\in\Cp} \max_{i=1}^N \absb{m_i(z)-m(\1i/N;z)}\,\to\, 0\,,
\] 
where $m:[0,1]\times \overline{\Cp} \to  \overline{\Cp}$ is the solution of the continuous QVE,
\[
-\frac{1}{m(x;z)}\,=\, z + \int_0^1 \!f(x,y)\1m(y;z)\2\dif y\,,\qquad x \in [0,1]\,,\; z \in \overline{\Cp}
\,.
\]
The continuous QVE such as this one fall into the class of general QVEs thoroughly analyzed in the companion paper \cite{AEK1}. In particular, the stability analysis applies and the density of states converges to a limit
\[
\rho^{(N)}(\tau)\,\to\, \frac{1}{\pi}\int_0^1\Im \2m(x;\tau)\1\dif x\,.
\]

We introduce a notion for expressing that events hold with high probability in the limit as $N$ tends to infinity.

\NDefinition{Overwhelming probability}{
Suppose $N_0: (0,\infty)\to \N$ is a given function, depending only on the model parameters $p$, $P$, $L$ and $\ul{\mu}$, as well as on the tolerance exponent $\gamma$.
For a sequence $A=(A^{(N)})_N$ of random events we say that $A$ hold {\bf asymptotically with overwhelming probability} (a.w.o.p.), if for all $D>0$:
\bels{}{
\PP(A^{(N)})\,\geq\, 1-N^{-D},\qquad N\2\geq\2 N_0(D)\,.
}
}

There is a simple connection between the notions of stochastic domination and asymptotically overwhelming probability.
 For two sequences $A = A^{(N)}$
and $B = B^{(N)}$ the statement '$A$ implies $B$ a.w.o.p.' is equivalent to $\Ind_A \prec \Ind_B$, where the threshold $N_0$, implicit in the stochastic domination, does not depend on $\eps$, i.e., $N_0(\eps,D)=N_0(D)$.

We denote by $\lambda_1\leq \dots\leq \lambda_N$ the eigenvalues of the random matrix $\vect{H}$. The following corollary shows that the eigenvalue distribution converges to the density of states as $N$ tends to infinity. 

\begin{corollary}[Convergence of cumulative eigenvalue distribution]
\label{crl:Convergence of cumulative eigenvalue distribution}
Assume \emph{(A)-(D)}. 
Then uniformly for all $\tau \in \R$ the cumulative empirical eigenvalue distribution approaches the integrated density of states,
\bels{difference between distributions in the spectrum}{
\absbb{\2\#\{\2i: \lambda_i \2\leq\2 \tau\}- N\!\int_{-\infty}^{\tau}\msp{-10}\rho(\omega)\2\dif \omega\,}
\,\prec\, 
\min\setbb{
\frac{1}{\Delta(\tau)^{1/3}\!+\rho(\tau)}
\2,\,
N^{1/5} 
}
\,.
}
Furthermore, for an arbitrary tolerance exponent $\gamma\in(0,1)$ there are no eigenvalues away from the support of the density of states,
\bels{gaps in spectrum}{
\max_{k=0}^K\,\#\setb{\,i: \beta_k+\delta_k\2<\2 \lambda_i \2<\2\alpha_{k+1}\!-\delta_k}
\;=\,0
\qquad \text{a.w.o.p.}
\,,
}
where we interpret $\beta_0:=-\infty$, $\alpha_{K+1}:=+\infty$ and $\delta_k$ is defined as $\delta_0:=\delta_K:=N^{\gamma-2/3}$, as well as
\bels{definition of deltak}{
\delta_k\,:=\,\frac{N^{\gamma}}{(\alpha_{k+1}\!-\beta_k)^{1/3}N^{2/3}\!}\;, \qquad k=1,\dots,K-1\,.
}
\end{corollary}

Based on \eqref{interpretation DOS} we define the index, $ i(\tau) $, of an eigenvalue that we expect to be located close to the spectral parameter $ \tau $ by
\nc
\bels{def of i(tau)}{
i(\tau)\,:=\, \bigg\lceil  N \!\int_{-\infty}^\tau\msp{-10}\rho(\omega)\2\dif \omega\bigg\rceil
\,.
}
Here, $\lceil\omega\rceil$ denotes the smallest integer  that is bigger or equal to $\omega$ for any $\omega \in \R$. 

\begin{corollary}[Rigidity of eigenvalues]
\label{crl:Rigidity of eigenvalues}
Assume \emph{(A)-(D)}, and let $\gamma\in(0,1)$ be an arbitrary tolerance exponent. Denote
\bels{definition of epsk}{
\eps_k \,:=\, 
N^{\gamma}\min\setbb{  
\frac{1}{N^{3/5}\!}\;,
\frac{1}{(\alpha_{k+1}-\beta_k)^{1/9}N^{2/3}\msp{-7}}\,}
\,,\qquad k\2=\2 1,\dots,K-1\,,
}
and $\eps_0:=\eps_K:=N^{\gamma-2/3} $.
Then uniformly for every
\bels{tau well inside the support}{
\tau \,\in\, \bigcup_{k=1}^K\bigl[\1\alpha_k + \eps_{k-1}\1,\2\beta_k-\eps_k\bigr] \,,
}
the eigenvalues satisfy the rigidity
\bels{rigidity of eigenvalues in support}{
\abs{\1\lambda_{i(\tau)}-\2\tau\2}
\,\prec\, 
\min\setbb{
\frac{1}{(\1\Delta(\tau)^{1/3}+\rho(\tau)\1)\2\rho(\tau)\1N\2}\,,\frac{1}{\2N^{\13/5}\msp{-7}
}\,}
\,.
}
Furthermore, if $ \tau $ is close to the extreme edge, $ \tau \in( \alpha_1,\alpha_1 +\eps_0) $ or  $ \tau \in (\beta_K-\eps_K,\beta_K\1] $, then 
\bels{extreme edge fluctuation}{
\abs{\1\lambda_{i(\tau)}-\2\tau\2} \,\prec\, N^{-2/3}
.
}
Finally, if $\tau \in (\beta_k-\eps_k,\alpha_{k+1}+\eps_k)$ for some $ 1 \leq k \leq K-1$, then the corresponding eigenvalue is close to an internal edge in the sense that
\bels{edge eigenvalue rigidity}{
\lambda_{i(\tau)}\,\in\, \bigl[\1\beta_k-2\1\eps_k\1,\2\beta_k+\delta_k\bigr]
\cup 
\bigl[\1\alpha_{k+1}-\delta_k\1,\2\alpha_{k+1}+2\1\eps_k\bigr]
\qquad 
\text{a.w.o.p.}
\,,
}
where $ \delta_k $ is defined in \eqref{definition of deltak}.
\end{corollary}

\begin{wrapfigure}{r}{0.5\textwidth}
\centering
\vspace{-0.6cm}
\includegraphics[width=0.5\textwidth]{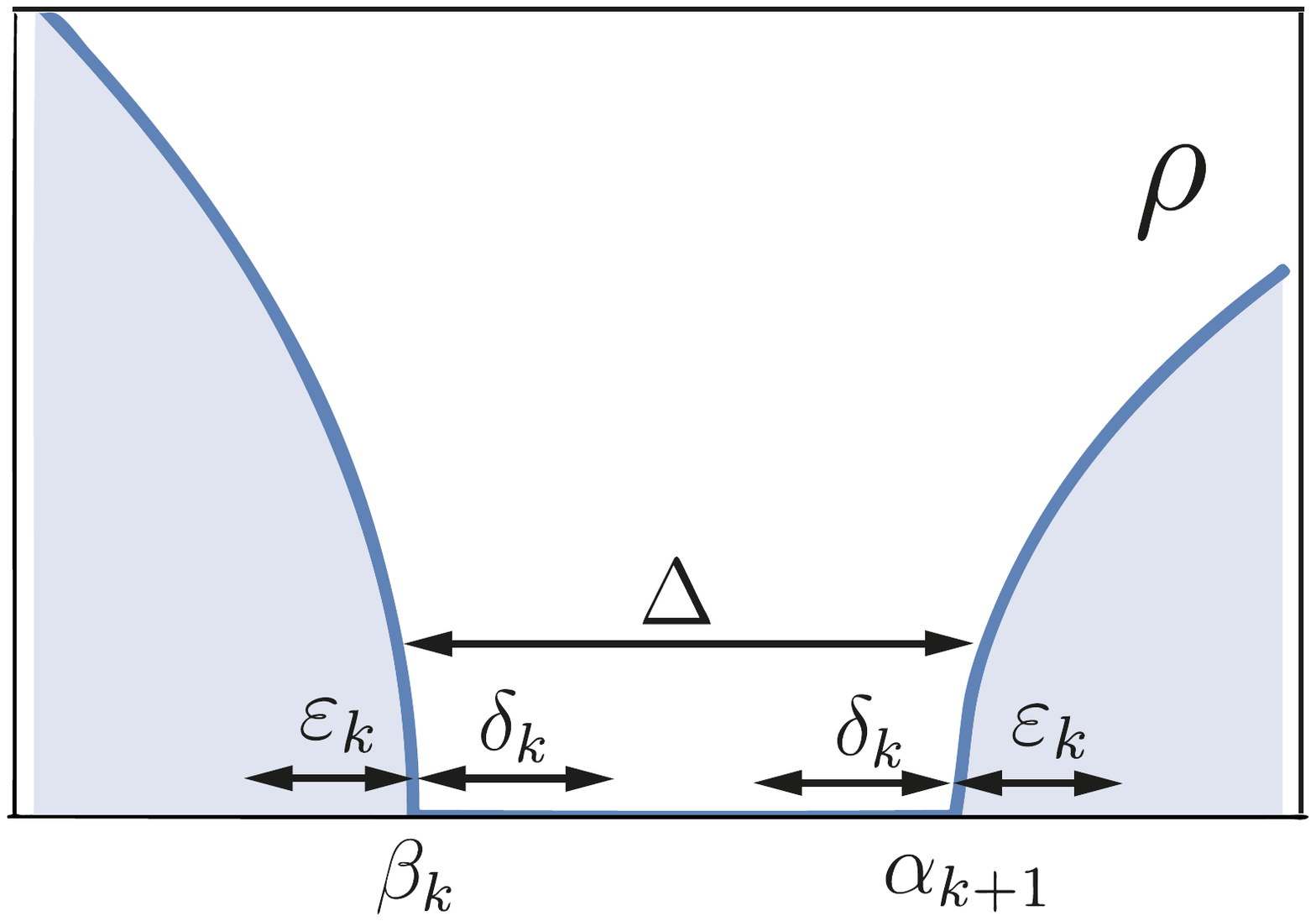}
\caption{
Notations of Corollary~\ref{crl:Rigidity of eigenvalues}: 
At the edges of a gap of length $ \Delta $ in $ \supp \rho$ the bound on the eigenvalue fluctuation is $ \delta_k $ inside the gap and $\eps_k$ inside the support.
}
\label{GapNotationFigure}
\vspace{-0.2cm}
\end{wrapfigure}

\begin{remark}[Eigenvalues outside $ \supp \rho $] 
The statements \eqref{extreme edge fluctuation} and \eqref{edge eigenvalue rigidity}
are an immediate consequence of \eqref{rigidity of eigenvalues in support} and \eqref{gaps in spectrum}. They simply express the fact that the small number of $\Ord(N^\eps)$ eigenvalues, very close to the edges, are found in the space that is left for them by the other eigenvalues for which the rigidity statement  \eqref{rigidity of eigenvalues in support} applies. 
For an illustration see Figure \ref{GapNotationFigure}. 
We also note that results of this type date back to at least \cite{bai1998} (in the sample covariance context). 
\end{remark}

\NTheorem{Anisotropic law}{
Assume \emph{(A)-(D)} and fix arbitrary $ \gamma > 0 $. 
Then uniformly for all $z= \tau +\cI\1\eta \in \Cp$ with $\eta \geq N^{\gamma-1}$, and for any two deterministic $ \ell^2$-unit vectors $ \vect{w},\vect{v} $ we have  
\bels{isoeq}{
\absbb{
\sum_{i,j=1}^N \overline{w}_i\1 G_{i j}(z)v_j- \sum_{i=1}^N m_i(z)\2\overline{w}_i\2v_i\,
}
\;\prec\;
\sqrt{\frac{\rho(z)}{N\1\eta}}+\frac{1}{N\1\eta} + \min\setbb{\!\frac{1}{\sqrt{\msp{-2}N\eta}}\2, \frac{\kappa(z)}{N\eta}}
\,,
}
where $\kappa$ is the function from Theorem \ref{thr:Local law}. 
}

\begin{corollary}[Delocalization of eigenvectors]
\label{crl:Delocalization of eigenvectors}
Assume \emph{(A)-(D)} and fix arbitrary $ \gamma > 0 $. 
Let $ \brm{u}^{(i)} \in \C^N $ be the $ \ell^2$-normalized eigenvector of $\brm{H}$ corresponding to the eigenvalue $\lambda_i$. 
All eigenvectors are delocalized in the sense that for any deterministic unit vector $\brm{b}\in \C^N$ we have
\bels{delocalization of eigenvectors}{ 
\absb{\2\brm{b}\cdot \brm{u}^{(i)}} \,\prec\, 
\frac{1}{\!\sqrt{N\2}}
\,.
}
In particular, the eigenvectors are completely delocalized, i.e., $\norm{\1\vect{u}^{(i)}}_\infty = \max_j \abs{\1u^{(i)}_j}\prec N^{-1/2}$.
\end{corollary}

\begin{definition}[$q$-full random matrix]
\label{def:q-full}
We say that $\vect{H}$ is {\bf $q$-full } for some $q>0$ (independent of $ N$) if either of the following applies: 
\begin{itemize}
\item $\vect{H}$ is real symmetric and $\EE \2h_{ij}^2 \geq q/N$ for all $i,j=1,\dots,N$;
\item 
$\vect{H}$ is complex hermitian and for all $i,j=1,\dots,N$ the real symmetric $2 \times 2$-matrix,
\[
\brm{\sigma}_{i j}\,:=\, 
\left(
\begin{array}{cc}
\EE \2(\Re\2 h_{i j})^2&\EE \2(\Re\2 h_{i j})(\Im\2 h_{i j})
\\
\EE \2(\Re\2 h_{i j})(\Im\2 h_{i j}) &\EE \2(\Im\2 h_{i j})^2
\end{array}
\right)\,,
\]
is strictly positive definite such that $ \brm{\sigma}_{i j}\geq q/N$.
\end{itemize}
\end{definition}
If $ \brm{H} $ is real symmetric, then the $ q$-fullness of $ \brm{H} $ is equivalent to the property (B) with $ L = 1 $ and $ q = p $. On the other hand, in the complex hermitian case the $ q $-fullness condition is stronger than a lower bound on $ \EE\,\abs{h_{ij}}^2 = s_{ij}$,  and it can not be captured by the matrix $ \brm{S} $ alone. 

\begin{theorem}[Universality]
\label{thr:Universality}
Suppose \emph{(A)} and \emph{(D)} hold, and $ \brm{H} $ is $ q $-full. 
Then for all $\eps>0$, $n \in \N$ and all smooth compactly supported observables $ F:\R^n \to \R$, there are two positive constants $ C $ and $ c $, depending on $ \eps,q $ and $ F $ in addition to the model parameters,  such that for any $\tau \in \R$ with $\rho(\tau)\geq \eps $ the local eigenvalue distribution is universal,
\[
\absbb{\,\EE\2F\Bigr( \bigl(\2N\msp{-1}\rho(\1\lambda_{i(\tau)}\msp{-1})\2(\1\lambda_{i(\tau)}-\lambda_{i(\tau)+j})\2\bigr)_{j=1}^n \Bigr)
-\,
\EE_{\rm G}\2F\Bigl( \bigl(\2N\msp{-1}\rho_{\rm sc}(0)\2(\1\lambda_{\ceil{N/2}}-\lambda_{\ceil{N/2}+j})\2\bigr)_{j=1}^n \Bigr)
}
\,\leq\,
C\2N^{-c}.
\]
Here, $\EE_{\rm G}$ denotes the expectation with respect to the standard Gaussian ensemble, i.e., with respect to GUE and GOE in the cases of complex Hermitian and real symmetric $\vect{H}$, respectively, and $\rho_{\rm sc}(0)=1/(2\1\pi)$ is the value of Wigner's semicircle law at the origin. 
\end{theorem}

This theorem concerns the universality in the bulk. With the help of our local law one may also prove a weaker version of the universality at the edges (including the internal edges). Since our local law, Theorem~\ref{thr:Local law}, is optimal at the  edges, a direct application of the Green function comparison theorem from Section 6 of \cite{EYYrigi} (with straightforward adjustments)
shows edge universality in the sense that the edge statistics may depend only on the second moments
encoded in the matrices  $ \brm{\sigma}_{ij} $. In particular, it is the same as the edge statistics of a Wigner-type matrix  with centered Gaussian entries with coinciding second moments. This argument holds for the extreme edges as well as for the  internal edges. However, it  does not yet prove the \emph{Tracy-Widom law}, i.e. that  the edge statistics is independent even of the variances $\brm{S}$.
 
\begin{convention}[Constants and comparison relation]
\label{conv:Constants and comparison relation}
We use the convention that every positive constant with a lower star index, such as $\delta_*$, $c_*$ and $\lambda_*$, explicitly depends only on the model parameters $P$, $p$ and $L$ from \emph{(B)--(D)}. These dependencies can be reconstructed from the proofs, but
we will not follow them. Constants $c, c_1, c_2, \dots, C, C_1, C_2, \dots$ also depend only on $P$, $p$ and $L$. They will have a local meaning within a specific proof. 

For two non-negative functions $\varphi$ and $\psi$ depending on a set of parameters $u \in U$, we use the {\bf comparison relation} 
\bels{}{
\varphi \;\gtrsim\; \psi\,,
}
if there exists a positive constant $c$, depending explicitly on $P$, $p$ and $L$ such that $\varphi(u) \geq c \2\psi(u)$ for all $u \in U$.  
The notation 
 $\psi \sim \varphi$ means that both $\psi \lesssim \varphi$ and $\psi \gtrsim \varphi$ hold true. In this case we say that $\psi$ and $\varphi$ are {\bf comparable}. We also write $\psi= \varphi +\Ord(\vartheta)$, if $|\psi-\varphi|\lesssim \vartheta$.
\end{convention}

We denote the normalized scalar product between two vectors $\vect{u},\vect{w} \in \C^N$ and the average of a vector by 
\bels{defs of inner product and averaging}{
\avg{\1\vect{u},\vect{w}} \,:=\, \frac{1}{N}\sum_{i=1}^N \overline{u}_i\2 w_i
\,,\qquad\text{and}\qquad
\avg{\vect{w}} \,:=\, \frac{1}{N}\sum_{i=1}^N w_i
\,,
}
respectively. 
Note that with this convention $\abs{\avg{\1\brm{u},\brm{u}\1}} = N^{-1}\| \brm{u}\|^2_{\ell^2} $.

\section{Bound on the random perturbation of the QVE}

We will make the following standing assumptions for the rest of this paper,
\begin{itemize}
\item \emph{The assumptions \emph{(A)--(D)} hold true and an arbitrary tolerance exponent $ \gamma \in (0,1) $ is fixed;}
\end{itemize}
which are always assumed to hold unless explicitly otherwise stated. 

We introduce the notation $\vect{G}^{(V)}$ for the resolvent of the matrix 
$\vect{H}^{(V)}$, which is identical to $\vect{H}$ except for the removal of the rows and columns corresponding to the indices $V \subseteq \{1,\dots,N \}$. The enumeration of the indices is kept, even though $\vect{G}^{(V)}$ has a lower dimension. 

The diagonal elements of the resolvent, $\vect{g}:=(G_{11},\dots,G_{NN})$, satisfy the perturbed quadratic vector equation
\bels{perturbed SCE}{
-\,\frac{1}{g_i(z)\!} \,=\, z + \sum_{j=1}^N s_{ij}\1 g_j(z) + d_i(z)\,,
}
for all $z \in \Cp$ and $i=1,\dots,N$.
The random perturbation $\vect{d}=(d_1,\dots,d_N)$ is given by
\bels{perturbation d}{
d_k 
\;:=\;
\sum^{(k)}_{i\neq j} h_{ki}G^{(k)}_{ij}h_{jk} 
\,+\,
\sum^{(k)}_i (\1\abs{h_{ki}}^2-s_{ki})\2 G^{(k)}_{ii}
\,-\,
\sum^{(k)}_i s_{ki} \frac{G_{ik}G_{ki}}{g_k}
\,-\,
h_{kk}
\,-\,
s_{kk}\1g_k	
\,.
}
Here and in the following, the upper indices on the sums indicate which indices are not summed 
over.
For the proof of this simple identity as well as \eqref{off-diagonal expansion} below via the Schur complement formula we refer to \cite{EKYY}. As in \eqref{perturbation d} we will often omit the dependence on the spectral parameter $z$ in our notation, i.e., $G_{i j}=G_{i j}(z)$, $d_k=d_k(z)$, etc.. 

We will now derive an upper bound on $\norm{\vect{d}}_\infty = \max_i|d_i|$, provided $|g_i-m_i|$ is bounded by a small constant. At the same time we will control the off-diagonal elements $G_{k l}$ of the resolvent.
These satisfy the identity
\bels{off-diagonal expansion}{   
G_{kl} \,=\, G_{kk} G^{(k)}_{ll} \sum_{i,j}^{(k l)} h_{ki}\1G^{(kl)}_{ij}h_{jl} \,-\,G_{kk}G^{(k)}_{ll}h_{kl} 
\,,
} 
for  $k\neq l$. The strategy in what follows below is that \eqref{perturbation d} and 
\eqref{off-diagonal expansion} are used to improve a rough bound on the entries of the resolvent $\vect{G}$ to 
get the correct bounds on the random perturbation and the off-diagonal resolvent elements. Later, in Section \ref{result in the bulk section}, the stability of the QVE under the small perturbation, $\vect{d}$, will provide the improved bound on the diagonal elements, $G_{ii}-m_i=g_i-m_i$. 

We introduce a short notation for the difference between $\vect{g}$ and the solution $\vect{m}$ of the unperturbed equation 
\eqref{SCE},
\bels{}{
\Err_{\mathrm{d}}(z) \,&:=\, \max_i|G_{ii}(z)-m_i(z)|\,,
\\
\Err_{\mathrm{o}}(z) \,&:=\, \max_{i\neq j } \abs{G_{ij}(z)}\,,
\\
\Err(z)\,&:=\, \max \big\{\Err_{\mathrm{d}}(z) ,\Err_{\mathrm{o}}(z) \big\}
\,.
}
The following lemma is analogous to Lemma 5.2 in \cite{EKYY} with minor modifications. For the completeness of this paper, we repeat these arguments. One small modification is that our estimates also deal with the regime where $|z|$ is large. 
To keep the formulas short we denote
\[
[z]\,:=\, 1+|z|
\,.
\]
The dependence of the upcoming error bounds on $[z]$ is not always optimal and this dependence is not kept in the statement of our main result Theorem \ref{thr:Local law}, either. In fact, the regime $[z]\sim 1$ is the most interesting, since our results show that the spectrum of $\vect{H}$ lies a.w.o.p. inside a compact interval (cf. Corollary \ref{crl:Convergence of cumulative eigenvalue distribution}). For the first reading we therefore recommend to think of $[z]=1$ in most of our proofs. The $[z]$-dependence is used mainly in order to propagate a bound from the regime of very large imaginary part of the spectral parameter ($\Im \1z \geq N^{\15}$) to the entire
domain, on which Theorem \ref{thr:Local law} holds.

\NLemma{Bound on perturbation}{
There is a small positive constant \cob $\lambda_* \sim 1$\nc, such that 
uniformly for all spectral parameters $z=\tau+\cI\1\eta \in \Cp$ with $\eta \geq N^{\gamma-1}$:
\begin{subequations}
\label{bulk result}
\begin{align}
\label{}
|d_k(z)| \,\Ind \big( \Err(z) \leq \lambda_*/[z] \2\big)
\;&\prec\;
[z]^{-2} \sqrt{\frac{\Im \la \1\vect{g}(z)\ra}{N \eta}}+\frac{1}{\sqrt{N}}
\,,
\\
\label{bound on Err_o}
\Err_{\mathrm{o}}(z) \,\Ind  \big( \Err(z) \leq \lambda_*/[z]\2 \big)
\;&\prec\;[z]^{-2}
\Biggl(
\sqrt{\frac{\Im \la \1\vect{g}(z)\ra}{N \eta}}+\frac{1}{\sqrt{N}}
\Biggr)
\,.
\end{align}
\end{subequations}
}

For the proof of this lemma we will need an additional property of the solution of the QVE that is a corollary of Theorem \ref{thr:Solution of the QVE}, where all properties of $\vect{m}$ taken from \cite{AEK1} are summarized.

\NCorollary{Bounds on solution}{The absolute value of the solution of the QVE satisfies
\bels{}{
\abs{\1m_i(z)}
\,\sim\, 
[z]^{-1}
\,,
\qquad  z\in \Cp\,,\; i=1,\dots,N
\,.
}
}

\begin{Proof}[Proof of Lemma \ref{lmm:Bound on perturbation}]
Here we use the three large deviation estimates, 
\begin{subequations}
\label{large deviation bounds}
\begin{align}
\label{large deviation off diagonal}
\absbb{\,\sum^{(k)}_{i\neq j} h_{ki}\1G^{(k)}_{ij}h_{jk}}
\;&\prec\;
\biggl(\, \sum^{(k)}_{i\neq j} s_{ki}\1 s_{jk} \big|G^{(k)}_{ij}\big|^2 \biggr)^{\!1/2}
,
\\
\label{large deviation off diagonal 2}
\absbb{\,\sum^{(k l)}_{i,j} h_{ki}\1G^{(k l)}_{ij}h_{jl}}
\;&\prec\;
\biggl(\, \sum^{(k l)}_{i,j} s_{ki}\1 s_{jl} \big|G^{(k l)}_{ij}\big|^2 \biggr)^{\!1/2}
,
\\
\label{large deviation diagonal}
\absbb{\,\sum^{(k)}_i \bigl(\2\abs{h_{ki}^2}-s_{ki}\bigr)\1 G^{(k)}_{ii}}
\,&\prec\, 
\biggl(\,  \sum^{(k)}_i s_{ki}^2\2 \big|G^{(k)}_{ii}\big|^2 \biggr)^{\!1/2}
.
\end{align}
\end{subequations}
Since $\vect{G}^{(V)}$ is independent of the rows and columns of $\vect{H}$ with indices in $V$, these estimates follow directly from the large deviation bounds in Appendix C of \cite{EKYY}.
Furthermore, we use 
\bels{H entry bound}{
\abs{\1h_{ij}}\,\prec\, N^{-1/2}\,,\qquad s_{i j}\,\leq\, N^{-1}
\,,
}
where latter the inequality is just assumption \eqref{flat} and the bound on $h_{i j}$ follows from \eqref{bounded moments}. 
We remark that the stochastic domination in \eqref{large deviation bounds} and \eqref{H entry bound} is uniform in $k,l$ and $ i,j $, respectively, i.e., the threshold function $ N_0 $ in Definition~\ref{def:Stochastic domination} does not depend on $ i,j,k,l $.

We will now show that the removal of a few rows and columns in $\vect{H}$ will only have a small effect on the entries of the resolvent.
The general resolvent identity,
\bels{resolvent identity}{
G_{i j} \,=\, G_{i j}^{(k)} + \frac{G_{i k}G_{k j}}{G_{kk}}\,,\qquad k \not \in \{i,j\}\,,
}
leads to the bound
\bels{removing upper index}{
\big|G_{i j}^{(k)} - G_{i j}\big|\,\Ind\big( \Err \leq \lambda_*/[z] \big)
\,=\,
\frac{| G_{i k} G_{k j} |}{|g_k|} \,\Ind\big( \Err \leq \lambda_*/[z] \big)
\,\lesssim\,
[z]\Err_{\mathrm{o}}^2\,.
}
In the inequality we used that $|m_k(z)|\sim [z]^{-1}$ (cf. Corollary \ref{crl:Bounds on solution}), $|g_k| = |m_k|+\Ord(\Err)$ and that $\lambda_*$ is chosen to be small enough.
We use \eqref{removing upper index} in a similar calculation for $G_{i j}^{(l)}$
and find that on the event where $\Err \leq \lambda_*/[z]$,
\bels{}{
\big|G_{i j}^{(k l)} - G_{i j}^{(l)}\big|
\,=\,
\frac{\big| G_{i k}^{(l)} G_{k j}^{(l)} \big|}{\big|G_{k k}^{(l)}\big|} 
\,\lesssim\,
\frac{
(\,| G_{i k}| + \Ord(\1[z]\1\Err_{\mathrm{o}}^2\2)\,)
(\,|G_{k j} | + \Ord(\1[z]\1\Err_{\mathrm{o}}^2\2)\,)
}{\abs{\1g_k} + \Ord(\1[z]\1\Err_{\mathrm{o}}^2\1)
} 
\,.
}
Again using \eqref{removing upper index} and that the denominator of the last expression is comparable to $[z]^{-1}$, we conclude
\bels{G^T-G is small}{
|G_{ij}^{(kl)} - G_{i j}|\,\Ind\big( \Err \leq\lambda_*/[z] \big)
\,\lesssim\, [z]\2\Err_{\mathrm{o}}^2\,,
}
provided $\lambda_*$ is small. 
Therefore, we see that it is possible to remove one or two upper indices from $G_{i j}$ for the price of a term, whose size is at most $[z]\Err_{\mathrm{o}}^2$.

We have now collected all necessary ingredients and use them to estimate all the terms in \eqref{perturbation d} one by one. We start with the first summand.
By \eqref{large deviation off diagonal} we find
\bels{first summand large deviation}{
\absbb{\,\sum^{(k)}_{i\neq j} h_{ki}G^{(k)}_{ij}h_{jk}}^2
\,\prec\, 
\sum^{(k)}_{i\neq j} s_{ki} s_{jk} \big|G^{(k)}_{ij}\big|^2 
\,\leq\,
\frac{1}{N^2}\sum^{(k)}_{i\neq j} \big|G^{(k)}_{ij}\big|^2 
\,.
}
With the help of \eqref{removing upper index} we remove the upper index from $G_{ij}^{(k)}$ and get
\bels{first summand bound}{
\absbb{\,\sum^{(k)}_{i\neq j} h_{ki}G^{(k)}_{ij}h_{jk}}^2
\Ind\big( \Err \leq \lambda_*/[z] \big)
\,\prec\,
\bigl(\1\Err_{\mathrm{o}}^2 \,+\,[z]^2\Err_{\mathrm{o}}^4\2\bigr)
\Ind\bigl( \Err \leq \lambda_*/[z]\1 \bigr)
\,\lesssim\,\Err_{\mathrm{o}}^2 
\,.
}

For the second summand in \eqref{perturbation d} we use the large deviation bound for the diagonal, \eqref{large deviation diagonal}, and find that
\bels{second summand large deviation}{
\absbb{\,\sum^{(k)}_i (\abs{h_{ki}}^2-s_{ki})\2 G^{(k)}_{ii}}^2
\,\prec\; 
 \sum^{(k)}_i s_{ki}^2\2 \big|G^{(k)}_{ii}\big|^2 
\,\leq\,
\frac{1}{N^2}\sum^{(k)}_i \big|G^{(k)}_{ii}\big|^2\,.
}
By removing the upper index again we estimate
\bels{}{
\big|G^{(k)}_{ii}\big|\,\Ind\big( \Err \leq \lambda_*/[z] \big)
\,\lesssim\, |m_i|+\Err_{\mathrm{d}}+ [z]\Err_{\mathrm{o}}^2\,.
}
We use this in \eqref{second summand large deviation} and for sufficiently small $\lambda_*$ we arrive at
\bels{second summand bound}{
\absbb{\,\sum^{(k)}_i (\abs{h_{ki}}^2-s_{ki})\2 G^{(k)}_{ii}}^2
\Ind\big( \Err \leq \lambda_*/[z] \big)
\;\prec\;
\frac{1}{[z]^2 N}\,.
}

The third summand in \eqref{perturbation d} is estimated directly by
\bels{third summand bound}{
\absbb{\,\sum^{(k)}_i s_{ki} \frac{G_{ik}G_{ki}}{g_k}}\,\Ind\big( \Err \leq \lambda_*/[z] \big)
\,\leq\,
\frac{\Err_{\mathrm{o}}^2}{|g_k|}
\,\Ind\big( \Err \leq \lambda_*/[z] \big)
\,\lesssim\, 
\Err_{\mathrm{o}}
\,.
}

We combine the estimates for the individual terms \eqref{first summand bound}, \eqref{second summand bound}, \eqref{third summand bound} and \eqref{H entry bound}.
Altogether we conclude that 
\bels{bound on d}{ 
|d_k| \,\Ind\big( \Err \leq \lambda_*/[z] \big)
\,\prec\,
\Err_{\mathrm{o}}(z)+\frac{1}{\sqrt{N}}
 \,.
}

We will now derive in a similar fashion a stochastic domination bound for the off-diagonal error term $\Err_{\mathrm{o}}$. Afterwards, we will combine the two bounds and infer the claim of the lemma. 
For the off-diagonal error term we proceed along the same lines as for $|d_k|$, using \eqref{off-diagonal expansion} instead of 
\eqref{perturbation d}. 
For $k \neq l$ we find  
\bels{bound on off diagonal error start}{
|G_{kl}|^2
\,\prec\, 
|g_k|^2\absb{G_{ll}^{(k)}}^2
\Biggl(
\frac{1}{N^2}
\sum^{(k l)}_{i, j}  \big|G^{(k l)}_{ij}\big|^2+\frac{1}{N}
\Biggr)
\,.
}
Here, we applied the large deviation bound \eqref{large deviation off diagonal 2}.
Using the Ward identity for the resolvent $\vect{G}^{(kl)}$,
\bels{Ward identity}{
\sum^{(kl)}_{j} \big|G^{(kl)}_{ij}\big|^2 \,=\, \frac{\Im \,G_{i i}^{(kl)}}{\eta}\,,
} 
 and \eqref{removing upper index} for removing the upper index of $G_{ll}^{(k)}$ we get
\bels{}{
|G_{kl}|^2\,\Ind\big( \Err \leq \lambda_*/[z] \big)
\,\prec\,
[z]^{-4}
\Biggl(\frac{1}{N^2\eta}\sum_{i}^{(k l)} \Im \2G_{ii}^{(k l)}+\frac{1}{N} \Biggr)
\,.
}
We remove the upper indices from $G_{ii}^{(k l)}$ and end up with
\bels{bound on off diagonal error}{
\Err_{\mathrm{o}} \,\Ind\big( \Err \leq \lambda_*/[z] \big)
\,\prec\, 
[z]^{-2}\Biggl( 
\sqrt{\frac{\Im \,\avg{\1\vect{g}\1}}{N \eta}}+ \sqrt{\frac{[z]}{N\eta}}\;\Err_{\mathrm{o}}+\frac{1}{\sqrt{N}}
\Biggr)
\,.
}
The bound remains true without the summand containing $\Lambda_{\mathrm{o}}$ on the right hand side, since this term can be absorbed into the left hand side, as its coefficient is bounded by $N^{-\gamma/2}$, while on the left $\Lambda_{\mathrm{o}}$ is not multiplied by a small coefficient. 
Putting \eqref{bound on d} and \eqref{bound on off diagonal error} together yields the desired result \eqref{bulk result}.
\end{Proof}

\section{Local law away from local minima}
\label{result in the bulk section}

In this section we will use the stability of the QVE to establish the main result away from the local minima of the density of states inside its own support, i.e. away from the set
\bels{definition of M}{
\mathbb{M}\,:=\,\setb{\tau \in \supp \rho: \;\tau \text{ is the location of a local minimum of } \rho\, }
\,.
}
The case where $z$ is close to $\mathbb{M}$ requires a more detailed analysis. 
This is given is Section \ref{Main result at the edges}. 
At the end of this section we will also sketch the proof of Corollary \ref{crl:Local law in the bulk}. 
In this section we prove the following.

\NProposition{Local law away from local minima}{ 
Let $\delta_*$ be any positive constant, depending only on the model parameters $p$, $P$ and $L$. 
Then, uniformly for all $ z=\tau+\cI\1\eta$ with $\eta \geq N^{\gamma-1}$ and $\dist(z,\mathbb{M})\geq \delta_*$, we have 
\begin{subequations}
\label{bulk error bound}
\begin{align}
\label{bulk error bound:diag}
[z]^2\Err_{\mathrm{d}}(z)+\norm{\vect{d}(z)}_\infty
\,&\prec\,
[z]^{-2}\sqrt{\frac{\rho(z)}{N \eta}}+\frac{\;[z]^{-6}\msp{-5}}{N\eta}\,+\frac{1}{\!\sqrt{N}}\;,
\\
\label{bulk error bound:off-diag}
\Err_{\mathrm{o}}(z)
\,&\prec\,
[z]^{-2}\sqrt{\frac{\rho(z)}{N\eta}}
+ \frac{\;[z]^{-4}\msp{-5}}{N\eta}\,+\frac{\;[z]^{-2}\msp{-5}}{\!\sqrt{N}}\,.
\end{align}
\end{subequations}
Furthermore, on the same domain, for any  sequence of deterministic vectors $\vect{w}=\vect{w}^{(N)} \in \C^N$  with the uniform bound, $\norm{\vect{w}}_\infty\leq 1$, we have
\bels{bulk average error bound}{
\abs{\avg{\2\vect{w}\1,\vect{g}(z) -\vect{m}(z)\1}}\;\prec\; [z]^{-3}\2\frac{\rho(z) }{N\eta}+\frac{\,[z]^{-7}\msp{-3}}{(N\eta)^2\,}+\frac{\;[z]^{-2}\msp{-5}}{N}
\;.
} 
}
This proposition, combined with the properties of $ \rho $ given in Theorem~\ref{thr:Solution of the QVE} later, yields the local law (Theorem \ref{thr:Local law}) away from the set $\mathbb{M}$.
Indeed, using $ \rho(z) \gtrsim [z]^{-2}\eta $ (cf. relations \eqref{qualitative size of rho} below)  and $ \kappa(z) \ge 0 $ we see that \eqref{bulk error bound} implies  \eqref{pointwise local law}.

In order to see that also the averaged local law \eqref{averaged local law} follows from \eqref{bulk average error bound} we split the domain $ \sett{z\in\Cp:\dist(z,\mathbb{M})} \ge \delta_\ast $ into three subdomains that are considered separately.
To this end, let $ B_0 $ and $ B_1 $ be the upper bounds on $ \kappa $ from \eqref{kappa bound} and \eqref{improved kappa bound}, respectively.

First we consider the regime $ \eta \ge \delta_\ast/2 $. Using $ \Delta^{1/3}+\rho \lesssim 1 $ we see that $ B_0 \gtrsim 1 $. 
Similarly, we get  $ B_1 \gtrsim \eta\2[z]^{-1} $. Since $ (N\eta)^{-1}B_k $, $ k=0,1$, are both bigger than the right hand side of \eqref{bulk average error bound}, we obtain \eqref{averaged local law} for $ \eta \ge \delta_\ast/2 $.

Now we consider the regime $ \eta \leq \delta_\ast/2  $, which is split into two cases depending on whether $ \dist(\Re\,z,\supp \rho) = 0 $, or not. 
In the former case $[z] \lesssim 1 $ and  $ \dist(z,\supp \rho) = \eta $, and \eqref{lower bound for rho when bounded away from MM} implies $ \rho(\Re\,z) \sim 1 $. Feeding these estimates into \eqref{kappa bound} and \eqref{improved kappa bound} yields $ B_0 \sim 1 $ and $ B_1 \gtrsim 1 $. These imply  \eqref{averaged local law}.

Finally, suppose $ \dist(\Re\,z,\supp \rho) \ge \delta_\ast/2  $ and $ \eta \leq \delta_\ast/2 $. 
In this regime $\Delta \sim 1 $ (cf. \eqref{def of Delta_delta}), while \eqref{rho away from supp} implies $\rho \sim \eta\,[z]^{-2} $.
Hence, $ B_0 \sim 1 $ and $ B_1 \gtrsim \eta\,[z]^{-1} $, and  $ (N\eta)^{-1}\min\sett{B_0,B_1} \ge [z]^{-1}N^{-1} $.
By comparing with the right hand side of \eqref{bulk average error bound} we conclude that \eqref{averaged local law} applies for all $ \dist(z,\mathbb{M}) \ge \delta_\ast $.

The proof of Proposition~\ref{prp:Local law away from local minima} uses a continuity argument in $z$. In particular, continuity of the solution of the QVE is needed. The statement of the following corollary is part of the properties of $\vect{m}$ listed in Theorem \ref{thr:Solution of the QVE} below. 

\NCorollary{Stieltjes-transform representation}{
For every $i=1, \dots, N$ there is a probability density $p_i:\R \to [0,\infty)$ with support in $[-2,2]$ such that $m_i$ is the Stieltjes-transform of this density, i.e.,
\bels{Stieltes-transform representation}{
m_i(z)\,=\, \int_\R\frac{p_i(\tau)\2\dif \tau}{\tau-z}\,,\qquad z \in \Cp\,.
}
The solution of the QVE is uniformly H\"older-continuous,
\bels{Hoelder continuity of m}{
\norm{\vect{m}(z_1)-\vect{m}(z_2)}_\infty\,\lesssim\,|z_1-z_2|^{1/3}\,,\qquad 
z_1,z_2 \in \overline \Cp\,. 
}
}
Since the solution can be extended to the real line, it is the harmonic extension to the complex upper half plane of its own restriction to the real line. Therefore, $\Im\2m_i(\tau)=\pi\2p_i(\tau)$ for $\tau \in \R$. The density of states is the average of the probability densities $p_i$, i.e., $\rho=\la \vect{p} \ra$. 

Since we will estimate the difference, $\vect{g}-\vect{m}$, we start by deriving an equation for this quantity. Using the QVE for $\vect{m}$ and the perturbed equation \eqref{perturbed SCE} for $\vect{g}$ we find
\bels{}{
g_i-m_i
&\,=\, 
-\,\frac{1}{z+(\vect{S}\vect{g})_i+d_i}+\frac{1}{z+(\vect{S}\vect{m})_i}
\\
&\,=\, 
\frac{(\vect{S}(\vect{g}-\vect{m}))_i+d_i}{(z+(\vect{S}\vect{g})_i+d_i)(z+(\vect{S}\vect{m})_i)}
\\
&\,=\,
m_i^2(\vect{S}(\vect{g}-\vect{m}))_i+m_i\1(g_i-m_i)(\vect{S}(\vect{g}-\vect{m}))_i+m_i\, g_i\, d_i\,.
}
Rearranging the terms leads to 
\bels{quadratic SCE}{
\big( (\Id-\diag(\vect{m})^2 \vect{S})(\vect{g}-\vect{m})  \big)_i
\,=\,
m_i\1(g_i-m_i)\1(\1\vect{S}(\vect{g}-\vect{m}))_i+m_i^2\, d_i+m_i\,( g_i - m_i)\, d_i\,.
}
In the proof of Proposition \ref{prp:Local law away from local minima}
we will view \eqref{quadratic SCE} as a quadratic equation for $\vect{g}-\vect{m}$ and we use its stability 
to bound $\Err_{\mathrm{d}}$ in terms of $\norm{\vect{d}}_\infty$.
We will now demonstrate this effect in the case when $z$ is far away from the support of the density of states. 

\begin{lemma}[Stability far away from support]
\label{lmm:Stability far away from support}
For $z \in \Cp$ with $|z|\geq 10$, we have
\bels{QVE stability away from support}{
\Err_{\mathrm{d}}(z)\2 \Ind(\Err_{\mathrm{d}}(z) \leq 4 \2|z|^{-1} )
\;\lesssim\;
\abs{z}^{-2}\norm{\vect{d}(z)}_\infty
\,.
}
Furthermore, there is a matrix valued function $\vect{T}:\Cp \to \C^{N \times N}$, depending only on $\brm{S} $ and satisfying the uniform bound
$\norm{\vect{T}(z)}_{\infty \to \infty}\lesssim 1$, such that for all $\vect{w} \in \C^N$ and $|z|\geq 10$ the averaged difference between $\vect{g}$ and $\vect{m}$ satisfies the improved bound
\bels{averaged QVE stability away from support}{
\absb{\avgb{\1\vect{w},\vect{g}(z)-\vect{m} (z)}}
\1\Ind\big( \Err_{\mathrm{d}}(z) \leq 4 \2|z|^{-1} \big)
\;\lesssim\;
\abs{z}^{-2} 
\big(\2\norm{\vect{w}}_\infty \norm{\vect{d}(z)}_{\infty}^2+ \abs{\avg{\1\vect{T}(z) \vect{w},\vect{d}(z)}}\2\big)
\,.
}
\end{lemma}

For a matrix $ \brm{A} $ we denote by $ \norm{\brm{A}}_{\infty\to\infty} $ the operator norm of $ \brm{w} \mapsto \brm{A}\brm{w} $ on $ \ell^\infty $. 

\begin{Proof}
Since the matrix $\vect{S}$ is flat (cf. \eqref{flat}), it satisfies the norm bound $\norm{\vect{S}}_{\infty\to\infty}\leq 1$.
We also have the trivial bound $|m_i(z)|\leq 1/\dist(z,\supp \rho) \leq 2\1|z|^{-1}\leq 1/5$ 
at our disposal. This follows directly from the Stieltjes transform representation \eqref{Stieltes-transform representation}. In particular, 
\bels{}{
\norm{(\2\Id-\diag(\vect{m})^2 \vect{S}\1)^{-1}}_{\infty \to \infty}
\,\leq\, 2
\,,
}
from the geometric series. By inverting the matrix $\Id-\diag(\vect{m})^2 \vect{S}$ and using the trivial bound on $\vect{m}$ in \eqref{quadratic SCE} we find
\bels{}{
\Err_{\mathrm{d}}(z) \,\leq\, 
4\2\Bigl(\,
|z|^{-1}\Err_{\mathrm{d}}(z)^2+|z|^{-1}\Err_{\mathrm{d}}(z)\1\norm{\vect{d}(z)}_\infty+2\2|z|^{-2}\norm{\vect{d}(z)}_\infty  
\Bigr)
\,.
}
Using the bound inside the indicator function from \eqref{QVE stability away from support} and $|z|\geq 10$ the assertion \eqref{QVE stability away from support} of the lemma follows. 

The bound for the average, \eqref{averaged QVE stability away from support}, follows by taking the inverse of $\Id-\diag(\vect{m})^2 \vect{S}$ on both sides of \eqref{quadratic SCE} and using \eqref{QVE stability away from support} and $|m_i|\sim |z|^{-1}$. 
\end{Proof}

For the proof of Proposition \ref{prp:Local law away from local minima}
we use the stability of \eqref{quadratic SCE} also close to $\supp\rho$. This requires more care and is carried out in detail in \cite{AEK1}. The result of that analysis is Theorem~\ref{thr:Stability around support}. Here we will only need the following consequence of that theorem and \eqref{lower bound for rho when bounded away from MM}.

\begin{corollary}[Stability away from minima]
\label{crl:Stability away from minima}
Suppose $\delta_*$ is an arbitrary positive constant, depending only on the model parameters $p$, $P$ and $L$.
Let $\vect{d}: \Cp \to \C^N$, $\vect{g}: \Cp \to (\C\backslash\sett{0})^N$ be arbitrary vector valued functions on the complex upper half plane that satisfy 
\bels{}{
-\frac{1}{g_i(z)} \,=\, z + \sum_{j=1}^N s_{i j} \2g_j(z) +d_i(z)\,, \quad z \in \Cp\,.
}
There exists a positive constant $\lambda_* \sim 1$, such that the QVE is stable away from $\mathbb{M}$,
\bels{rough stability}{
\norm{\1\vect{g}(z)-\vect{m}(z)}_\infty\,\Ind\big(\2\norm{\vect{g}(z)-\vect{m}(z)}_\infty \leq \lambda_* \big)
\,\lesssim\, \,\norm{\vect{d}(z)}_\infty\,,\quad z \in \Cp\,,\; 
\dist(z,\mathbb{M})\geq\delta_*
\,.
} 
Furthermore, there is a matrix valued function $\vect{T}:\Cp \to \C^{N \times N}$, depending only on  $\brm{S}$ and satisfying the uniform bound
$\norm{\vect{T}(z)}_{\infty \to \infty}\lesssim 1$, such that for all $\vect{w} \in \C^N$,
\bels{bound on average in the bulk}{
|\la \1\vect{w},\vect{g}(z)-\vect{m} (z)\ra |
\1\Ind\big(\2 \norm{\vect{g}(z)-\vect{m}(z)}_\infty \leq \lambda_* \big)
\,\lesssim\,\norm{\vect{w}}_\infty \norm{\1\vect{d}(z)}_{\infty}^2+ |\la\1 \vect{T}(z) \vect{w},\vect{d}(z)\1 \ra|\,,
}
for $z \in \Cp$ with $\dist(z,\mathbb{M})\geq\delta_*$.
\end{corollary}

Furthermore, the following fluctuation averaging result is needed. It was first established for generalized Wigner matrices with Bernoulli distributed entries in \cite{EYYber}. 

\NTheorem{Fluctuation Averaging}{
For any $z \in \Cp $, with $ \Im\,z \ge N^{\gamma-1} $,  and any sequence of deterministic vectors $\vect{w}=\vect{w}^{(N)} \in \C^N$ with the uniform bound, $\norm{\vect{w}}_\infty\leq 1$ the following holds true: If $\Err_{\mathrm{o}}(z) \prec \Phi/[z]^2$ \nc for some deterministic (N-dependent) $\Phi\leq N^{-\gamma/3}$ and $\Err(z)\prec N^{-\gamma/3}/(1+|z|)$ , then 
\bels{fluctuation averaging result}{
\abs{\avg{\1\vect{w},\vect{d}(z)}} \,\prec\, [z]^{-1}\2\Phi^2 + \frac{1}{N}
\,,
}
where $ \brm{d}(z) $ is defined in \eqref{perturbation d}.
}

\begin{Proof} The proof directly follows the one given in \cite{EKYY}. We only mention some minor necessary modifications. 
Let $Q_kX:=X-\EE[X|\vect{H}^{(k)}]$ be the complementary projection to the conditional expectation of a random variable $X$ given the matrix $\vect{H}^{(k)}$, in which the $k$-th row and column are removed. From the definition of $\vect{d}$ in \eqref{perturbation d} and Schur's complement formula in the form,
\bels{diagschur}{
\frac{1}{G_{kk}}\,=\, h_{kk}-z-\sum_{i,j}^{(k)}h_{k i }G_{i j}^{(k)}h_{j k}\,,
}
we infer the identity
\[
d_k\,=\, -\2Q_k\frac{1}{G_{kk}}-s_{kk}G_{kk}- \sum_{i}^{(k)}s_{k i}\frac{G_{i k}G_{k i}}{G_{kk}}\,.
\]
In particular, we have that a.w.o.p.
\[
\Big|d_k + Q_k \frac{1}{G_{kk}}\Big|\,\lesssim \;\frac{\;[z]^{-1}\msp{-7}}{N} \,+ [z]\Err_{\mathrm{o}}^2\,.
\]
Thus, proving \eqref{fluctuation averaging result} reduces to showing 
\[
\absbb{\2\frac{1}{N}\sum_{k=1}^N \overline{w}_k\2 Q_k \frac{1}{G_{kk}}\,}
\,\prec\; 
[z]^{-1}\2\Phi^2 + \frac{1}{N}\,.
\]
In the setting where $\vect{H}$ is a generalized Wigner matrix and $|z|\leq 10$ this bound is precisely the content of Theorem 4.7 from \cite{EKYY}. 

The a priori bound used in the proof of that theorem is replaced by
\bels{FA a priori bound}{
\absbb{\1Q_k \frac{1}{G^{(V)}_{kk}}}
\,\prec\, 
\Err_{\mathrm{o}} + \frac{1}{\!\sqrt{N}}
\,,
}
for any $V \subseteq \{1, \dots,N\}$ with $N$-independent size. This bound is proven in the same way as \eqref{bound on d}. Here, the $N_0$ hidden in the stochastic domination depends on the size $|V|$ of the index set. Following the proof of Theorem 4.7 given in \cite{EKYY} with \eqref{FA a priori bound} and tracking the $z$-dependence, 
\[
\frac{1}{|G_{kk}^{(V)}(z)|}\,\prec\, [z]\,,
\]
yields the fluctuation averaging, Theorem \ref{thr:Fluctuation Averaging}.
\end{Proof}

\begin{Proof}[Proof of Proposition \ref{prp:Local law away from local minima}]
Let us show first that \eqref{bulk average error bound} follows directly from \eqref{bulk error bound} by applying the fluctuation averaging, Theorem \ref{thr:Fluctuation Averaging}. 
Indeed, \eqref{bulk error bound} provides a deterministic bound on the off-diagonal error, $\Err_{\mathrm{o}}$, which is needed to apply the fluctuation averaging to the second terms on the right hand sides of \eqref{bound on average in the bulk} and \eqref{averaged QVE stability away from support}.
It also shows that the indicator functions on the left hand sides of \eqref{bound on average in the bulk} and \eqref{averaged QVE stability away from support} are a.w.o.p. nonzero. 
The stability bound \eqref{averaged QVE stability away from support} valid in the large $ \abs{z} $ regime is necessary to get the correct $ [z] $-factors in \eqref{bulk average error bound}. 
Thus, \eqref{bulk average error bound} is proven, provided \eqref{bulk error bound} is true. 

The proof of \eqref{bulk error bound} is split into the consideration of two different regimes.
In the first regime the absolute value of $z$ is large, $|z|\geq N^{\25}$. In this case we make use only of weak a priori bounds on the resolvent elements and the entries of $\vect{d}$. Together with Lemma \ref{lmm:Stability far away from support} they will suffice to prove \eqref{bulk error bound} in this case.
 In the second regime, $|z|\leq N^{\25}$, we use a continuity argument. We will establish a gap in the possible values that the continuous function, $z\mapsto [z]\Err(z)$, might have. Here, the stability result Corollary \ref{crl:Stability away from minima} is used. We use this gap to propagate the bound with the help of Lemma \ref{lmm:Bound propagation} in the appendix from $|z|=N^{\25}$ to the whole domain where $|z|\leq N^{\25}$, $\eta \geq N^{\gamma-1}$ and we stay away from $\mathbb{M}$.  
 
{\sc Regime 1: } Let $|z|\geq N^{\15}$. We show that the indicator functions in the statement of Lemma~\ref{lmm:Bound on perturbation} are a.w.o.p. not vanishing. We start by showing that the diagonal contribution, $\Err_{\mathrm{d}}$, to $\Err$ is sufficiently small. The reduced resolvent elements for an arbitrary $V \subseteq \{1, \dots,N\}$ satisfy 
\bels{reduced resolvent a priori bound}{
|G^{(V)}_{i j}(z)|\,\leq\,\eta^{-1}\,\leq\, N^{1-\gamma}\,.
}
From this and the definition of $\vect{d}$ in \eqref{perturbation d} we read off the  a priori bound,
\bels{d a priori bounds}{
\norm{\vect{d}(z)}_\infty\,\prec\,N^{\22-\gamma}.
}
Here, we used the general resolvent identity \eqref{resolvent identity} in the form
$G_{i k}G_{k i}=g_k(g_i-G_{ii}^{(k)})$.
Since $\vect{g}$ satisfies the perturbed QVE \eqref{perturbed SCE} and $|\sum_{j=1}^N s_{ij}\1 g_j(z) + d_i(z)|\prec N^{2-\gamma}$ from \eqref{d a priori bounds} and \eqref{reduced resolvent a priori bound} we conclude that uniformly for $|z|\geq N^{\12}$ we have 
\bels{diagonal bound}{
|g_k(z)|\,\leq\, 2\2|z|^{-1},\qquad \text{a.w.o.p.}\,.
}
With the trivial bound $|m_i(z)|\leq 1/\dist(z, \supp \rho)$ on the solution of the QVE we infer that on this domain the indicator function in \eqref{QVE stability away from support} is a.w.o.p. non-zero and therefore uniformly for $|z|\geq N^{\22}$. Lemma~\ref{lmm:Stability far away from support} yields 
\bels{diagonal error bound}{
\Err_{\mathrm{d}}(z)
\,\lesssim\, 
|z|^{-2}\2\norm{\vect{d}(z)}_\infty
\,\leq\, 
N^{-\gamma/2}\2|z|^{-1}
,
\qquad \text{a.w.o.p.}
\,.
}
In the last inequality we have used \eqref{d a priori bounds} in the form $ \norm{\brm{d}}_\infty \leq N^{-\gamma/2}N^2 $ a.w.o.p. (cf. Definitions~\ref{def:Stochastic domination} and~\ref{def:Overwhelming probability}) and the extra factor $ [z]^{-2} $ on the right hand side of \eqref{QVE stability away from support}. 
Thus, for $|z|\geq N^{\22}$ the diagonal contribution to $\Err$ does not play a role in the indicator function in the statement of Lemma~\ref{lmm:Bound on perturbation}. 

Now we derive a similar bound for the off-diagonal contribution $\Err_{\mathrm{o}}$.
Using the resolvent identity \eqref{resolvent identity}  for $i=j$ again, the bound $|h_{i j}|\prec N^{-1/2}$ on the entries of the random matrix and the a priori bound on the reduced resolvent elements, \eqref{reduced resolvent a priori bound}, in the expansion formula \eqref{off-diagonal expansion} yields
\bels{}{
|G_{k l}(z)|\,\prec\, \big(\,|g_k(z)\1g_l(z)|+|G_{k l}(z)G_{l k}(z)|\,\big)\2N^{\22-\gamma}
\,,\qquad 
|G_{k l}(z)|\,\prec\,|g_k(z)|\2N^{\23-\gamma}\,,
}
for $k \neq l$. With the bound \eqref{diagonal bound} we conclude that
\bels{}{
\Err_{\mathrm{o}}(z)
\;\prec\; |z|^{-2}\2N^{2-\gamma}+ |z|^{-1}\2N^{\25-2\gamma}\Err_{\mathrm{o}}(z)\,,\qquad |z|\2\geq\2 N^{\12}\,.
}
Thus, $\Err_{\mathrm{o}}\prec N^{-3}|z|^{-1}$ on the domain where $|z|\geq N^{\25}$. We conclude that Lemma \ref{lmm:Bound on perturbation} applies in this regime even without the indicator functions in the formulas \eqref{bulk result}. We use the bound from this lemma for the norm of $\vect{d}$ and the off-diagonal contribution, $\Err_{\mathrm{o}}$, to $\Err$, while we use the first inequality in \eqref{diagonal error bound} for the diagonal contribution, $\Err_{\mathrm{d}}$. In this way, we get
\bels{Regime 1 first error bounds}{
|z|^2\Err_{\mathrm{d}}+\norm{\vect{d}}_\infty\;
&\prec\; 
|z|^{-2}\sqrt{\frac{\rho}{N \eta}}\,+|z|^{-2}\sqrt{\frac{\Err_{\mathrm{d}}}{N \eta}}+\frac{1}{\!\sqrt{N}}
\,,
\\
|z|^2\Err_{\mathrm{o}}\,
&\prec\,\sqrt{\frac{\rho}{N \eta}}\,+\,\sqrt{\frac{\Err_{\mathrm{d}}}{N \eta}}+\frac{1}{\!\sqrt{N}}
\,,
}
where we also used $g_k = m_k + \Ord(\Err_{\mathrm{d}})$. 
Applying the weighted Cauchy-Schwarz inequality, $ \sqrt{\alpha\2\beta\1}\leq \theta\,\alpha+\theta^{-1}\beta $,  we find for any $\eps \in (0,\gamma)$ that the right hand side of the first inequality can be estimated further by
\[
|z|^2\Err_{\mathrm{d}}+\norm{\vect{d}}_\infty\,\prec\, 
|z|^{-2}\sqrt{\frac{\rho}{N \eta}}+N^{-\eps} |z|^2\Err_{\mathrm{d}}+ |z|^{-6}\frac{N^\eps}{N \eta}+\frac{1}{\!\sqrt{N}}\,.
\]
The term $N^{-\eps}|z|^2\Err_{\mathrm{d}}$ can be absorbed into the left hand side and by the definition of the stochastic domination and since $\eps$ is arbitrarily small the remaining $N^\eps$ on the right hand side can be replaced by $1$ without changing the correctness of this bound  (cf. (i) and (ii) of Lemma~\ref{lmm:Basic facts about stochastic domination}).
In this way we arrive at
\[
|z|^2\Err_{\mathrm{d}}+\norm{\vect{d}}_\infty\,\prec\, 
|z|^{-2}\sqrt{\frac{\rho}{N \eta}}+\frac{\; |z|^{-6}\msp{-5}}{N \eta}\,+\frac{1}{\!\sqrt{N}}
\,.
\]
For the bound on the off-diagonal error term we plug this result into \eqref{Regime 1 first error bounds} and get
\[
\Err_{\mathrm{o}}\,\prec\,|z|^{-2}\sqrt{\frac{\rho}{N\eta}}
+ \frac{\,|z|^{-6}\msp{-5}}{N\eta}+\frac{\,|z|^{-3}\msp{-5}}{N^{1/4}}\sqrt{\frac{1}{N\eta}}+\frac{\,|z|^{-2}\msp{-5}}{\!\sqrt{N}}\;.
\]

{\sc Regime 2: } Now let $|z|\leq N^{\15}$ and  
suppose that $\delta_*$ is a positive constant, depending only on the model parameters $p$, $P$ and $L$. 
The diagonal contribution, $\Err_{\mathrm{d}}$, satisfies
\bels{bulk bound}{
\Err_{\mathrm{d}}(z)\,\Ind\big(\Err_{\mathrm{d}}(z) \2\leq\2 \lambda_*/[z]\big)
\;\lesssim\;
[z]^{-2}\norm{\vect{d}(z)}_\infty\,,
}
according to \eqref{QVE stability away from support} in Lemma \ref{lmm:Stability far away from support} (for $|z|\geq 10$) and \eqref{rough stability} from Corollary \ref{crl:Stability away from minima} (for $|z|\leq 10$), where $\lambda_*$ is a sufficiently small positive constant.

We will now establish a gap in the possible values of $\Err(z) $ by showing (cf. \eqref{Err gap bound} below) that the right hand side of \eqref{bulk bound} is much less than $ \lambda_\ast/[z] $. 
To this end we  estimate the norm of $\vect{d}$ in \eqref{bulk bound} by Lemma \ref{lmm:Bound on perturbation} and also use the bound on the off-diagonal contribution, $\Err_{\mathrm{o}}$, from the same lemma,
\bels{Regime 2 first error bounds}{
\big(\1[z]^2\Err_{\mathrm{d}}+\norm{\vect{d}}_\infty\big)\,\Ind\big(\Err \2\leq\2\lambda_*/[z]\big)
\;&\prec\;
[z]^{-2}\sqrt{\frac{\Im \la \vect{g}\ra}{N \eta}}+\frac{1}{\sqrt{N}}
\,,
\\
[z]^2\Err_{\mathrm{o}}\,\Ind\big(\Err \2\leq\2\lambda_*/[z]\big)
\;
&\prec\;
\sqrt{\frac{\Im \la \vect{g}\ra}{N \eta}}+\frac{1}{\sqrt{N}}\,.
}
Now we use $ \Im\2\avg{\1\brm{g}\1} = \pi\1\rho  + \Im \2\avg{\1\vect{g}-\brm{m}\1} \lesssim \rho + \Err_{\mrm{d}} $ to estimate the first terms on the right hand side of \eqref{Regime 2 first error bounds}:
\[
\sqrt{\frac{\Im\2\avg{\1\brm{g}\1}}{N\eta}} 
\,\lesssim\, 
\sqrt{\frac{\pi\1\rho}{N\eta}} 
+ 
\sqrt{
\frac{1}{N\eta}\,\Err_{\mrm{d}}}
\,.
\]
Using again the weighted Cauchy-Schwarz inequality in the second term yields
\[
\big(\1[z]^2\Err_{\mathrm{d}}+\norm{\vect{d}}_\infty\big)\,\Ind\big(\Err \2\leq\2\lambda_*/[z]\big)
\,\prec\,
[z]^{-2}\sqrt{\frac{\rho}{N \eta}}+[z]^{-6}\frac{N^\eps}{N\eta}+\frac{1}{\sqrt{N}}+N^{-\eps}[z]^2\Err_{\mathrm{d}}
\,.
\]
The term $N^{-\eps}[z]^2\Err_{\mathrm{d}}$ can be absorbed  (cf. (ii) of Lemma~\ref{lmm:Basic facts about stochastic domination})  into the left hand side and we arrive at
\bels{Gap bound}{
\big(\1[z]^2\Err_{\mathrm{d}}+\norm{\vect{d}}_\infty\big)\,\Ind\big(\Err \2\leq\2\lambda_*/[z]\big)\,\prec\,
[z]^{-2}\sqrt{\frac{\rho}{N \eta}}+\frac{\;[z]^{-6}\msp{-6}}{N\eta}\,+\frac{1}{\!\sqrt{N\2}}\,.
}
For the off-diagonal error terms we plug this into the second bound of \eqref{Regime 2 first error bounds} after using $\Im \2\avg{\1\vect{g}\1} \lesssim \rho + \Err_{\mathrm{d}}$ and get
\bels{Gap bound off diagonal}{
\Err_{\mathrm{o}}
\;\prec\;
[z]^{-2}\sqrt{\frac{\rho}{N\eta}}
+\frac{\; [z]^{-6}\msp{-6}}{N\eta}\,+\frac{\;[z]^{-3}\msp{-6}}{N^{1/4}}\sqrt{\frac{1}{N\eta}}+\frac{\;[z]^{-2}\msp{-6}}{\!\sqrt{N}}\,.
}
In particular, we combine \eqref{Gap bound} and \eqref{Gap bound off diagonal} to establish a gap in the values that $\Err$ can take,
\bels{Err gap bound}{
\Err
\,\Ind\big(\Err \2\leq\2\lambda_*/[z]\big)\,\prec\,[z]^{-1}N^{-\gamma/2}
\,.
}
Here we used  $\eta \geq N^{\gamma-1}$.
This shows that either $ \Err \ge \lambda_\ast/[z] $ or $ \Err \leq N^{-\gamma/4}/[z] $ a.w.o.p.

Now we apply  Lemma \ref{lmm:Bound propagation}   on the connected domain 
\[
\setB{z \in \Cp:\2 \Im\1z \geq N^{\gamma-1},\2 \dist(z,\mathbb{M})\geq \delta_*,\2 |z|\leq N^5}
\,,
\]
 with the choices
\bels{}{
\varphi(z)\,:=\,[z]\1\Err(z)\,,\qquad \Phi(z)\,:=\, N^{-\gamma/3}\,, \qquad z_0\,:=\, \cI\2N^{\25}\,.
}
The continuity condition \eqref{weak continuity condition} of the lemma for these two functions follows from the H\"older-continuity, \eqref{Hoelder continuity of m}, of the solution of the QVE and the weak continuity of the resolvent elements, 
\bels{}{
|G_{i j}(z_1) - G_{i j}(z_2)| \,\leq \, \frac{|z_1-z_2|}{(\Im \2z_1) (\Im \2z_2)} \,\leq\, N^{\22} |z_1-z_2|\,.
}
The condition \eqref{initial value condition} holds since by \eqref{bulk error bound} on the first regime we have a.w.o.p. $\varphi(z_0)\leq \Phi(z_0)$. Finally, \eqref{Err gap bound} implies a.w.o.p. $\varphi\,\Ind (\1\varphi \in [\Phi-N^{-1}\!,\Phi\1]\2)< \Phi-N^{-1}$ and thus \eqref{gap condition}. We infer that a.w.o.p. $\varphi \leq \Phi$. In particular, the indicator function in \eqref{Gap bound} and \eqref{Gap bound off diagonal} is non-zero a.w.o.p..
Thus, \eqref{Gap bound} and \eqref{Gap bound off diagonal} imply \eqref{bulk error bound} in the second regime. 
\end{Proof}

We will now sketch the proof of Corollary \ref{crl:Local law in the bulk}. The set-up in this corollary differs slightly from the one used in the rest of this paper, because the uniform bound (assumption (C)) on the solution of \eqref{SCE} is not assumed.   We therefore use additional information from \cite{AEK1} about $\vect{m}$ in this more general setting. 

\begin{Proof}[Proof of Corollary \ref{crl:Local law in the bulk}] Since the boundedness assumption (C) on the solution of the QVE is dropped in this corollary, its proof starts by showing that nevertheless for some constant $P>0$ we have
\bels{boundedness on interval}{
|m_i(z)|\,\leq\, P\,,\qquad i=1,\dots,N\,,\; z \in I+ \cI\1(0,\infty)\,.
}
 In this setting the solution $\vect{m}(z)$ is not guaranteed to be extendable as a H\"older-continuous function with $N$-independent H\"older-norm to $z \in \overline{\Cp}$. The density of states, defined by \eqref{def of rho},  however still has a H\"older-norm with H\"older-exponent $1/13$ that is independent of $N$ (cf. (i) of Proposition 7.1 and (i) of Theorem 6.1 in \cite{AEK1}). 
Here, we used $L=1$ for the model parameter from assumption (B).  Furthermore, \eqref{boundedness on interval} follows from the lower bound on the density of states and (i) of Lemma~5.4 and (i) of Theorem~6.1 in \cite{AEK1}. 
For the proof of Proposition \ref{prp:Local law away from local minima} we only used the properties of the solution of \eqref{SCE}, valid for $z$ in the entire complex upper half plane, that are listed in Corollaries~\ref{crl:Bounds on solution}, \ref{crl:Stieltjes-transform representation} and \ref{crl:Stability away from minima}. 
These properties remain true for $\Re \2z \in I $ (cf. Theorem~2.1, (i) of Theorem~2.12, (i) of Proposition~5.3 and Proposition~7.1 \nc in \cite{AEK1}) if only \eqref{boundedness on interval} instead of \eqref{boundedness of m} is satisfied. Thus, \eqref{bulk error bound:diag}, \eqref{bulk error bound:off-diag} and \eqref{bulk average error bound} hold for $z \in I+ \cI\1[N^{\gamma-1},\infty)$ and Corollary \ref{crl:Local law in the bulk} is proven.
\end{Proof}

\section{Local law close to local minima}
\label{Main result at the edges}

\subsection{The solution of the QVE}
In this subsection we state a few facts about the solution $\vect{m}$ of the QVE \eqref{SCE} and about the stability of this equation against perturbations. These facts are summarized in two theorems that are taken from the companion paper \cite{AEK1}. The first theorem contains regularity properties of $\vect{m}$. Furthermore, it provides lower and upper bounds on the imaginary part, $\Im \la\vect{m}\ra=\pi \2\rho$, by explicit functions. It is a combination of the statements from Theorem~2.1, Theorem~2.4, Theorem~2.6 and Corollary A.1 of \cite{AEK1}.

\begin{theorem}[Solution of the QVE]
\label{thr:Solution of the QVE}
Let the sequence $\vect{S} = \vect{S}^{(N)}$ satisfy the assumptions \emph{(A)-(C)}. Then for every component, $m_i: \Cp \to \Cp$, of the unique solution, $\vect{m}=(m_1, \dots, m_N)$, of the QVE there is a probability density $p_i:\R\to [0,\infty)$ with support in the interval $[-2,2]$, such that
\bels{}{
m_i(z)\,=\, \int_\R \frac{p_i(\tau)\1\dif \tau}{\tau-z}\,,\qquad z \in \Cp\,,\; i=1,\dots,N\,.
}
The probability densities are comparable,
\bels{}{
p_i(\tau)\,\sim\,p_j(\tau)\,,\qquad \tau \in \R\,, \; i,j=1,\dots,N\,.
}
The solution $\vect{m}$ has a uniformly H\"older-continuous extension (denoted again by $\vect{m}$) to the closed complex upper half plane $\overline \Cp=\Cp\cup\R$,
\bels{}{
\norm{\1\vect{m}(z_1)-\vect{m}(z_2)}_\infty\,\lesssim\, |z_1-z_2|^{1/3}\,, \qquad z_1,z_2 \in \overline\Cp\,.
}
Its absolute value satisfies
\[
|m_i(z)|\,\sim\, [z]^{-1}\,,\qquad z \in \overline{\Cp}\,,\; i=1,\dots,N\,.
\]
Let $\rho: \R \to [0,\infty),\tau \mapsto \avg{\1\brm{p}(\tau)}$ be the density of states, defined in \eqref{def of rho}. 
Then there exists a positive constant $\delta_* \sim 1 $ such that the following holds true. 
The support of the density consists of $K \sim 1 $ disjoint intervals of lengths at least $2\1\delta_*$, i.e., 
\bels{defining property of of alpha_i and beta_i}{
\supp \rho \,=\, \bigcup_{i=1}^{K}\, [\1\alpha_i,\beta_i]
\,,\quad\text{where}\quad
\beta_i -\alpha_i \ge 2\1\delta_*
\,,\quad
\text{and}\quad 
\alpha_i < \beta_i< \alpha_{i+1}
\,.
}
The size of the harmonic extension \eqref{harmonic extension of rho} of $\rho$, up to constant factors, is given by explicit functions as follows. 
 Let $\eta \in [0,\delta_*]$.
\begin{subequations}
\label{qualitative size of rho}
\begin{itemize}
\item {\bf Bulk:} Close to the support of the density of states but away from the local minima in $\mathbb{M}$ (cf. \eqref{definition of M}) the function $\rho$ is comparable to $1$, i.e.,
\bels{lower bound for rho when bounded away from MM}{
\rho(\tau+\cI\1\eta)\,\sim\, 1\,,\qquad \tau \in \supp \rho\,,\; \dist(\tau,\mathbb{M})\geq \delta_*\,.
}
\item {\bf At an internal edge:} 
At the edges $\alpha_i, \beta_{i-1}$ with $i =2, \dots, K$ in the direction where the support of the density of states continues the size of $\rho$ is
\bels{rho at an internal edge}{
\rho(\alpha_i+\omega + \cI\2\eta)
\,\sim\,
\rho(\beta_{i-1}-\omega + \cI\2\eta)
\,\sim\, 
\frac{(\omega+\eta)^{1/2}}{(\alpha_i-\beta_{i-1}+\omega+\eta)^{1/6}}\,,\qquad \omega \in [\10\1,\delta_*]\,.
}
\item 
{\bf Inside a gap:} 
Between two neighboring edges $\beta_{i-1} $ and $\alpha_i$ with $i =2, \dots, K$,
the function $\rho$ satisfies
\bels{rho inside a gap}{
\rho(\beta_{i-1} + \omega+ \cI\2\eta)\,\sim\,\rho(\alpha_{i} - \omega+ \cI\2\eta)
\,\sim\,\frac{\eta}{(\alpha_i-\beta_{i-1}+\eta)^{1/6}(\omega+\eta)^{1/2}}\,, 
}
for all $\omega \in [0,(\alpha_i-\beta_{i-1})/2]$.
\item 
{\bf Around an extreme edge:} At the extreme points $ \alpha_1 $ and $ \beta_{K} $ of $ \supp\rho $ the density of states grows like a square root ,
\bels{rho around extreme edge}{
\rho(\alpha_1+\omega + \cI\2\eta)
\,\sim\,
\rho(\beta_{K}-\omega + \cI\2\eta)
\,\sim\,
\begin{cases}
\displaystyle
(\omega+\eta\1)^{1/2}\,,\quad	
&\omega \in [\10,\delta_*\1]\,,
\\
\displaystyle
\frac{\eta}{(\1\abs{\omega}+\eta)^{1/2}\!}\;,	
&\omega \in [-\delta_*,\10\1]\,.
\end{cases}
}
\item {\bf Close to a local minimum:} 
In a neighborhood of a local minimum in the interior of the support of the  density of states, i.e., for $\tau_0 \in \mathbb{M} \cap {\rm int}\supp \rho$, we have
\bels{rho close to a local minimum}{
\rho(\tau_0+\omega + \cI\2\eta)
\,\sim\, 
\rho(\tau_0)+ (|\omega|+\eta)^{1/3}\,,\qquad \omega \in [-\delta_*,\2\delta_*\1]\,.
}
\item {\bf Away from the support:} Away from the interval in which $\supp \rho$ is contained
\bels{rho away from supp}{
\rho(z)\,\sim\, \frac{\Im\2z}{|z|^2\!}\,,\qquad z \in \overline \Cp\,,\;\dist(z,[\alpha_1,\beta_K])\,\geq\,\delta_*\,.
}
\end{itemize}
\end{subequations}
\end{theorem}

The next theorem shows that the QVE is stable under small perturbations, $\vect{d}$, in the sense that once a solution of the perturbed QVE \eqref{perturbed QVE} is sufficiently close to $\vect{m}$, then the difference between the two can be estimated in terms of $\norm{\vect{d}}_\infty$. In \cite{AEK1} it is stated as  Proposition~10.1.

\begin{theorem}[Stability]
\label{thr:Stability around support}
There exists a scalar function $\sigma: \overline{\Cp} \to [\10\1,\infty)$, three vector valued functions $\vect{s}, \vect{t}^{(1)}, \vect{t}^{(2)}: \overline \Cp \to \C^N$, a matrix valued function $\vect{T}:\overline \Cp\to \C^{N\times N}$, all depending only on $\vect{S}$, and a positive constant $\lambda_*$, depending only on the model parameters $p$, $P$ and $L$,
such that for two arbitrary vector valued functions $\vect{d}: \Cp \to \C^N$ and $\vect{g}: \Cp \to (\C\backslash\sett{0})^N$ that satisfy 
\bels{perturbed QVE}{
-\frac{1}{g_i(z)} \,=\, z + \sum_{j=1}^N s_{i j} g_j(z) +d_i(z)\,, \quad z \in \Cp\,,
}
the difference between $\vect{g} = \brm{g}(z)$ and $\vect{m} = \brm{m}(z)$ is bounded in terms of
\bels{Def of Theta}{
\Theta = \Theta(z)\,:=\, \absb{\avgb{\2\vect{s}(z)\1,\2 \vect{g}(z)-\vect{m}(z)}}
\,,\qquad z \in \Cp\,,
}
in the following two ways. 
On the whole complex upper half plane
\begin{align}
\label{Bound of Err in terms of Theta}
\norm{\1\vect{g}-\vect{m}\1}_\infty
\2\Ind\big(\2\norm{\vect{g}-\vect{m}}_\infty\leq \lambda_*\big)
\;&\lesssim\; \Theta\,+\,\norm{\vect{d}}_\infty
\,,
\\
\label{bound on averages against bounded vectors}
|\avg{\1\vect{w},\vect{g}-\vect{m}}|
\2\Ind\big(\2\norm{\vect{g}-\vect{m}}_\infty\leq \lambda_*\big)
\;&\lesssim\; 
\norm{\vect{w}}_\infty\2\Theta
+\norm{\vect{w}}_\infty\norm{\vect{d}}_\infty^2+|\avg{\1\vect{T}\vect{w},\vect{d}\1}|
\,,
\end{align}
for any non-random $\brm{w} \in \C^N $.
The scalar function $\Theta : \overline{\Cp} \to [\10\1,\infty) $ satisfies a cubic equation
\bels{cubic for Theta}{
\absb{\2\Theta^3+\pi_2\1\Theta^2+\pi_1\1\Theta\2}\2\Ind\big(\2\norm{\vect{g}-\vect{m}}_\infty\leq \lambda_*\big)
\;\lesssim\; \norm{\vect{d}}_\infty^2+\abs{\avg{\1\vect{t}^{(1)}\!,\vect{d}\1}}\,+\,\abs{\avg{\1\vect{t}^{(2)}\!,\vect{d}\1}}
\,.
}
The coefficients $\pi_1,\pi_2:\Cp \to \C$ may depend on $\brm{S}$ and $\brm{g}$. They satisfy
\begin{subequations}
\label{comparison relation coefficients}
\begin{align}
\label{comparison relation for pi_1}
|\pi_1(z)|
\,&\sim\, 
\frac{\Im\2z}{\rho(z)}+\rho(z)\2(\1\rho(z)+\sigma(z))
\,,
\\
\label{comparison relation for pi_2}
|\pi_2(z)|
\,&\sim\, 
\rho(z)+\sigma(z)
\,,
\end{align}
\end{subequations}
for all $ z \in \Cp $. 
Moreover, the functions $\sigma$, $\vect{s}$, $\vect{t}^{(1)}\!, \vect{t}^{(2)}$ and $\vect{T}$ are regular in the sense that
\begin{align}
\label{sigma and s continuity}
\abs{\1\sigma(z_1)-\sigma(z_2)} 
\,+\, 
\norm{\1\vect{s}(z_1)-\vect{s}(z_2)} 
\;\lesssim\; 
\abs{z_1-z_2}^{1/3}
\,,
\qquad 
z_1,z_2 \in \overline{\Cp}
\,,
\\
\label{bounds on the auxiliary functions}
\sigma(z)
\,+\,
\norm{\1\vect{s}(z)}_\infty
\,+\,
\norm{\1\vect{t}^{(1)}(z)}_\infty
\,+\,
\norm{\1\vect{t}^{(2)}(z)}_\infty
\,+\,
\norm{\vect{T}(z)}_{\infty\to \infty}
\;\lesssim\; 
1\,,
\qquad z \in \overline{\Cp}
\,.
\end{align}
Furthermore, the function $ \sigma $ is related to the density of states by
\begin{subequations}
\label{sigma at minima}
\begin{align}
\sigma(\alpha_i)
\,&\sim\, 
\sigma(\beta_{i-1}) \,\sim\,(\alpha_i-\beta_{i-1})^{1/3}\,,
\qquad i=2,\dots,K\,, 
\\
\sigma(\alpha_1)\,&\sim\, \sigma(\beta_{K}) \,\sim\,1\,, 
\\
\sigma(\tau_0)\,&\lesssim\; \rho(\tau_0)^2\,, \qquad\qquad \tau_0 \in \mathbb{M}\backslash\{\alpha_i,\beta_i\}
\,.
\end{align}
\end{subequations}
\end{theorem}

We warn the reader that in this paper $ \Theta $ and $ \sigma $ denote  the absolute values of the quantities denoted by the same symbols in Proposition~10.1 of \cite{AEK1}. 
The function $\sigma$ appears naturally in the analysis of the QVE. Analogous to the more explicitly constructed function $\Delta$ from Definition~\ref{def:Local gap size}, at an edge the value of $\sigma^{\13}$ encodes the size of the corresponding gap in $\supp\rho$. At the local minima in $\mathbb{M}\backslash\{\alpha_i,\beta_i\}$ the value of $\sigma^{\13}$ is small, provided the density of states has a small value at the minimum. In this sense it is again analogous to $\Delta$, which vanishes at these internal minima.

\subsection{Coefficients of the cubic equation}

The stability of QVE near the points in $\mathbb{M}$ requires a careful analysis of the cubic equation \eqref{cubic for Theta} for $\Theta$ from 
Theorem \ref{thr:Stability around support}. For this, we will provide a more explicit description of the upper and lower bounds from \eqref{comparison relation coefficients} on the coefficients, $\pi_1$ and $\pi_2$, of the cubic equation.

\begin{proposition}[Behavior of the coefficients]
\label{prp:Behavior of the coefficients}
There exist $ \delta_\ast,c_\ast \sim 1 $ such that for all $\eta \in [0,\delta_*]$ the coefficients, $\pi_1$ and $\pi_2$, of the cubic equation \eqref{cubic for Theta} satisfy the following bounds. 
\begin{subequations}
\begin{itemize}
\item {\bf Around an internal edge:} 
At the edges $\alpha_i, \beta_{i-1}$ of the gap with length $\Delta:=\alpha_i-\beta_{i-1}$ for $i =2, \dots, K$, we have 
\bels{coefficients around an internal edge}{
\abs{\1\pi_1(\alpha_i+\omega + \cI\2\eta)}
\,&\sim\,
\abs{\1\pi_1(\beta_{i-1}-\omega + \cI\2\eta)}
\,\sim\, 
(\1\abs{\omega}+\eta)^{1/2}(\1\abs{\omega}+\eta+\Delta)^{1/6},
\\
\abs{\1\pi_2(\alpha_i+\omega + \cI\2\eta)}
\,&\sim\,
\abs{\1\pi_2(\beta_{i-1}-\omega + \cI\2\eta)}
\,\sim\, 
(\1\abs{\omega}+\eta+\Delta\1)^{1/3},
\quad \omega \in [-c_*\Delta,\delta_*]
\,.
}
\item 
{\bf Well inside a gap:} 
Between two neighboring edges $\beta_{i-1} $ and $\alpha_i$ of the gap with length $\Delta:=\alpha_i-\beta_{i-1}$ for $i =2, \dots, K$, the first coefficient, $\pi_1$, satisfies
\bels{pi1 well inside a gap}{
\abs{\1\pi_1(\alpha_i-\omega + \cI\2\eta)}
\,\sim\,
\abs{\1\pi_1(\beta_{i-1}+\omega + \cI\2\eta)}
\,\sim\, 
(\1\eta+\Delta)^{2/3}\,,
\qquad\omega \in \biggl[c_*\Delta\1,\2\frac{\Delta}{2}\biggr]
\,.
}
The second coefficient, $\pi_2$, satisfies the upper bounds,
\bels{pi2 well inside a gap}{
\begin{array}{c}
\msp{10}\abs{\1\pi_2(\alpha_i-\omega + \cI\2\eta)}
\;\lesssim\;
(\1\eta+\Delta)^{1/3}\,,
\vspace{0.2cm}
\\
\abs{\1\pi_2(\beta_{i-1}+\omega + \cI\2\eta)}
\;\lesssim\;
(\1\eta+\Delta)^{1/3}\,,
\end{array}
\qquad\omega \in \biggl[c_*\Delta\1,\2\frac{\Delta}{2}\biggr]
\,.
}
\item 
{\bf Around an extreme edge:} Around the extreme points $ \alpha_1 $ and $ \beta_{K} $ of $ \supp \rho $, we have 
\bels{coefficients around an extreme edge}{
\begin{array}{l}
\abs{\1\pi_1(\alpha_1+\omega + \cI\2\eta)}
\,\sim\,
\abs{\1\pi_1(\beta_{K}-\omega + \cI\2\eta)}
\,\sim\,
(\omega+\eta\1)^{1/2}
\vspace{0.2cm}
\\
\abs{\1\pi_2(\alpha_1+\omega + \cI\2\eta)}
\,\sim\,
\abs{\1\pi_2(\beta_{K}-\omega + \cI\2\eta)}
\,\sim\,
1\,,
\end{array}
\qquad\omega \in [-\delta_*,\delta_*]\,.
}
\item {\bf Close to a local minimum:} 
In a neighborhood of the local minimum in the interior of the support of the  density of states, i.e. for $\tau_0 \in \mathbb{M} \cap {\rm int}\supp \rho$, we have
\bels{coefficients close to a local minimum}{
\begin{array}{l}
\abs{\1\pi_1(\tau_0+\omega + \cI\2\eta)}
\,\sim\, 
\rho(\tau_0)^2+ (|\omega|+\eta)^{2/3}\,,
\vspace{0.2cm}
\\
\abs{\1\pi_2(\tau_0+\omega + \cI\2\eta)}
\,\sim\, 
\rho(\tau_0)+ (|\omega|+\eta)^{1/3}\,,
\end{array}
\qquad \omega \in [-\delta_*,\2\delta_*\1]\,.
}
\end{itemize}
\end{subequations}
\end{proposition}
\begin{Proof} The proof is split according to the cases above. In each case we combine the general formulas \eqref{comparison relation coefficients} with the knowledge about the harmonic extension, $\rho$, of the density of states from Theorem \ref{thr:Solution of the QVE} and about the behavior of the positive H\"older-continuous function, $\sigma$, at the minima in $\mathbb{M}$ from \eqref{sigma at minima}. The positive constant $\delta_*$ is chosen to have at most the same value as in Theorem \ref{thr:Solution of the QVE}. We start with the simplest case.
\\[0.3cm]
{\sc Around an extreme edge: } By the H\"older-continuity of $\sigma$ (cf. \eqref{sigma and s continuity}) and because $\sigma$ is comparable to $1$ at the points $\alpha_1$ and $\beta_K$ (cf. \eqref{sigma at minima}), this function is comparable to $1$ in the whole $\delta_*$-neighborhood of the extreme edges. Thus, using \eqref{comparison relation coefficients} inside this neighborhood, we find 
\[
|\1\pi_1(z)|\,\sim\, \frac{\Im\2z}{\rho(z)}+\rho(z)\,,\qquad |\1\pi_2(z)|\,\sim\, 1\,.
\]
The claim now follows from the behavior of $\rho$, given in Theorem \ref{thr:Solution of the QVE}, inside this domain. 
\\[0.3cm]
{\sc Close to a local minimum: } In this case $\rho+\sigma$ is comparable to $\rho$. In fact, using the $1/3$-H\"older-continuity of $\sigma$ (cf. \eqref{sigma and s continuity}) and its bound at the minimum, $\tau_0\in \mathbb{M}$, (cf. \eqref{sigma at minima}) we find 
\bels{}{
\rho(z)\,\leq\, \rho(z)+\sigma(z) \,\lesssim\, \rho(z)+\rho(\tau_0)^2+ |z-\tau_0|^{1/3}\,\sim\,  \rho(z)\,,\qquad |z-\tau_0|\leq \delta_*\,.
}
In the last relation we used the behavior \eqref{rho close to a local minimum} of $\rho$ from Theorem \ref{thr:Solution of the QVE}. By \eqref{comparison relation coefficients} we conclude that inside the $\delta_*$-neighborhood of $\tau_0$,
\bels{}{
|\1\pi_1(z)|\,\sim\, \frac{\Im\2z}{\rho(z)}+\rho(z)^2 \,,\qquad |\pi_2(z)|\,\sim\, \rho(z)\,.
}
Using the upper and lower bounds on $\rho(z)$ again, gives the desired result, \eqref{coefficients close to a local minimum}.
\\[0.3cm]
{\sc Around an internal edge: } First we prove the bounds on $|\pi_2|$, starting from \eqref{comparison relation coefficients}. The upper bound simply uses the $1/3$-H\"older-continuity and the behavior at the edge points of $\sigma$,
\bels{upper bound on second coefficient}{
|\pi_2(z)|\,\sim\, \rho(z)+\sigma(z)\,\lesssim\, \rho(z)+ \Delta^{\!1/3}+|z-\tau_0|^{1/3},
}
where $\tau_0$ is one of  the edge points $\alpha_i$ or $\beta_{i-1}$. The claim follows from plugging in the size of $\rho$ from the two corresponding domains in Theorem \ref{thr:Solution of the QVE}, i.e., the domain close to an edge, \eqref{rho at an internal edge}, and the domain inside a gap, \eqref{rho inside a gap}.

For the lower bound we consider two different regimes. In the first case $z$ is close to the edge point, $|z-\tau_0| \leq c\1\Delta$, for some small positive constant $c$, depending only on the model parameters $p$, $P$ and $L$. We find
\[
|\1\pi_2(z)|\,\sim\, \rho(z)+\sigma(z)\,\gtrsim\, \rho(z)+ \Delta^{\!1/3} -C\1|z-\tau_0|^{1/3}\,\sim\,\rho(z)+ \Delta^{\!1/3}\,,
\]
provided $c$ is small enough. 
This bound coincides with the lower bound on $\pi_2$ in \eqref{coefficients around an internal edge}, once the size of $\rho$ from \eqref{rho at an internal edge} is used.

In the second regime, $|z-\tau_0| \geq c\1\Delta$, we simply use $|\pi_2(z)|\gtrsim\rho(z)$ from \eqref{comparison relation coefficients}. If $\Re\1z \in \supp \rho$, then
the size of $\rho$ from \eqref{rho at an internal edge}
yields the desired lower bound. If, on the other hand, $\Re\1z$ lies inside a gap of $\supp \rho$, then we use the freedom of choosing the constant $c_*$ in Proposition \ref{prp:Behavior of the coefficients}. Suppose $c_*\leq c/2$. Then $|z-\tau_0| \geq c\1\Delta$ and $|\Re \1z -\tau_0|\leq c_*\Delta$ imply $\Im\1z \gtrsim \Delta$ and 
\[
\rho(z)\,\sim\,(\Im\1z)^{1/3}\,\gtrsim\,\Delta^{\!1/3}+|z-\tau_0|^{1/3}.
\]

This finishes the proof of the upper and lower bound on $|\pi_2|$ on this domain. 
For the claim about $|\pi_1|$ we plug the result about $|\pi_2|$ and the size of $\rho$ into
\bels{comparison of pi1 and pi2}{
|\1\pi_1|\,\sim\, \frac{\Im\1z}{\rho(z)}+ \rho(z)\2|\pi_2(z)|
\,.
}
{\sc Well inside a gap: } For the upper bound on $|\pi_2|$ we simply use \eqref{upper bound on second coefficient} again, which follows from \eqref{sigma and s continuity} and \eqref{sigma at minima}.  
The comparison relation for $|\pi_1|$ now follows from \eqref{comparison of pi1 and pi2} again. For the lower bound, $|\pi_1|\gtrsim \Im\1z /\rho$ and \eqref{rho inside a gap} from Theorem \ref{thr:Solution of the QVE} are sufficient.  
This finishes the proof of the proposition.
\end{Proof}

\subsection{Rough bound on $\Err$ close to local minima}

In the following lemma we will see that a.w.o.p. $\Lambda \leq c$ for some arbitrarily small constant $c>0$. Since the local law away from $\mathbb{M}$ is already shown in Proposition \ref{prp:Local law away from local minima} we may restrict to bounded $z$ in the following. From here on until the end of Section \ref{Main result at the edges} we assume $|z|\leq 10$.

\NLemma{Rough bound}{Let $\lambda_*$ be a positive constant. Then, uniformly for all $z= \tau + \cI\1\eta \in \Cp$ with $\eta \geq N^{\gamma-1}$, the function $\Err$ is uniformly small,
\bels{}{
\Err(z)\, \leq\, \lambda_*\qquad \text{a.w.o.p.}\,.
}
}
\begin{Proof}
Away from the local minima in $\mathbb{M}$ the claim follows from \eqref{bulk error bound} in Proposition \ref{prp:Local law away from local minima}. We will therefore prove that $\Err$ is smaller than any fixed positive constant in some $\delta$-neighborhood of $\mathbb{M}$. We will use the freedom to choose the size $\delta \sim 1 $ of these neighborhoods as small as we like. 

Let us sketch the upcoming argument.
Close to the points in $\mathbb{M}$ we make use of Theorem \ref{thr:Stability around support}. Using Lemma \ref{lmm:Bound on perturbation}, we will see that the right hand side of the cubic equation in $\Theta$, \eqref{cubic for Theta}, is smaller than a small negative power, $N^{-\eps}$, of $N$, provided $\Err$ is bounded by a small constant, $\Err\leq \lambda_*$. This will imply that $\Theta$ itself is small and through \eqref{Bound of Err in terms of Theta} that the bound on $\Err$ can be improved to $\Err\leq \lambda_*/2$. In this way we establish a gap in the possible values that the continuous function $\Err$ can take. Lemma \ref{lmm:Bound propagation} in the appendix is then used to propagate the bound on $\Err$ from Proposition \ref{prp:Local law away from local minima} into the $\delta$-neighborhoods of the points in $\mathbb{M}$.

Now we start the detailed proof from the fact that $\Theta$ satisfies the cubic equation \eqref{cubic for Theta}, whose right hand side is bounded by $C\norm{\vect{d}}_\infty$ for some constant $C$, depending only on the model parameters. Note that $\norm{\vect{d}}_\infty\lesssim 1$ as long as $\Err\leq \lambda_*$ because in this case $|m_i|\sim 1$, $|g_i|\sim 1$ and $\vect{g}$ satisfies the perturbed QVE with perturbation $\vect{d}$.
From the definition of $\Theta$ in \eqref{Def of Theta} and the uniform bound on $\vect{s}$ from \eqref{bounds on the auxiliary functions}, we get $\Theta\lesssim \Err$. 
Since the coefficient $|\pi_2|$ is uniformly bounded (cf. \eqref{comparison relation coefficients}), the cubic equation for $\Theta$ implies the three bounds
\begin{subequations}
\label{Theta cubic equation bounds}
\begin{align}
\label{first Theta cubic equation bound}
\Theta\,\Ind(\2\Err\leq \eps_1, \2|\pi_1|\geq C_1\eps_1)
\;&\lesssim\; 
\frac{\2\norm{\vect{d}}_\infty\!}{|\1\pi_1|}\;,
\\
\label{second Theta cubic equation bound}
\Theta\,\Ind(\2\Err\leq \eps_2, \2|\pi_2|\geq C_2\eps_2)
\;&\lesssim\;\frac{|\pi_1|}{|\pi_2|}+\frac{\norm{\vect{d}}_\infty^{1/2}\!}{|\pi_2|^{1/2}\!}\;,
\\
\label{}
\Theta\,\Ind(\2\Err\leq \lambda_*)
\;&\lesssim\; 
|\pi_2|+ \sqrt{|\pi_1|}+\norm{\vect{d}}_\infty^{1/3}
\,.
\end{align}
\end{subequations}
Here, $\eps_1,\eps_2 \in (0,\lambda_\ast) $, with $ \lambda_\ast \in (0,1) $ from Theorem~\ref{thr:Stability around support}, are arbitrary constants and $ C_1,C_2>0 $ depend only on the model parameters. 
We prove \eqref{second Theta cubic equation bound}; the other two bounds are obtained similarly.
%
%
First we show that under the assumptions $ \Err\leq \eps_2 $ and  $\abs{\pi_2} \ge C_2\eps_2 $ the second order term $ \pi_2\Theta^2 $ is at least three times larger than $ \Theta^3 $ provided $ C_2 \sim 1 $ is chosen to be sufficiently large.
Indeed, since $ \Theta \leq \norm{\brm{s}}_\infty \Err \leq \norm{\brm{s}}_\infty\eps_2 $ and $ \abs{\pi_2} \ge C_2\eps_2 $, it suffices to choose $ C_2 \ge 3\1\norm{\brm{s}}_\infty \sim 1 $. Here we have also used  \eqref{bounds on the auxiliary functions}.
Next we compare the second order term to the linear term $ \pi_1\Theta $. 
We may assume that $\Theta \ge 3\2\abs{\pi_1/\pi_2} $, otherwise \eqref{second Theta cubic equation bound} holds trivially.
%
%
%
Together with $ \abs{\pi_2}\Theta^2 \ge 3\2\Theta^3 $ proved above this implies that the second order term $ \pi_2\Theta^2 $ dominates the left hand side of \eqref{cubic for Theta}. 
Combining this with 
$ \abs{\1\avg{\2\brm{t}^{(j)},\brm{d}\1}} \lesssim \norm{\brm{d}}_\infty  $ (cf. \eqref{bounds on the auxiliary functions}) on the right hand side of \eqref{cubic for Theta}, hence yields
\bels{getting second Theta cubic equation bound}{
\frac{1}{3}\abs{\pi_2}\Theta^2 \,\leq\, 
\abs{\1\Theta^3+\pi_2\Theta^2+\pi_1\Theta\1}
\,\lesssim\, \norm{\brm{d}}_\infty 
\,.
}
In order to satisfy the constraint of \eqref{cubic for Theta} we have also used  $ \eps_2 \leq \lambda_\ast $.
%
%
This together with \eqref{getting second Theta cubic equation bound} yields \eqref{second Theta cubic equation bound}.\nc

Let $\delta \in (0,1)$ be another constant to be chosen later which  depends only on the model parameters $p$, $P$, and $L$. We split $\mathbb{M}$ into four subsets, which are treated separately,
\begin{equation*}
\begin{split}
&\mathbb{M}_1(\delta)
\,:=\,
\setb{\tau_0 \in \mathbb{M}\backslash\partial \supp \rho :  \rho(\tau_0)> \delta^{1/3}}
\,,
\quad
\mathbb{M}_2(\delta)
\,:=\, 
\setb{\tau_0 \in \partial \supp \rho:  \Delta(\tau_0)> \delta^{1/2}}
\,,
\\
&\mathbb{M}_3(\delta)
\,:=\,
\setb{\tau_0 \in \mathbb{M}\backslash \partial \supp \rho:  \rho(\tau_0)\leq \delta^{1/3}}\,,
\quad
\mathbb{M}_4(\delta)
\,:=\,
\setb{\tau_0 \in \partial \supp \rho :   \Delta(\tau_0)\leq \delta^{1/2} }
\,.
\end{split}
\end{equation*}
The function $\Delta$ is from Definition \ref{def:Local gap size} and its value is simply the length of the gap at the point $\tau_0 \in \partial \supp \rho$ where it is evaluated. We also define the $\delta$-neighborhoods of these subsets,
\[
\D_k(\delta)\,:=\, \setb{z \in \Cp:\; \dist(z,\mathbb{M}_k(\delta))\leq \delta\2}\,,\qquad 
k=1,2,3,4\,.
\]
As an immediate consequence of the upper and lower bounds on the coefficients, $\pi_1$ and $\pi_2$, presented in Proposition \ref{prp:Behavior of the coefficients}, we see that
\begin{subequations}
\label{coefficient bounds in terms of delta}
\begin{eqnarray}
\label{coefficient bounds on D1}
&|\pi_1(z)|\,\gtrsim \,\delta^{2/3}\,,\qquad   & z \in \D_1(\delta)\,,
\\
\label{coefficient bounds on D2}
& |\pi_1(z)|\,\lesssim\, \delta^{1/2}\,,\qquad |\pi_2(z)|\,\gtrsim \,\delta^{1/6}\,,\qquad & z \in \D_2(\delta)\,,
\\
&|\pi_1(z)|\,\lesssim \,\delta^{1/2}\,,\qquad |\pi_2(z)|\,\lesssim\,\delta^{1/6}\,,\qquad & z \in\D_3(\delta)\cup\D_4(\delta)\,.
\end{eqnarray}
\end{subequations}
On $\D_2(\delta)$ only the regimes around an internal edge, \eqref{coefficients around an internal edge}, and around an extreme edge, \eqref{coefficients around an extreme edge}, are relevant. The case well inside the gap, \eqref{pi1 well inside a gap} and \eqref{pi2 well inside a gap}, does not apply for small enough $\delta$, since $\Delta(\tau_0)>\delta^{1/2}$ but $|z-\tau_0|\leq \delta$. 

Now we make a choice for the two constants $\eps_1$ and $\eps_2$. We express them in terms of $\delta$ as
\[
\eps_1\,:=\,\delta\,,\qquad \eps_2\,:=\, \delta^{1/5}\,.
\]
We pair the bounds on $\Theta$ from \eqref{Theta cubic equation bounds} with the corresponding bounds from \eqref{coefficient bounds in terms of delta} on the coefficients of the cubic equation. For small enough $\delta$ the conditions on $\pi_1$ in \eqref{first Theta cubic equation bound} and $\pi_2$ in \eqref{second Theta cubic equation bound} are automatically satisfied by the choice of $\eps_1$ and $\eps_2$, as well as the upper and lower bounds from \eqref{coefficient bounds on D1} and \eqref{coefficient bounds on D2}. Thus, for small enough $\delta$ we end up with 
\begin{align*}
\Theta(z)\1\Ind(\2\Err(z)\leq \delta\2)
\;&\lesssim\;
\delta^{-2/3}\2\norm{\1\vect{d}(z)}_\infty\,, \msp{82} z \in \D_1(\delta)\,,
\\
\Theta(z)\1\Ind(\2\Err(z)\leq \delta^{1/5}\1)
\;&\lesssim\; 
\delta^{1/3}+\delta^{-1/12}\norm{\1\vect{d}(z)}_\infty^{1/2}
\,, 
\quad z \in \D_2(\delta)\,,
\\
\Theta(z)\1\Ind(\2\Err(z)\leq \lambda_*)
\;&\lesssim\; 
\delta^{1/6}+\norm{\1\vect{d}(z)}_\infty^{1/3}
\,,\msp{66}
z \in\D_3(\delta)\cup \D_4(\delta)
\,.
\end{align*}
At this stage we use Lemma \ref{lmm:Bound on perturbation} in the form of $\norm{\vect{d}}_\infty \prec N^{-\gamma/2}$ on the set where $\Err \leq \lambda_*/10$, say, and \eqref{Bound of Err in terms of Theta} from Theorem \ref{thr:Stability around support}. We may choose $\lambda_*$ to be sufficiently small compared to the constants with the same name from these two statements. Furthermore, we choose $\delta$ so small that $\delta^{1/5}\leq \lambda_*$. 
Since $ \norm{\brm{d}}_\infty \leq N^{-\gamma/2+c} $ a.w.o.p. for an arbitrary $ c > 0 $ we obtain 
\begin{subequations}
\label{Err gap for rough bound}
\begin{align}
\label{}
\Err(z) \,\Ind(\Err(z)\leq \delta\1)\,
&\lesssim\; 
\delta^{-2/3}N^{-\gamma/3}
\,,\msp{77} 
z \in \D_1(\delta)
\,,
\\
\text{a.w.o.p.}
\quad 
\Err(z)\,\Ind(\Err(z)\leq \delta^{1/5})
\,&\lesssim\; 
\delta^{1/3}+\delta^{-1/12}N^{-\gamma/5}
\,, \quad z \in \D_2(\delta)
\,,
\\
\Err(z)\,\Ind(\Err(z)\leq \lambda_*)
\,&\lesssim\; 
\delta^{1/6}+N^{-\gamma/7}
\,,\msp{66} 
z \in\D_3(\delta)\cup\D_4(\delta)
\,.
\end{align}
\end{subequations}
The right hand sides, including the constants from the comparison relation, can be made smaller than any given constant $\lambda_*$ by choosing $\delta=\delta_*$, depending only on the model parameters, small enough and $N$ sufficiently large. 
Furthermore, \eqref{Err gap for rough bound} establish a gap in the possible values that $\Err$ can take on the $\delta_*$-neighborhood of any point in $\mathbb{M}$. By Proposition \ref{prp:Local law away from local minima} we have the bound $\Err \prec N^{-\gamma/2}$ outside these $\delta_*$-neighborhoods and thus also for at least one point in the boundary of each neighborhood. Now we apply Lemma \ref{lmm:Bound propagation} to each neighborhood and in this way we propagate the bound $\Err\leq \lambda_*$ to every point $z$ in the $\delta_*$-neighborhood of $\mathbb{M}$ with $\Im\1z \geq N^{\gamma-1}$.
\end{Proof}

\subsection{Proof of Theorem \ref{thr:Local law}}
\label{ssec:Proof of Theorem LL}

According to Proposition \ref{prp:Local law away from local minima} the local law, Theorem \ref{thr:Local law}, holds outside the $\delta_*$-neighborhoods of the points in $\mathbb{M}$. It remains to show that it is true inside these neighborhoods as well. 
From here on we assume that $z\in \Cp$ satisfies $\dist(z,\mathbb{M})\leq \delta_*$ and $\Im\1z \geq N^{\gamma-1}$. Let $\tau_0\in \mathbb{M}$ be one of the closest points to $z$ in $\mathbb{M}$, i.e.,
\[
|z-\tau_0|\,=\, \dist(z,\mathbb{M})\2.
\]
When $ \tau_0 \in \partial \supp \rho $ we denote by $ \theta = \theta(\tau_0) \in \sett{\pm 1} $ the direction that points towards the gap in $ \supp \rho $ at $ \tau_0 $. 
In case $\tau_0 \notin \partial\supp\rho$  we make the arbitrary choice  $\theta:=+1$, i.e.,
\[
\theta
:= 
\begin{cases}
-1 & \text{ if } \; \tau_0 \in \{\alpha_i\}\,,
\\ 
+1 & \text{ if } \; \tau_0 \in \{\beta_i\}\,,
\\
+1 & \text{ if } \; \tau_0 \in \mathbb{M} \backslash \partial \supp \rho\,.
\end{cases}
\]
The minimum $ \tau_0 $ will be considered fixed in the following analysis. We parametrize $ z $ as follows in the neighborhood of $\tau_0 \in \mathbb{M}$:
\bels{}{
z \,=\, \tau_0 + \theta\2\omega + \cI\1\eta
\,,
}
where $ \eta \in (0,\delta_*]$ and $\omega \in [-\delta_*, \delta_*]$.
We will then prove the local law in the form
\begin{subequations}
\label{local law with error function}
\begin{align}
\Err(z)\,
&\prec\, \sqrt{\frac{\rho(z)}{N \eta}} + \frac{1}{N \eta}+\2\Bound(\omega,\eta)
\,,
\\
\absb{\avgb{\1\vect{w}, \vect{g}(z)-\vect{m}(z)}}
\,&\prec\;\, \Bound(\omega,\eta)
\,,
\end{align}
\end{subequations}
where the positive error function $\Bound: [-\delta_*, \delta_*] \times (0,\delta_*] \to (0,\infty)$ is given as the unique solution of an explicit cubic equation in \eqref{definition of error function} below.  

To define $\Bound$ we introduce explicit auxiliary functions $\widetilde \pi_1$, $\widetilde \pi_2$ and $\widetilde \rho$ that are comparable in size to the corresponding functions $\pi_1$, $\pi_2$ and $\rho$. The reason for using these auxiliary quantities for the definition of $\Bound$ instead of the original ones is twofold. Firstly, in this way $\Bound$ will be an explicit function  instead of one that is implicitly defined through the solution of the QVE. The function $\Bound$ is explicit in the sense that there is a formula for the solution of the cubic equation that defines it and the coefficients are given by the explicit functions $\widetilde \pi_1$, $\widetilde \pi_2$ and $\widetilde \rho$. Secondly, $\Bound$ will be monotonic of its second variable, $\eta$. This property will be used later. The definition of the three auxiliary functions will be different, depending on whether $\tau_0$ is in the boundary of the support of the density of states or not. 
Recall the definition \eqref{def of Delta_delta} of $ \Delta_\delta(\tau) $.

\begin{itemize}
\item {\bf Edge: } If $\tau_0 \in \partial \supp \rho$, i.e. $\tau_0$ is an edge of a gap of size $\Delta:=\Delta_{0}(\tau_0)$ in the support of the density of states or an extreme edge. Then we define the three explicit functions
\begin{subequations}
\begin{align}
\widetilde{\rho}(\omega,\eta)
\,&:=\, 
\begin{cases}
\displaystyle\frac{(|\omega|+\eta)^{1/2}}{(\Delta+|\omega|+\eta)^{1/6}}
\,,\qquad
&\omega \in \bigl[-\delta_*,\20\1\bigr]\,,
\\
\displaystyle\frac{\eta^{\left.\right.}}{(\Delta+\eta)^{1/6}(\omega+\eta)^{1/2}}
\,,\qquad
&\omega \in \bigl[\10\1,c_*\Delta\bigr]\,,
\\
\displaystyle
\frac{\eta^{\left.\right.}}{(\Delta+\eta)^{2/3}}
\,,\qquad
&
\displaystyle
\omega \in \biggl[c_*\Delta\1,\2\frac{\Delta}{2}\biggr]
\,.
\end{cases}
\\
\widetilde{\pi}_1(\omega,\eta)
\,&:=\, 
\begin{cases}
(|\omega|+\eta)^{1/2}(|\omega|+\eta+\Delta)^{1/6}
\,,\qquad
&\omega \in \bigl[-\delta_*,\20\1\bigr]\,,
\\
(\omega+\eta)^{1/2}(\Delta+\eta)^{1/6}
\,,\qquad
&\omega \in \bigl[\10\1,c_*\Delta\bigr]\,,
\\
(\Delta+\eta)^{2/3}
\,,\qquad
\displaystyle
&\omega \in \biggl[c_*\Delta\1,\2\frac{\Delta}{2}\biggr]
\end{cases}
\\
\label{pi2 tilde at an edge}
\widetilde{\pi}_2(\omega,\eta)
\,&:=\,
\begin{cases}
(|\omega|+\eta+\Delta)^{1/3}
\,,\qquad
&\omega \in \bigl[-\delta_*,\20\1\bigr]\,,
\\
(\Delta+\eta)^{1/3}
\,,\qquad
&\omega \in \bigl[\10\1,c_*\Delta\bigr]\,,
\\
(\Delta+\eta)^{1/3}
\,,\qquad
&
\displaystyle
\omega \in \biggl[c_*\Delta\1,\2\frac{\Delta}{2}\biggr]
\end{cases}
\end{align}
\end{subequations}
Here, \cob $ c_* \sim 1 $ \nc is the constant from Proposition \ref{prp:Behavior of the coefficients}.
\item {\bf Internal minimum: }
If $\tau_0 \in \mathbb{M} \backslash \partial \supp \rho$, then we define for $ \omega \in [-\delta_*,\delta_*]$ the three functions
\begin{subequations}
\begin{align}
\widetilde{\rho}(\omega,\eta)
\,&:=\
\rho(\tau_0)+ (|\omega|+\eta)^{1/3}
\,,
\\
\widetilde{\pi}_1(\omega,\eta)
\,&:=\,
\rho(\tau_0)^2+ (|\omega|+\eta)^{2/3}
\,,
\\
\widetilde{\pi}_2(\omega,\eta)
\,&:=\, 
\rho(\tau_0)+ (|\omega|+\eta)^{1/3}
\,,
\end{align}
\end{subequations}
\end{itemize}
By design (cf. Proposition \ref{prp:Behavior of the coefficients} and Theorem \ref{thr:Solution of the QVE}) these functions satisfy 
\bels{tilde has same size as no tilde}{
\rho(\tau_0 + \theta\2\omega + \cI\1\eta) 
\;\sim\; \widetilde{\rho}(\omega,\eta)
\,,\qquad\text{and}\qquad 
\abs{\1\pi_k(\tau_0 + \theta\2\omega + \cI\1\eta)}
\;\sim\; 
\widetilde{\pi}_k(\omega,\eta)
\,,
}
except in one special case where the second bound does not hold, namely when $k=2$, $\tau_0 \in \partial\supp \rho$ and $\omega \in [c_*\Delta,\Delta/2]$. In this case only the direction $|\pi_2| \lesssim \widetilde \pi_2$ is true (cf. \eqref{pi2 well inside a gap}).  

We fix a positive constant $\widetilde \eps\in (0, \gamma /16)$.
The value of the function $\Bound$ at $(\omega, \eta)$ is then defined to be the unique positive 
solution of the cubic equation
\bels{definition of error function}{
\Bound(\omega, \eta)^3 + \widetilde\pi_2(\omega, \eta) \2\Bound(\omega, \eta)^2 + \widetilde\pi_1(\omega, \eta) \2\Bound(\omega, \eta) \;=\, N^{8\1\widetilde\eps}\;\frac{\Bound(\omega, \eta)}{N\eta} +\frac{\widetilde{\rho}(\omega, \eta)}{N\eta}+\frac{1}{(N\eta)^2}\,,
}

With the choices \eqref{kappa bound} and \eqref{improved kappa bound} for $ \kappa = \kappa(z) $ we have
\bels{Bound kappa relationship}{
\Bound \,\leq\, N^{9\1\widetilde \eps}\min\setbb{\!\frac{1}{\sqrt{N\1\eta}}\,,\2 \frac{\kappa}{N\eta}\!}
\,,
}
for any $N \geq N_0$, where the threshold $N_0$ here depends on $\widetilde \eps$ in addition to $p$, $P$, $L$, $\ul{\mu}$ and $\gamma$. 
The inequality \eqref{Bound kappa relationship} is verified by plugging its right hand side into \eqref{definition of error function} in place of $\Bound$ and checking that on each regime the resulting expression on the right hand side of \eqref{definition of error function} is smaller than the resulting expression on the left hand side of \eqref{definition of error function}. The factor of $N^{9\1\widetilde \eps}$ in \eqref{Bound kappa relationship} can be absorbed in the stochastic domination in \eqref{local law with error function}.  Thus  \eqref{local law with error function} becomes \eqref{pointwise local law} and \eqref{averaged local law} of Theorem~\ref{thr:Local law}. 

Before we start the proof of the local law \eqref{local law with error function}, let us motivate the definition of $\Bound$.
As a consequence of Lemma~\ref{lmm:Rough bound} the indicator function equals one a.w.o.p. in the statement of Lemma~\ref{lmm:Bound on perturbation}. Thus, uniformly in the $\delta_*$-neighborhood of $\tau_0$ we have 
\bels{d and off diagonal error bound}{
\norm{\vect{d}}_{\infty} +\Err_{\mathrm{o}} \,\prec\, \sqrt{\frac{\1\rho+\abs{\avg{\1\brm{g}-\brm{m}\1}}}{N\eta}}
+
\frac{1}{\!\sqrt{\msp{-1}N\1}}
\,.
}
Here we used $\Im \avg{\1\brm{g}} \lesssim \rho + \abs{\avg{\1\vect{g}-\vect{m}\1}}$.
Since at the end the local law implies $\abs{\avg{\1\brm{g}-\brm{m}\1}}\prec \Bound$, heuristically we may replace $\abs{\avg{\1\brm{g}-\brm{m}\1}}$ in \eqref{d and off diagonal error bound} by $\Bound$. In this case, from the fluctuation averaging, Theorem \ref{thr:Fluctuation Averaging}, we would be able to conclude that for any  deterministic vector $\vect{w}$ with bounded entries,
\bels{heuristic d bound}{
\norm{\vect{d}}_{\infty}^2 + |\la \1\vect{w},\vect{d}\1\ra|
\;\prec\; 
\frac{\Bound}{N\eta}+\frac{\rho}{N\eta}+\frac{1}{(N\eta)^2}\,.
}
Up to the technical factor of $N^{8\1\eps}$ the right hand side coincides with the right hand side of the cubic equation defining $\Bound$. On the other hand, 
the right hand side of the cubic equation \eqref{cubic for Theta} for the quantity $\Theta$ from Theorem \ref{thr:Stability around support} is of the same form as the left hand side of \eqref{heuristic d bound}. Therefore, we infer
\bels{bound on cubic in Theta}{
\abs{\2\Theta^3+\pi_2\1\Theta^2+\pi_1\1\Theta\2}\,\prec\,\frac{\Bound}{N\eta}+\frac{\rho}{N\eta}+\frac{1}{(N\eta)^2}\,.
}
We will argue that on appropriately chosen domains out of the three summands in the cubic expression in $\Theta$ always one is the biggest by far. 
Therefore, the error function $\Bound$, defined by \eqref{definition of error function}, is essentially the best bound on $\Theta$ that one may hope to deduce from \eqref{bound on cubic in Theta}. Indeed, since $\Theta$ is by definition an average of $\vect{g}-\vect{m}$, we expect $\Theta \prec \Bound$.

We will now prove \eqref{local law with error function}. To this end we gradually improve the bound on $\Theta$. Fix some $\eps\in (0,\widetilde\eps)$.
The sequence of deterministic bounds on this quantity is defined as
\bels{definition of Phik}{
\Phi_0\,:=\, 1\,, \qquad \Phi_{k+1}\,:=\, \max\setb{\2N^{-\1\eps}\1\Phi_k\2, \, N^{\19\1 \eps}\Bound\2}
\,.
}
From here on until the end of this section the threshold function $N_0$ from the definition of the stochastic domination (cf. Definition \ref{def:Stochastic domination}) as well as the definition of 'a.w.o.p.' (cf. Definition \ref{def:Overwhelming probability}) may depend on $\eps$ in addition to $p$, $P$, $L$, $\ul{\mu}$ and $\gamma$. At the end of the proof we will remove this dependence. 
The following lemma is essential for doing one step in the upcoming iteration.

\NLemma{Improving bound through cubic}{Suppose that for all $z\in \tau_0 + [-\delta_*,\delta_*]+\cI\1[N^{\gamma-1},\delta_*]$ and some $k \in \N$ the quantity $\Theta(z)$ from \eqref{Def of Theta} fulfills
\bels{cubic inequality for Omega}{
\absb{\,\Theta(z)^3 +\pi_2(z)\2\Theta(z)^2+\pi_1(z)\2\Theta(z)\,}
\,\prec\, 
\frac{\rho(z)+\Phi_k(\omega,\eta)}{N\eta}+\frac{1}{(N\eta)^2}\,.
}
Then $\Theta(z)\prec\2 \Phi_{k+1}(\omega,\eta)$.
}

We will postpone the proof of this lemma until the end of this section. First we show how to use this result to prove the local law (Theorem \ref{thr:Local law}). Fix an integer $k\geq 0$ and assume that $\Theta+\abs{\avg{\1\brm{g}-\brm{m}\1}}\prec \Phi_k$ is already proven. 
For $k=0$ this follows from the rough bound on $\Err$ in Lemma \ref{lmm:Rough bound}, $\Err \prec 1=\Phi_0$.
As an induction step we show that $\Theta+\abs{\avg{\1\brm{g}-\brm{m}\1}}\prec \Phi_{k+1} $.

From \eqref{d and off diagonal error bound} we see that
\bels{d and Lambdao bound}{
\norm{\vect{d}}_{\infty} +\2 \Err_{\mathrm{o}} \,\prec\, \sqrt{\frac{\rho+\Phi_k}{N\eta}}+\frac{1}{N\eta}\,.
}
The right hand side is a deterministic bound on the off-diagonal error $\Err_{\mathrm{o}}$.
Therefore the fluctuation averaging (Theorem \ref{thr:Fluctuation Averaging}) is applicable to  $ \avg{\2\brm{t}^{(1)}\!,\brm{d}\1} $ and $ \avg{\2\brm{t}^{(2)}\!,\brm{d}\1} $ on right hand side of the cubic equation \eqref{cubic for Theta} 
\[
\absb{\avg{\2\brm{t}^{(j)}\!,\brm{d}\1}}
\,\prec\, 
\biggl(\sqrt{\frac{\rho+\Phi_k}{N\eta}}+\frac{1}{N\eta}\biggr)^2
\,,
\]
where $ N^{-1} $ from \eqref{fluctuation averaging result} has been neglected since $ \rho \gtrsim \eta $. \nc
In this way we see that the hypothesis \eqref{cubic inequality for Omega} of Lemma \ref{lmm:Improving bound through cubic} is satisfied. 
Using the lemma the bound on $\Theta$ is improved to 
\bels{}{
\Theta(z)\,\prec\, \Phi_{k+1}(\omega,\eta)\,.
}

In order to improve the bound on $\abs{ \avg{\1 \brm{g} - \brm{m} \1} }$ as well, we 
use the bound \eqref{bound on averages against bounded vectors} from Theorem \ref{thr:Stability around support} for averages of $\vect{g}-\vect{m}$ against bounded vectors. Since by Lemma \ref{lmm:Rough bound} the deviation function $\Lambda$ is bounded by a small constant, the indicator function in \eqref{bound on averages against bounded vectors} is a.w.o.p. non-zero.
Choosing $\vect{w}=(1, \dots,1)$, we find that
\bels{}{
\abs{\avg{\1\vect{g}-\vect{m}\1}}
\,\lesssim\, 
\Theta+ \norm{\vect{d}}_\infty^2
+
|\la\1  \widetilde{\vect{ w}}, \vect{d}\1\ra|\,,\qquad \text{a.w.o.p.}\,,
}
where $\widetilde{\vect{w}}=\vect{T}\vect{w}$ is a bounded, $\norm{\widetilde{\vect{w}}}_\infty\lesssim 1$, deterministic vector. Together with the bound \eqref{d and Lambdao bound} we apply the fluctuation averaging (Theorem \ref{thr:Fluctuation Averaging}) again, 
\bels{}{
\abs{\avg{\1\vect{g}-\vect{m}\1}}
\;\prec\; 
\Phi_{k+1}
+ \frac{\rho+\Phi_k}{N\eta}+\frac{1}{(N\eta)^2}
\;\lesssim\; 
N^{-\eps}\Phi_k+\Phi_{k+1}
\;\lesssim\; 
\Phi_{k+1}
\,.
}
This concludes one step in the iteration, i.e., we have shown $\Theta+\abs{\avg{\1\brm{g}-\brm{m}\1}}\prec\2 \Phi_{k+1}$.

We repeat this step finitely many times and each time improve $\Phi_k$ by a factor of $N^{-\eps}$ until it reaches its target value $N^{9 \eps }\Bound$ and is not improved anymore. 
Note that all constants in our estimates, explicit and hidden, depend only on the model parameters and $ \eps $.  In particular, the number of steps needed is uniform in $ (\omega,\eta)$.
At that stage we have 
\[
\Theta+\abs{\avg{\1\vect{g}-\vect{m}\1}}
\,\prec_{\eps}\, 
N^{\19\1\eps}\Bound
\,,
\]
where the subindex $\eps$ indicates that the threshold $N_0$ from the  stochastic domination may depend on $\eps$. But since $\eps > 0 $ was arbitrary, we infer (cf. (i) of Lemma~\ref{lmm:Basic facts about stochastic domination}) that  $ \Theta+\abs{\avg{\1\brm{g}-\brm{m}\1}}\prec\Bound$, where now and until the start of the proof of Lemma \ref{lmm:Improving bound through cubic} below the stochastic domination is $\eps$-independent. By \eqref{d and off diagonal error bound}  we conclude
\bels{stochastic domination bound on d and off-diagonal error}{
\norm{\vect{d}}_{\infty} + \Err_{\mathrm{o}}\,\prec\, \sqrt{\frac{\rho}{N\1\eta}}+\frac{1}{N\1\eta}+ \Bound\,.
}

For the bound on the diagonal contribution, $\Err_{\mathrm{d}}$, we use \eqref{Bound of Err in terms of Theta} to get
\[
\Err_{\mathrm{d}}\,\lesssim\, \Theta + \norm{\vect{d}}_{\infty} \,\prec\, 
\sqrt{\frac{\rho}{N\1\eta}}+\frac{1}{N\1\eta} + \Bound\,.
\]
Finally, with the help of \eqref{bound on averages against bounded vectors}, \eqref{stochastic domination bound on d and off-diagonal error} and the fluctuation averaging, we prove the bound on averages of $\vect{g}-\vect{m}$ against any bounded, $\norm{\vect{w}}_\infty\leq 1$, deterministic vector,
\[
\abs{\avg{\2\vect{w},\vect{g}-\vect{m}\1}}
\,\prec\,  
\frac{\rho}{N\1\eta}+\frac{1}{(N\1\eta)^2}+\Theta
\,\prec\,\frac{\rho}{N\1\eta}+\frac{1}{(N\1\eta)^2}+\Bound\,.
\]
This finishes the proof of Theorem \ref{thr:Local law} apart from the proof of Lemma \ref{lmm:Improving bound through cubic} which we will tackle now.

\begin{Proof}[Proof of Lemma \ref{lmm:Improving bound through cubic}]
The spectral parameter $z=\tau_0 + \theta\2\omega + \cI\1\eta$ lies inside 
the $\delta_*$-neighborhood of $\tau_0$. We fix $\omega \in [-\delta_*,\delta_*]$ and show that the claim holds for any choice of $\eta \in [N^{\gamma-1},\delta_*]$.
We split the interval of possible values of $\eta$ into two or three regimes, depending on the case we are treating. 
\begin{itemize}
\item {\bf Edge: } If $\tau_0 \in \partial \supp \rho$ is an edge of a gap of size $\Delta:=\Delta_{0}(\tau_0)$, then we define 
\begin{equation*}
\begin{split}
I_1(\omega)\,&:=\, \setbb{\eta\in [N^{\gamma-1},\delta_*]:\, \frac{(|\omega|+\eta)^{1/2}}{(|\omega|+\eta+\Delta)^{1/6}}\2\geq\2 N^{-5\1\eps}\Phi_k(\omega,\eta) }
\,,
\\
I_2(\omega)\,&:=\, \setbb{
\eta\in [N^{\gamma-1},\delta_*]:\, N^{5\1\eps}\frac{(|\omega|+\eta)^{1/2}}{(|\omega|+\eta+\Delta)^{1/6}}\2\leq\2 \Phi_k(\omega,\eta)\2\leq\2N^{2\1\eps}(|\omega|+\eta+\Delta)^{1/3} 
}
\,,
\\
I_3(\omega)\,&:=\, \setbb{
\eta\in [N^{\gamma-1},\delta_*]:\, (|\omega|+\eta+\Delta)^{1/3}\2\leq\2 N^{-2\1\eps}\Phi_k(\omega,\eta)
}
\,.
\end{split}
\end{equation*}
If any of the two regimes $I_l(\omega)$ with $l=2,3 $ consists of a single point only, then we set $I_l(\omega):=\emptyset$.
\item {\bf Internal minimum: }
If $\tau_0 \in \mathbb{M} \backslash \partial \supp \rho$, then we set $I_2(\omega):=\emptyset$ and define
\begin{equation*}
\begin{split}
I_1(\omega)\,&:=\,\setB{ \eta\in [N^{\gamma-1},\delta_*]: \, \rho(\tau_0)+(|\omega|+\eta)^{1/3}\2\geq\2 N^{-2\1\eps}\Phi_k(\omega,\eta)}
\,,
\\
I_3(\omega)\,&:=\,\setB{ \eta\in [N^{\gamma-1},\delta_*]: \, \rho(\tau_0)+(|\omega|+\eta)^{1/3}\2\leq\2 N^{-2\1\eps}\Phi_k(\omega,\eta)}
\,.
\end{split}
\end{equation*}
If $I_3(\omega)$ consists of a single point only, then we set $I_3(\omega):=\emptyset$.
\end{itemize}
\begin{figure}
	\centering
	\includegraphics[width=0.85\textwidth]{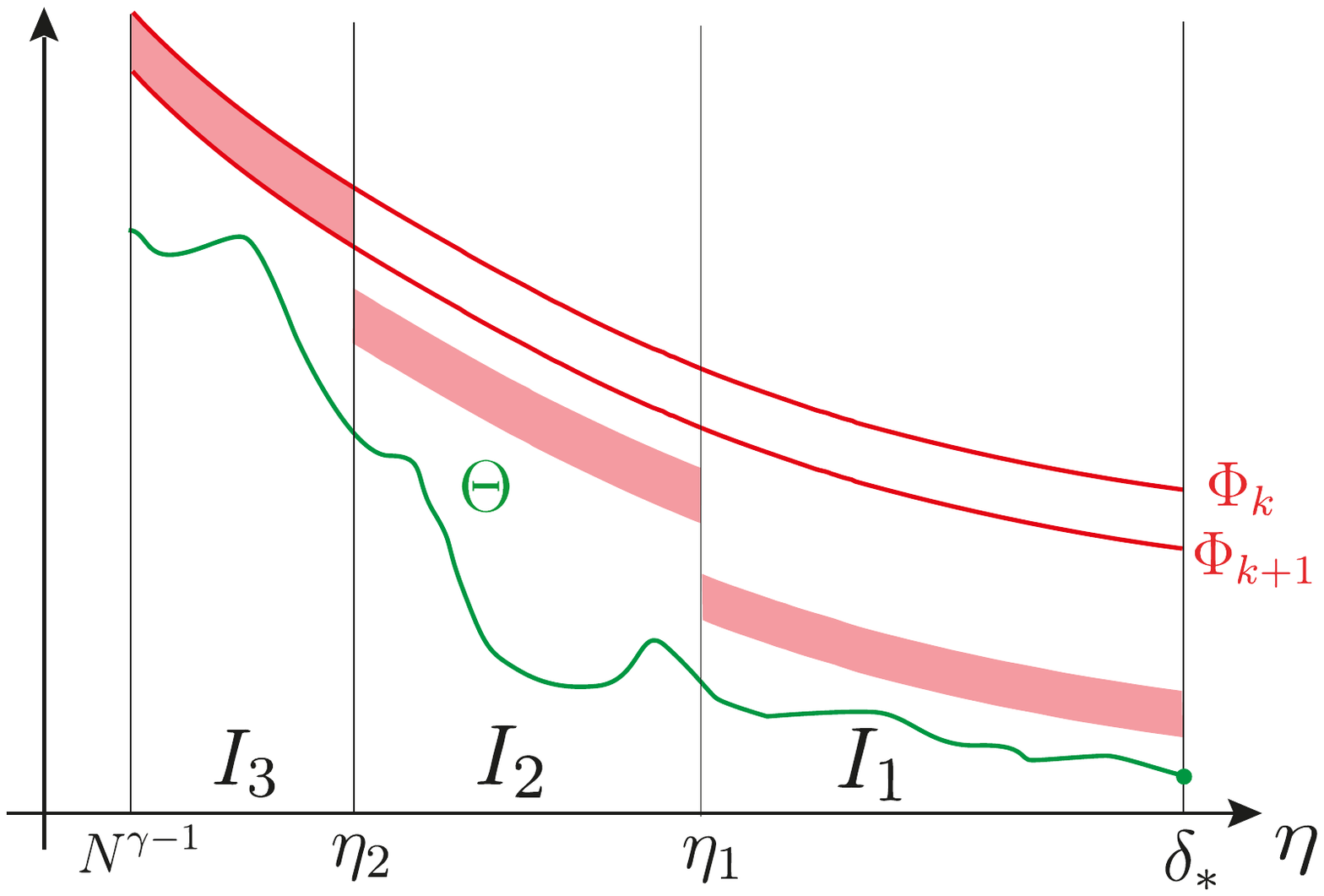}
	\caption{
	The shaded area is forbidden for $ \Theta $. Since the continuous function $ \Theta $ lies below this region at $ \eta = \delta_\ast $ it stays below it for any $ \eta \ge N^{\gamma-1} $, hence proving $ \Theta \leq \Phi_{k+1} $.
	}
	\label{GapArgumentFigure}
\end{figure}
In the cubic equation \eqref{definition of error function}, used to define the error function $\Bound$, the coefficients $\widetilde \pi_1$ and $\widetilde \pi_2$ on the left hand side are monotonously increasing functions of $\eta$. The linear and the constant coefficient of $\Bound$ on the right hand side are monotonously decreasing in $\eta$. Thus, $\Bound$ itself is a monotonously decreasing function of $\eta$. From this fact and the definition of the regimes $I_1$, $I_2$ and $I_3$ we see that $I_1=[\eta_1,\delta_*]$, $I_2=[\eta_2,\eta_1]$ and $I_3=[N^{\gamma-1},\eta_2]$ for some $\eta_1, \eta_2 \in [N^{\gamma-1},\delta_*]$. Here, we interpret $I_2=\emptyset$ if $\eta_1\leq \eta_2$ and $I_3=\emptyset$ if $\eta_2 \leq N^{\gamma-1}$.

Now we define a $z$-dependent indicator function
\bels{}{
\chi(\omega,\eta)\,:=\, 
\begin{cases}
\Ind\Bigl(\2N^{-7\eps}\Phi_k(\omega,\eta)\2\leq\2\Theta(\tau_0 + \theta\2\omega + \cI\1\eta)\2\leq\2 N^{-6 \eps}\Phi_k(\omega,\eta)
\2\Bigr)
& \text{ if } \;\eta \in I_{1}(\omega)
\vspace{0.1cm}
\\
\Ind\Bigl(\2
N^{-4\eps}\Phi_k(\omega,\eta)\2\leq\2\Theta(\tau_0 + \theta\2\omega + \cI\1\eta)\2\leq\2 N^{-3\eps}\Phi_k(\omega,\eta)
\2\Bigr)
& \text{ if } \;\eta \in I_{2}(\omega)
\vspace{0.1cm}
\\
\Ind\Bigl(\2
N^{-\eps}\Phi_k(\omega,\eta)\2\leq\2\Theta(\tau_0 + \theta\2\omega + \cI\1\eta)\2\leq\2 \Phi_k(\omega,\eta)
\2\Bigr)
& \text{ if }\; \eta \in I_{3}(\omega)
\end{cases}\,.
}
This function fixes the values of $\Theta$ to a small interval just below the deterministic control parameter $\Phi_k$. We will prove that $\Theta$ cannot take these values, i.e. $\chi=0$ a.w.o.p.. Figure \ref{GapArgumentFigure} illustrates this argument. Compared to Figure 6.1 in \cite{EKYY} we see that instead of two there are now three domains, $I_1(\omega)$, $I_2(\omega)$ and $I_3(\omega)$, to be distinguished. The reason for this extra complication is that \eqref{cubic for Theta} is cubic in $\Theta$, compared to the quadratic equation for $[v]$ that appeared in the proof of Lemma 6.2 in \cite{EKYY}.
To see that $\chi=0$, first note that the choice of the domains, $I_l$, ensures that there is always one summand on the left hand side of the cubic equation 
\eqref{cubic for Theta} for $\Theta$ which dominates the two others by a factor $N^\eps$, whenever $\chi$ does not vanish. In fact, by construction we have: 
\\[0.3cm]
{\sc Claim:}\emph{ 
 The random functions $\Theta$ and $\chi$ satisfy a.w.o.p.}
\bels{compare with and without absolute value}{
\Bigl(\Theta(z)^3 +\2\widetilde \pi_2(\omega,\eta)\2\Theta(z)^2+\2\widetilde\pi_1(\omega,\eta)\2\Theta(z) \Bigr)\,\chi(\omega, \eta)\quad&
\\
\lesssim\;
\absB{\Theta(z)^3 +\pi_2(z)\2\Theta(z)^2+\pi_1(z)\2\Theta(z)}&
\;.
}

We will verify this fact at the end of the proof of this lemma. Now we will simply use it.
First we combine the assumption \eqref{cubic inequality for Omega} of the lemma and \eqref{compare with and without absolute value} to obtain
\[
N^{-\eps}(\2\Theta^3 + \wti{\pi}_2\Theta^2 + \wti{\pi}_1\Theta\2)\;\chi 
\,\leq\,   
\frac{\2\wti{\rho}+\Phi_k}{N\eta} + \frac{1}{(N\eta)^2}
\qquad
\text{a.w.o.p.}
\,.
\]
Here we also gave up a factor of $N^\eps$ to get an inequality instead of the stochastic domination, and replaced $ \rho $ by the comparable quantity $ \wti{\rho} $.
By the definition of the indicator function $\chi$ we have $\Theta \chi \geq N^{-7\1\eps} \Phi_k $. 
Using this to bound the left hand side, and that $ \eps \leq \widetilde \eps $, we obtain
\[
\big(\mathcal{R}^3 +\widetilde \pi_2\2\mathcal{R}^2+\widetilde\pi_1\2\mathcal{R} \big)\2\chi\,\leq\,N^{8\1\widetilde \eps}\frac{\mathcal{R}}{N\eta}+
\frac{\widetilde\rho}{N\eta}+\frac{1}{(N\eta)^2}
\,,\qquad 
\text{a.w.o.p.}\,, 
\qquad\mathcal{R}\2:=\2 N^{-8\1\eps}\Phi_k
\,.
\]
Comparing this with the defining equation \eqref{definition of error function} for $ \Bound $ we conclude that a.w.o.p. $ N^{-8\1\eps}\Phi_k\2\chi\leq \Bound $. 

On the other hand, by the definition of $\Phi_k$ in \eqref{definition of Phik} we know that $\Phi_k> N^{8\1\eps}\Bound$. These two inequalities yield
\bels{chi vanishes}{
\chi(\omega,\eta)\,=\, 0 \,,\qquad  \eta \in [N^{\gamma-1},\delta_*]\,,\qquad \text{a.w.o.p.}\,.
}

Now we successively, for $l=1,2,3$, apply Lemma \ref{lmm:Bound propagation} on the connected domains $\tau_0 + \theta \1\omega + \cI\2I_l(\omega)$ with the choices 
$\varphi:=\Theta$ and
\[
\Phi(\tau_0 + \theta\1\omega +\cI\1\eta)\,:=\, 
\begin{cases}
N^{-6\1\eps}\Phi_k(\omega,\eta) & \text{ if }\; l=1\,,
\\
N^{-3\1\eps}\Phi_k(\omega,\eta) & \text{ if }\; l=2\,,
\\
\Phi_k(\omega,\eta) & \text{ if }\; l=3\,,
\end{cases}
\quad 
z_0\,:=\, 
\begin{cases}
\tau_0 + \theta\1\omega +\cI\1\delta_* & \text{ if }\; l=1\,,
\\
\tau_0 + \theta\1\omega +\cI\1\eta_1 & \text{ if }\; l=2\,,
\\
\tau_0 + \theta\1\omega +\cI\1\eta_2 & \text{ if }\; l=3\,,
\end{cases}
\]
where as explained after the definition of $I_1$, $I_2$ and $I_3$ above we have $I_1=[\eta_1,\delta_*]$, $I_2=[\eta_2,\eta_1]$ and $I_3=[N^{\gamma-1},\eta_2]$.
The condition \eqref{weak continuity condition} of the lemma is satisfied by the definition of $\Theta$ in \eqref{Def of Theta}, the H\"older-continuity of the solution of the QVE, the weak Lipschitz-continuity of $\vect{g}$ with Lipschitz-constant $N^{\12}$ and the H\"older-continuity of $\vect{s}$ from \eqref{sigma and s continuity}. 
The gap condition,
\eqref{gap condition}, holds because of \eqref{chi vanishes} and the definition of $\chi$ and $\Phi$ for an appropriate choice of the exponent $D_3$.

The condition, $\varphi(z_0) \leq \Phi(z_0)$ a.w.o.p., necessary for the application of Lemma \ref{lmm:Bound propagation} on the first domain, $\tau_0 + \theta \1\omega + \cI\2I_1(\omega)$, is obtained form Proposition \ref{prp:Local law away from local minima}. With Lemma \ref{lmm:Bound propagation} we propagate the bound to all $z \in \tau_0 + \theta \1\omega + \cI\2I_1(\omega)$.
Now we apply Lemma \ref{lmm:Bound propagation} on the second domain $\tau_0 + \theta \1\omega + \cI\2I_2(\omega)$, provided $I_2(\omega)$ is not empty. The bound \eqref{initial value condition} for the new $z_0=\tau_0 + \theta \1\omega + \cI\1\eta_1$ is obtained from the previous step. 
 Finally, we apply Lemma \ref{lmm:Bound propagation} to $\tau_0 + \theta \1\omega + \cI\2I_3(\omega)$, in case it is not empty, with the new choice $z_0=\tau_0 + \theta \1\omega + \cI\1\eta_2$. 
 Altogether,  we applied the lemma at most three times. Through this procedure we prove that a.w.o.p. $\Theta(z) \leq \Phi(z)$ for all $z \in  \tau_0 + \theta \1\omega + \cI\2[N^{\gamma-1},\delta_*]$. On the third domain, $\tau_0 + \theta \1\omega + \cI\2I_3(\omega)$, we use that a.w.o.p. $\chi=0$ (cf. \eqref{chi vanishes}) and thus a.w.o.p. $\Theta(z) \leq N^{-\eps}\Phi_k$. Altogether we showed that in the $\delta_*$-neighborhood of $\tau_0$,
\[
\text{a.w.o.p.}\qquad \Theta(z)\,\leq\, N^{-\eps}\Phi_k\,\leq\, \Phi_{k+1}\,.
\]
This finishes the proof of Lemma \ref{lmm:Improving bound through cubic} up to verifying the claim \eqref{compare with and without absolute value}. 
\\[0.3cm]
{\sc Proof of the claim: }
For the proof of \eqref{compare with and without absolute value} one verifies case by case that on $I_1$ the term $\widetilde \pi_1 \Theta \sim |\pi_1| \Theta$ is bigger than the two other terms, $\widetilde \pi_2 \Theta^2$ and $\Theta^3$ by a factor of $N^\eps$. If $I_3$ is not empty then the term $\Theta^3$ is the biggest. If $I_2$ is not empty, then $|\pi_2| \sim \widetilde \pi_2$ and $\widetilde \pi_2 \Theta^2$ is the biggest term by a factor of $N^\eps$. More specifically, when $ \eta \in I_j $ and $ \chi = \chi(\omega,\eta) = 1 $ we show   
\[
\absb{\Theta^3 + \pi_2\Theta^2 +\pi_1\Theta}
\,\sim\,  
\abs{\pi_j}\Theta^j 
\,\sim\, 
\widetilde\pi_j\Theta^j 
\,\sim\;
\Theta^3 + \widetilde
\pi_2\Theta^2 +\widetilde\pi_1\Theta
\,,
\]
where $ \pi_3 = \wti{\pi}_3 := 1 $. 
As an example we demonstrate these relations in a few cases:
\begin{itemize}
\item {\bf Well inside a gap: } If $\tau_0 \in \partial \supp \rho$ and $\omega \in [c_*\Delta, \Delta/2]$ then $I_2(\omega)=\emptyset$. We now check that on $I_1(\omega)$ the linear term in $\Theta$ is the biggest while on $I_3(\omega)$ the cubic term dominates. 
First, let $\eta \in I_1(\omega)$. Then the following chain of inequalities hold, 
\[
\widetilde{\pi}_1 \Theta \,\sim\, |\pi_1| \Theta\,\sim\, (\Delta+\eta)^{2/3} \Theta
\,\gtrsim\, N^{-5\1\eps} (\Delta+\eta)^{1/3} \2\Phi_k \2\Theta
\,\sim\, N^{-5\1\eps}\2\widetilde\pi_2 \2\Phi_k \2\Theta
\,\gtrsim\, N^{-10\1\eps} \2\Phi_k^2 \2\Theta\,.
\]
Here, we used \eqref{tilde has same size as no tilde}, \eqref{pi1 well inside a gap}, the definition of $I_1(\omega)$ and \eqref{pi2 tilde at an edge} in the form $\widetilde \pi_2 \sim (\Delta +\eta)^{1/3}$.
Now we can use $ \chi $ to replace $ \Phi_k $ by $ \Theta $. By definition of $\chi$ and since $\widetilde \pi_k \gtrsim |\pi_k|$ for $k=1,2$ we also get
\[
N^{-5\1\eps}\2\widetilde\pi_2 \2\Phi_k \2\Theta\2\chi\,\geq\, N^{\eps}
\2\widetilde\pi_2 \2\Theta^2\2\chi\,\gtrsim\, N^{\eps}
|\pi_2|\Theta^2\2\chi\,, \qquad
N^{-10\1\eps} \2\Phi_k^2 \2\Theta\2\chi\,\geq\,N^{2\1\eps}\2\Theta^3\2\chi\,.
\]
We conclude that on $I_1(\omega)$ the linear term in $\Theta$ dominates the others,
\[
\widetilde{\pi}_1 \Theta\2\chi \,\gtrsim\,N^{\eps}(\Theta^3+ \widetilde{\pi}_2 \Theta^2)\2\chi\,.
\]
Suppose now that $\eta \in I_3(\omega)$. In this case, using the choice of the indicator function $\chi$,
\[
\Theta^3 \2\chi\,\geq\,   N^{-\eps}\Phi_k\2\Theta^2 \2\chi
\,\geq\,   N^{-2\1\eps}\Phi_k^2\2\Theta \2\chi\,.
\]
By definition of $I_3(\omega)$ and \eqref{pi2 tilde at an edge} we find that
\[
N^{-\eps}\Phi_k\2\Theta^2 \,\gtrsim\, N^{\eps}(\Delta+\eta)^{1/3}\2\Theta^2
\,\sim\,  N^{\eps}\widetilde \pi_2\2\Theta^2
\,,\quad
 N^{-2\1\eps}\Phi_k^2\2\Theta \,\gtrsim\,N^{2\1\eps}(\Delta+\eta)^{2/3}\1\Theta 
 \,\sim\, N^{2\1\eps}\widetilde \pi_1 \1\Theta 
\,.
\]
Altogether we find that the cubic term dominates the two others,
\[
\Theta^3\2\chi \,\gtrsim\, N^\eps(\widetilde \pi_2 \1\Theta^2+\widetilde \pi_1 \1\Theta)\1\chi\,.
\]
\item {\bf Inside a gap close to an edge on $I_2$: } If $\tau_0 \in \partial \supp \rho$, $\omega \in [0,c_*\Delta]$ and $\eta \in I_2(\omega)$, then we will show the quadratic term in $\Theta$ dominates the two other terms. We have
\[
|\pi_2| \Theta^2\,\sim\,\widetilde{\pi}_2 \2\Theta^2\,\sim\, (\Delta+\eta)^{1/3}\2 \Theta^2\,\gtrsim\, N^{-2\1\eps}\Phi_k\2\Theta^2,
\]
where in the inequality we used the definition of $I_2(\omega)$. The choice of $\chi$ guarantees that $\Phi_k \2\chi \geq N^{3\1\eps}\1\Theta\2\chi$. Thus, the quadratic term is larger than the cubic term by a factor of $N^\eps$.  On the other hand
\bea{
(\Delta+\eta)^{1/3}\2 \Theta^2\2\chi\,&\gtrsim\,N^{-4\1\eps}(\Delta+\eta)^{1/3}\Phi_k\2 \Theta\,\gtrsim\, N^\eps(\omega+\eta)^{1/2}(\Delta+\eta)^{1/6}\2\Theta
\,\sim\, N^\eps\1\widetilde \pi_1\1\Theta
\\
&\sim\, N^\eps|\pi_1|\Theta
\,.
}
Here, in the first inequality we used the indicator function $\chi$ and in the second inequality the definition of $I_2(\omega)$. Altogether, we arrive at
\[
\widetilde \pi_2 \2\Theta^2\2\chi\,\gtrsim\, N^\eps(\Theta^3+ \widetilde\pi_1 \2\Theta)\2\chi\,.
\]
\item {\bf Internal minimum on $I_1$:} If $\tau_0 \in \mathbb{M} \backslash\partial \supp \rho$ and $\eta \in I_1(\omega)$, then the linear term is the biggest,
\[
|\pi_1|\Theta \,\sim\, \widetilde\pi_1\2\Theta\,\sim\, \big(\rho(\tau_0)^2 + (|\omega|+\eta)^{2/3}\big)\2\Theta\,\gtrsim \, N^{-2\1\eps}\big(\rho(\tau_0) + (|\omega|+\eta)^{1/3}\big)\Phi_k\2\Theta\,.
\] 
Here, we used \eqref{tilde has same size as no tilde} and the definitions of $\widetilde \pi_1$ and $I_1(\omega)$, respectively.
Since $\Phi_k \2\chi \geq N^{6\1\eps}\Theta \2\chi$ and by the definition of $\widetilde\pi_2 $ this shows that the linear term is larger than the quadratic term by a factor of $N^{4\1\eps}$. In order to compare the linear with the cubic term we estimate further. By definition of $I_1(\omega)$,
\[
N^{-2\1\eps}\bigl(\,\rho(\tau_0) + (|\omega|+\eta)^{1/3}\,\bigr)\2\Phi_k\2\Theta\,\geq\,
N^{-4\1\eps}\2\Phi_k^2\1\Theta\,.
\]
Again we use the lower bound on $\Phi_k \2\chi$ and get
\[
N^{-4\1\eps}\2\Phi_k^2\1\Theta\2\chi \geq\, N^{8\1\eps}\Theta^3\2\chi\,.
\]
Thus we showed that on the domain $I_1(\omega)$ 
\[
\widetilde \pi_1 \2\Theta \, \chi \,\gtrsim\, N^\eps\, (\1\Theta^3+\widetilde \pi_2 \2\Theta^2\2)\2\chi\,.
\]
\end{itemize}
The other cases are proven similarly. This completes the proof of \eqref{compare with and without absolute value}.
\end{Proof}

\section{Rigidity and delocalization of eigenvectors}

\subsection{Proof of Corollary \ref{crl:Convergence of cumulative eigenvalue distribution}}
Here we explain how the local law, Theorem \ref{thr:Local law}, is used to estimate the difference between the cumulative density of states and the eigenvalue distribution function of the random matrix $\vect{H}$. The following auxiliary result shows that the difference between two probability measures can be estimated in terms of the difference of their respective Stieltjes transforms. For completeness the proof is given in the appendix. It uses a Cauchy-integral formula that was also applied in the construction of the Helffer-Sj\"ostrand functional calculus (cf. \cite{DFunc}) and it appeared in different variants in \cite{EYYber}, \cite{ERSY} and \cite{EYY}.

\NLemma{Bounding measures by Stieltjes transforms}{There is a universal constant $C>0$, such that for any two probability measures, $\nu_1$ and $\nu_2$, on the real line and any three numbers $\eta_1, \eta_2,\eps \in (0,1]$ with $\eps\geq \max\{\eta_1,\eta_2\}$, the difference between the two measures evaluated on the interval $[\tau_1,\tau_2]\subseteq \R$, with $\tau_1< \tau_2$, satisfies
\bels{Bounding measures by Stieltjes transforms}{
\absb{\1\nu_1([\tau_1,\tau_2])\,-\,\nu_2([\tau_1,\tau_2])\1}
\;\leq\;
C\,\Bigl(\,
\nu_1([\tau_1-\eta_1,\tau_1])
\,+\,
\nu_1([\tau_2,\tau_2+\eta_2])
\,+\, J_1 +J_2+J_3
\Bigr)
\,.
}
Here, the three contributions to the error, $J_1$, $J_2$ and $J_3$, are defined as
\bels{definition of Js}{
J_1\,&:=\, 
\int_{\tau_1-\eta_1}^{\tau_1} \msp{-8}\dif \omega\,\biggl( 
\Im\2 m_{\nu_1}(\omega+\cI\1\eta_1)+|m_{\nu_1-\nu_2}(\omega+\cI\1\eta_1)|+\frac{1}{\eta_1} \int_{\eta_1}^{2\1\eps} \dif \eta\1|  m_{\nu_1-\nu_2}(\omega+\cI\1\eta)|
\biggr)
\,,
\\
J_2\,&:=\, 
\int_{\tau_2}^{\tau_2+\eta_2} \msp{-8}\dif \omega\,\biggl( 
\Im\2 m_{\nu_1}(\omega+\cI\1\eta_2)+|m_{\nu_1-\nu_2}(\omega+\cI\1\eta_2)|+\frac{1}{\eta_2} \int_{\eta_2}^{2\1\eps} \dif \eta\1|  m_{\nu_1-\nu_2}(\omega+\cI\1\eta)|
\biggr)
\,,
\\
J_3\,&:=\, 
\frac{1}{\eps}\int_{\tau_1-\eta_1}^{\tau_2+\eta_2}\msp{-8}\dif \omega \int_{\eps}^{2\1\eps}\dif \eta\1 | m_{\nu_1-\nu_2}(\omega+\cI\1\eta)|\,,
}
where $m_\nu$ denotes the Stieltjes transform of $\nu$ for any signed measure $\nu$. 
}

We will now apply this lemma to prove Corollary \ref{crl:Convergence of cumulative eigenvalue distribution} with the choices of the measures
\bels{choices of nus}{
\nu_1(\dif \omega)\,:=\, \rho(\omega)\dif \omega\,,\qquad\text{and}\qquad \nu_2(\dif \omega )\,:=\, \frac{1}{N}\sum_{i=1}^N\delta_{\lambda_i}(\dif \omega)
\,.
}

As a first step we show that a.w.o.p. there are no eigenvalues with an absolute value larger or equal than $10$, i.e.,
\bels{No evs beyond 10}{
\#\sett{\,i:\2 |\lambda_i|\geq 10\2}\,=\,0\qquad \text{a.w.o.p.}\,.
}
We focus on the eigenvalues $\lambda_i \geq 10$. The ones with  $\lambda_i \leq -10$ are treated in the same way. 
We will show first that there are no eigenvalues in a small interval around $\tau$ with $\tau\geq 10$. In fact, we prove that 
for $\gamma \in (0,1/3)$,
\bels{bound on number of ev around tau}{
\#\setb{\,i:\2 \tau\leq\lambda_i\leq  \tau+ N^{-1}}
\;\prec\, 
N^{-\gamma}
\,.
}
For this we apply  Lemma \ref{lmm:Bounding measures by Stieltjes transforms} with the same choices of the measures $\nu_1$ and $\nu_2$ as in \eqref{choices of nus} and with
\bels{}{
\eta_1\,:=\, \eta_2\,:=\,\eps\,:=\, N^{\gamma-1}\,,\qquad \tau_1\,:=\, \tau\,,\qquad \tau_2\,:=\, \tau+N^{-1}.
}
Theorem \ref{thr:Local law} takes the form
\bels{LL well outside the spectrum}{
\big| \la\1 \vect{g}(\omega+\cI\1\eta)\ra -\la\1 \vect{m}(\omega+\cI\1\eta) \ra \big|\,\prec\, \frac{1}{N} +N^{-2\1\gamma},
\qquad
(\omega,\eta) \in \Gamma
\,,
}
where $ \Gamma := [\tau-N^{\gamma-1},\tau + 2N^{\gamma-1}] \times [N^{\gamma-1},2N^{\gamma-1}] $. 
Here we used $ \kappa(\omega+\cI\1\eta\1) \lesssim \eta_1 + (N\eta)^{-1} $, that follows from the facts that we are well outside $ \supp \rho \subset [-2,2\1] $, and hence $ \Delta(\omega) = 1 $ by \eqref{def of Delta_delta} so the condition \eqref{condition for improved kappa bound} holds, and thus \eqref{improved kappa bound} is applicable. 

Using the definition of stochastic domination (Definition~\ref{def:Stochastic domination}), the basic union bound, and the part (iii) of Lemma~\ref{lmm:Basic facts about stochastic domination}  we see that the estimate \eqref{LL well outside the spectrum} holds even with supremum over $ (\omega,\eta) \in \wht{\Gamma} $, where $ \wht{\Gamma} := (N^{-10}\Z)^2 \cap \Gamma  $ is a fine grid of spacing $ N^{-10} $ with $ \abs{\1\wht{\Gamma}\1} \leq N^{20} $.
Using the Lipschitz-continuity of $z \mapsto \avg{\1\brm{g}(z)\1}$ with Lipschitz-constant bounded by $N^2$, as well as the uniform $1/3$-H\"older-continuity of $ z \mapsto \avg{\1\vect{m}(z)\1}$, we can extend the supremum  over $ \wht{\Gamma} $  to the entire domain $ \Gamma $, i.e., 
\[
\sup_{(\omega,\eta) \ins \Gamma} \big| \la\1 \vect{g}(\omega+\cI\1\eta)\ra -\la\1 \vect{m}(\omega+\cI\1\eta) \ra \big|\,\prec\, \frac{1}{N} +N^{-2\1\gamma}.
\]
\nc
Plugging this bound into the definitions of $J_1$, $J_2$ and $J_3$ from \eqref{definition of Js} and using \eqref{Bounding measures by Stieltjes transforms} and the fact that $\rho=0$ in this regime shows the validity of \eqref{bound on number of ev around tau}. 

We conclude that a.w.o.p. there are no eigenvalues in an interval of length $N^{-1}$ to the right of $\tau$. By using a union bound this implies that
\[
\#\{\,i:\2 10\leq\lambda_i\leq  N\}\,=\, 0\qquad \text{a.w.o.p.}\,.
\]
The eigenvalues larger than $N$ are treated by the following simple argument,  
\[
\max_{i=1}^N\lambda_i^2 \,\leq\, \sum_{i=1}^N \lambda_i^2\,=\,  \sum_{i,j=1}^N|h_{ij}|^2\,\prec\, N\,.
\]
Thus \eqref{No evs beyond 10} holds true.

Now we apply Lemma \ref{lmm:Bounding measures by Stieltjes transforms} to prove \eqref{difference between distributions in the spectrum}. 
In case $|\tau|\geq 10$ the bound \eqref{difference between distributions in the spectrum} follows because a.w.o.p. there are no eigenvalues of $\vect{H}$ with absolute value larger or equal than $10$.
Thus, we
fix $\tau \in (-10,10)$ and make the choices
\bels{}{
\eta_1\,:=\, \eta_2\,:=\, N^{\gamma-1}\,,\qquad \tau_1\,:=\, -10\,,\qquad \tau_2\,:=\, \tau\,,\qquad \eps\,:=\,1\,.
}
Again we use \eqref{averaged local law} from
Theorem \ref{thr:Local law}, the Lipschitz-continuity of $\la \vect{g}\ra$ and the H\"older-continuity of $\la \vect{m}\ra$ to see that uniformly for all $\eta \geq N^{\gamma-1}$,
\[
\sup_{\omega \ins [0,\eta_1]}\!\big| \la\2 \vect{g}(\tau_1-\omega+\cI\1\eta)\1\ra -\la\2 \vect{m}(\tau_1-\omega+\cI\1\eta)\1 \ra \big|
\prec\, \frac{1}{N} +\frac{1}{(N\eta)^2}\,.
\]
Here we evaluated $\Delta(\tau_1)=1$ and thus $\kappa\lesssim \eta + (N\eta)^{-1}$.
With $J_1$ defined as in \eqref{definition of Js} we infer $J_1 \prec N^{-1}$. Theorem \ref{thr:Local law} also implies the bound
\[
\sup_{\omega \ins [-20,20]}\sup_{\eta \ins [1,2]}\big| \la\2 \vect{g}(\omega+\cI\1\eta)\1\ra -\la\2 \vect{m}(\omega+\cI\1\eta)\1 \ra \big|
\prec\, \frac{1}{N} \,,
\]
since in this regime $\kappa \lesssim 1$,
thus showing that $J_3 \prec N^{-1}$. We are left with estimating the three terms constituting $J_2$.  The first and second of these terms are estimated trivially by using the boundedness of their integrands. 
Therefore, we conclude that
\bels{bound on integrated distribution difference}{
\absbb{\int_{-10}^\tau\rho(\omega)\dif \omega -\frac{\#\{\,i : \, -10 \leq \lambda_i\leq \tau\}}{N}}
\,\prec\, 
N^{\gamma-1}
+
R(\tau)\,,
}
where the error term, $R$, is defined as
\bels{definition of R}{
R(\tau)\,:=\, 
 N^{1-\gamma}\int_{0}^{N^{\gamma-1}} \msp{-15}\dif \omega
\int_{N^{\gamma-1}}^{2}\msp{-8} \dif \eta\, 
\min\setbb{\frac{1}{N \eta(\Delta(\tau+\omega)^{1/3}+ \rho(\tau+\omega +\cI\1\eta))}\,,\,\frac{1}{\,(N  \eta)^{1/2}\!}}\,.
}
This expression is derived by using the bound \eqref{kappa bound} on $\kappa$ for the integrand of the third contribution to $J_2$.  

To estimate $R$ further we distinguish three cases, depending on whether $\tau$ is away from $\mathbb{M}$, close to an edge or close to a local minimum in the interior of $\supp \rho$. In each of these cases we prove 
\bels{bound on R}{
R(\tau)\,\prec\, \min\setbb{ 
\frac{1}{N (\Delta(\tau)^{1/3}+ \rho(\tau))}\,,\, \frac{1}{N^{4/5\!}\,}
}\,.
}
{\sc Away from $\mathbb{M}$: } In case $\dist(\tau,\mathbb{M})\geq \delta_*$, with $\delta_*$ the size of the neighborhood around the local minima from Theorem \ref{thr:Solution of the QVE}, we have $\Delta^{\!1/3}+\rho\sim 1$ and thus the $\eta$-integral in \eqref{definition of R} yields a factor comparable to $N^{-1}\log N$. Thus, $R(\tau)\prec N^{-1}$.
\\[0.3cm]
{\sc Close to an edge: } Let $\dist(\tau,\{\alpha_k,\beta_k\})\leq \delta_*$. Then from the size of $\rho$ at an internal edge, at the extreme edges and inside the gap (cf. \eqref{rho at an internal edge}, \eqref{rho around extreme edge} and \eqref{rho inside a gap} from Theorem \ref{thr:Solution of the QVE}) we see that 
\[
\Delta(\tau+\omega)^{1/3}+ \rho(\tau+\omega +\cI\1\eta)\,\sim\,\big( \Delta(\tau) + \dist(\tau, \{\alpha_k,\beta_k\})+\eta\big)^{1/3}.
\]
for any $\omega \in [0,N^{\gamma-1}]$ and $\eta \in [N^{\gamma-1},2]$. With this the size of $R$ is given by
\[
R(\tau)\,\sim\, 
\int_{N^{\gamma-1}}^{2}\msp{-8} \dif \eta\1 
\min\setbb{
\frac{1}{N \eta(\Delta(\tau)+  \dist(\tau, \{\alpha_k,\beta_k\})+\eta)^{1/3}}\,,\frac{1}{\,(N  \eta)^{1/2}\!}
}\,.
\]
Integrating over $\eta$ yields that
\[
R(\tau)\,\lesssim\,\min\setbb{\frac{ \log N}{N(\Delta(\tau)+  \dist(\tau, \{\alpha_k,\beta_k\}))^{1/3}}\,,\frac{1}{N^{4/5}\!}}\,.
\]
Now \eqref{bound on R} follows by using the size of $\rho$ from Theorem \ref{thr:Solution of the QVE} again. 
\\[0.3cm]
{\sc Close to an internal local minimum: } Suppose $|\tau-\tau_0|\leq \delta_*$ for some $\tau_0 \in \mathbb{M} \backslash \partial \supp \rho$. Then using the size of $\rho$ from \eqref{rho close to a local minimum} of Theorem \ref{thr:Solution of the QVE} we see that
\[
R(\tau)\,\sim\, 
\int_{N^{\gamma-1}}^{2}\msp{-8} \dif \eta\1 
\min\setbb{
\frac{1}{N \eta(\rho(\tau_0)+ |\tau-\tau_0|^{1/3}+\eta^{1/3})}\,,\,\frac{1}{(N  \eta)^{1/2}\!}
}\,.
\]
The bound \eqref{bound on R} follows by performing the integration over $\eta$.
\\[0.3cm]

This finishes the proof of  \eqref{bound on R}. We insert this bound into \eqref{bound on integrated distribution difference} and use that $\gamma$ was arbitrary. Thus, we find
\[
\absbb{\,
\int_{-10}^\tau\!\rho(\omega)\2\dif \omega \;-\;\frac{\!\#\{\,i : \, -10 \leq \lambda_i\leq \tau\}}{N}\,
}
\;\prec\;  
\min\setbb{ \frac{1}{N (\Delta(\tau)^{1/3}+ \rho(\tau))}\,,\, \frac{1}{N^{4/5}\!}
}
\,.
\]
This finishes the proof of \eqref{difference between distributions in the spectrum} since there are no eigenvalues below $-10$.

Now we prove \eqref{gaps in spectrum}. Let $\tau \in \R \backslash \supp \rho$. 
Suppose that for some $k=1, \dots ,K$ we have $|\tau-\beta_k| = \dist(\tau, \partial\supp \rho)$. The case when $\tau$ is closer to the set  $\{\alpha_k\}$ than to $\{\beta_k\}$ is treated similarly. 
Suppose further that 
\[
\tau \,\geq\, \alpha_k + \delta_k\,, 
\]
where $\delta_k$ are defined as in \eqref{definition of deltak} and $\delta_0=N^{\gamma-2/3}$. Note that there is nothing to show if $k>1$ and the size of the gap, $\alpha_k-\beta_{k-1}$, is smaller than $2\1\delta_k$, i.e., if such a $\tau$ does not exist. In particular, we have $\alpha_k-\beta_{k-1}=\Delta(\tau)\gtrsim N^{-1/2}$. We will show that a.w.o.p. there are no eigenvalues in an interval of length $N^{-2/3}$ to the right of $\tau$, i.e.
\bels{no eigenvalues to the right of tau in the gap}{
\#\setb{\,i:\, \tau \leq \lambda_i \leq \tau+N^{-2/3}}\;=\, 0 \qquad \text{a.w.o.p.}\,.
}

We apply 
Lemma \ref{lmm:Bounding measures by Stieltjes transforms} 
with the same choices of the measures $\nu_1$ and $\nu_2$ as in \eqref{choices of nus}. Additionally, we set
\bels{}{
\eta_1\,:=\, \eta_2\,:=\,\eps\,:=\, N^{-2/3}\,,\qquad \tau_1\,:=\, \tau\,,\qquad \tau_2\,:=\, \tau+N^{-2/3}.
}
We use the local law, Theorem \ref{thr:Local law}, to estimate the differences between the Stieltjes transforms of the two measures for the integrands in the definition of the three error terms, $J_1$, $J_2$ and $J_3$ from \eqref{definition of Js}. By the definition of $\delta_k$ the condition \eqref{condition for improved kappa bound} is satisfied inside the integrals and we use the improved bound, \eqref{improved kappa bound}, on $\kappa$. Indeed, we find
\[
\sup \big| \la\2 \vect{g}(\omega+\cI\1\eta)\ra -\la\1 \vect{m}(\omega+\cI\1\eta) \ra \big|\,\prec\, 
\frac{1}{N\delta_k\Delta(\tau)^{1/3}}
+\frac{1}{N^{2/3} \1\2\delta_k^{1/2}\Delta(\tau)^{1/6}\!}
\;,
\]
where the supremum is taken over $\omega \in [\tau- N^{-2/3}, \tau+ 2N^{-2/3}]$ and $\eta \in [N^{-2/3}, 2N^{-2/3}]$.
With this, the definition of $\delta_k$ and the size of $\rho$ from \eqref{rho inside a gap} and \eqref{rho around extreme edge} we infer
\[
J_1+J_2 +J_3\,\prec\, 
N^{-1-\gamma/2}
.
\]
From this \eqref{no eigenvalues to the right of tau in the gap} follows. The claim, \eqref{gaps in spectrum}, is now a consequence of a simple union bound taken over the events in \eqref{no eigenvalues to the right of tau in the gap} with different choices of $\tau$. This finishes the proof of Corollary~\ref{crl:Convergence of cumulative eigenvalue distribution}.
\qed

\subsection{Proof of Corollary \ref{crl:Rigidity of eigenvalues}}
Here we show how we get the rigidity, Corollary \ref{crl:Rigidity of eigenvalues}, from Corollary \ref{crl:Convergence of cumulative eigenvalue distribution}. Fix a $\tau \in [\alpha_1,\beta_K]$. We define the random fluctuation to the left, $\delta_-$, and to the right, $\delta_+$, of the eigenvalue $\lambda_{i(\tau)}$ as 
\begin{subequations}
\begin{align}
\delta_+(\tau)
\,&:=\,
\inf \setbb{ 
\delta \2\geq\20\2:\, 
2+\Big| \#\setb{\,i:\, \lambda_i \leq \tau+\delta}\,-\,N\!\int_{-\infty}^{\tau+\delta}\msp{-10}\rho(\omega)\1\dif \omega \Big|
\2\leq\2 N\!
\int_{\tau}^{\tau+\delta}\msp{-10}\rho(\omega)\1\dif \omega
}
\\
\delta_-(\tau) \,&:=\,  
\inf\setbb{ 
\delta \2\geq\20\2:\, 
1+
\Big| \#\setb{\,i:\, \lambda_i \leq\tau-\delta}\,-\,N\!\int_{-\infty}^{\tau-\delta}\msp{-10}\rho(\omega)\1\dif \omega \Big|
\2\leq\2 N\!
\int_{\tau-\delta}^{\tau}\msp{-2}\rho(\omega)\1\dif \omega
}
\,.
\end{align}
\end{subequations}
We show now that with this definition,
\bels{left and right fluctuation}{
\lambda_{i(\tau)} \,\in\, \bigl[\1\tau-\delta_-(\tau),\tau+\delta_+(\tau)\bigr]
\,.
}

We start with the upper bound on $\lambda_{i(\tau)}$. By the definition of $i(\tau)$ we find the inequality
\[
\#\setb{\,i:\, \lambda_i \leq \lambda_{i(\tau)}}
\,=\, 
i(\tau)\,\leq\, 1+N\int_{-\infty}^{\tau}\rho(\omega)\dif \omega
\,=\, 1+ N\int_{-\infty}^{\tau+\delta_+}\msp{-15}\rho(\omega)\dif \omega- N\int_{\tau}^{\tau+\delta_+}\msp{-15}\rho(\omega)\dif \omega\,.
\]
The definition of $\delta_+=\delta_+(\tau)$ implies that 
\[
\#\setb{\,i:\, \lambda_i \leq \lambda_{i(\tau)}}
\;<\, 
\#\{\,i:\, \lambda_i \leq \tau +\delta_+\}\,.
\]
By monotonicity of the cumulative eigenvalue distribution, we conclude that $ \lambda_{i(\tau)}\leq\tau +\delta_+$. Thus, the upper bound is proven. 

Now we show the lower bound. We start similarly,
\[
\#\setb{\,i:\, \lambda_i \leq \lambda_{i(\tau)}}
\;=\, 
i(\tau)
\,\ge\, 
N\!\int_{-\infty}^{\tau}\rho(\omega)\dif \omega
\,=\, 
N\!\int_{-\infty}^{\tau-\delta_-}\msp{-15}\rho(\omega)\dif \omega 
\,+\,  
N\!\int_{\tau-\delta_-}^{\tau}\msp{-10}\rho(\omega)\dif \omega
\,.
\]
By definition of $\delta_-$ we get
\[
\#\setb{\,i:\, \lambda_i \leq \lambda_{i(\tau)}}
\;\,\geq\, 
1+\, \liminf_{\eps\downarrow 0}\#\setb{\,i:\, \lambda_i \leq \tau-\delta_-\!-\eps}
\,.
\]
Here the $\liminf$ is necessary, since the cumulative eigenvalue distribution is not continuous from the left. We conclude that
$\lambda_{i(\tau)} \geq \tau-\delta_--\eps$ for all $\eps>0$ and therefore the lower bound is proven.

Now we start with the proof of \eqref{rigidity of eigenvalues in support}. For this we show that for any $\tau$ that is well inside the support of the density of states, i.e., that satisfies 
\eqref{tau well inside the support}, we have
\bels{upper bound for delta+ and delta-}{
\delta_-(\tau)+\delta_+(\tau)\,\prec\, \delta\,,\qquad\delta\,:=\, 
\min\setbb{
\frac{1}{\rho(\tau)(\Delta(\tau)^{1/3}+\rho(\tau))N}\,,\frac{1}{N^{3/5}\!}
}
\,.
}
If $\tau$ is in the bulk, i.e., $\dist(\tau,\mathbb{M})\geq \delta_*$, then $\delta \sim N^{-1}$ and thus 
\eqref{upper bound for delta+ and delta-} follows from \eqref{difference between distributions in the spectrum}. 
We distinguish the two remaining cases, namely whether $\tau$ is close to an edge or to a local minimum inside the interior of $\supp \rho$.
\\[0.3cm]
{\sc Close to an edge: }  Suppose that $\tau \in [\beta_k-\delta_*,\beta_k-\eps_k]$. The case when $\tau$ is closer to $\{\alpha_k\}$ than to $\{\beta_k\}$ is treated similarly.  By the definition of $\eps_k$ in \eqref{definition of epsk} and by the size of $\rho$ from \eqref{rho around extreme edge} and \eqref{rho at an internal edge} in Theorem \ref{thr:Solution of the QVE} we see that $\eps_k \gtrsim N^\gamma \delta$. Using Corollary \ref{crl:Convergence of cumulative eigenvalue distribution} we find for any $\eps\in (0,\gamma/2)$ that 
\[
\absbb{\,\#\setb{\,i:\, \lambda_i \leq \tau+N^\eps\delta\2}\,-\,N\!\int_{-\infty}^{\tau+N^\eps\delta}\msp{-15}\rho(\omega)\1\dif \omega\,}
\,\prec\, 
\min\setB{  
\bigl(\Delta(\tau)+\beta_k-\tau\bigr)^{-1/3},\2
N^{1/5}
}
\,.
\]
On the other hand
\[
N
\int_{\tau}^{\tau+N^\eps\delta}\msp{-10}\rho(\omega)\dif \omega\,\sim\, \frac{N^{1+\eps}\delta\2(\beta_k-\tau)^{1/2}}{(\Delta(\tau)+\beta_k-\tau)^{1/6}}
\,\gtrsim\, 
N^\eps 
\min\setB{  
\bigl(\Delta(\tau)+\beta_k-\tau\bigr)^{-1/3},\2
N^{1/5}
}
\,.
\]
Here we used the size of $\rho$ from Theorem \ref{thr:Solution of the QVE}, the definition of $\delta$ and $\beta_k-\tau \geq \eps_k$.  Since $\eps$ was arbitrary we conclude that $\delta_+(\tau)\prec \delta$. The bound, $\delta_-(\tau)\prec\delta$, is shown in the same way. 
\\[0.3cm]
{\sc Close to internal local minima: } Suppose $|\tau-\tau_0|\leq \delta_*$ for some $\tau_0 \in \mathbb{M}\backslash \partial \supp \rho$.  Then by \eqref{rho close to a local minimum} with $\Delta(\tau_0)=0$ and the definition of $\delta$ in \eqref{upper bound for delta+ and delta-} we have 
\[
\delta\,\sim\, \min\setbb{
\frac{1}{(\rho(\tau_0)^3+|\tau-\tau_0|)^{2/3}N},\frac{1}{N^{3/5}}
}
\,.
\]
We apply \eqref{difference between distributions in the spectrum} from Corollary \ref{crl:Convergence of cumulative eigenvalue distribution} and, using \eqref{rho close to a local minimum} again, we get
\bels{internal minima bound 1}{
\absbb{ 
&\#\setb{\,i:\, \lambda_i \leq \tau+N^\eps\delta}\,-\,N\!\int_{-\infty}^{\tau+N^\eps\delta}\msp{-10}\rho(\omega)\1\dif \omega
\,
}
\\
&\msp{240}\prec\;
\min\setB{  
\bigl(\2\rho(\tau_0)^3+|\tau+N^\eps\delta-\tau_0|\bigr)^{-1/3}\!,\2
N^{1/5}
}
\,.
}
On the other hand, we find
\bels{internal minima bound 2}{
N
\int_{\tau}^{\tau+N^\eps\delta}\msp{-10}\rho(\omega)\1\dif \omega
\;\sim\; 
N^{1+\eps}\delta \,\Bigl(\,\rho(\tau_0)^3 + |\tau-\tau_0|+N^\eps\delta\,\Bigr)^{1/3}
.
}
We will now verify that for large enough $N$,
\bels{Estimate for delta+}{
&N^{\eps/2}\min\setB{  
\bigl(\,\rho(\tau_0)^3+|\tau+N^\eps\delta-\tau_0|\2\bigr)^{-1/3}\!,\2
N^{1/5}
}
\\
&\msp{300}\lesssim\;
N^{1+\eps}\delta \,\Bigl(\2\rho(\tau_0)^3+|\tau-\tau_0|+N^\eps\delta\2\Bigr)^{1/3}
.
}

We distinguish three cases. First let us consider the regime where $\rho(\tau_0)^3+|\tau-\tau_0|\leq N^{-3/5}$. Then we have $\delta = N^{-3/5}$ and 
\[
N^{1+\eps}\delta \2\bigl(\2\rho(\tau_0)^3+|\tau-\tau_0|+N^\eps\delta\,\bigr)^{1/3}\,\sim\, N^{4\eps/3}N^{1/5}\,.
\]

Now we treat the situation where, $N^{-3/5}<\rho(\tau_0)^3+|\tau-\tau_0|\leq N^{3\eps/2-3/5}$. In this case 
\[
N^{1+\eps}\delta \,\bigl(\,\rho(\tau_0)^3+|\tau-\tau_0|+N^\eps\delta\,\bigr)^{1/3}
\,\gtrsim \,  \frac{N^\eps}{(\1\rho(\tau_0)^3+|\tau-\tau_0|\1)^{1/3}}\,\geq\, N^{\eps/2}N^{1/5}\,.
\]

Finally, we consider $\rho(\tau_0)^3+|\tau-\tau_0|> N^{3\eps/2-3/5}$. Then for large enough $N$ we find on the one hand
\[
\min\setB{  
\bigl(\2\rho(\tau_0)^3+|\tau+N^\eps\delta-\tau_0|\bigr)^{-1/3},\2
N^{1/5}
}
\,\sim\, \frac{1}{(\1\rho(\tau_0)^3+|\tau-\tau_0|\1)^{1/3}}\,,
\]
and on the other hand
\[
N^{1+\eps}\,\delta \,\bigl(\2\rho(\tau_0)^3+|\tau-\tau_0|+N^\eps\delta\,\bigr)^{1/3}\,\gtrsim \,  \frac{N^\eps}{(\1\rho(\tau_0)^3+|\tau-\tau_0|\1)^{1/3}}\,.
\]
Thus, \eqref{Estimate for delta+} holds true and since $\eps$ was arbitrary, we infer from \eqref{internal minima bound 1} and \eqref{internal minima bound 2} that $\delta_+(\tau)\prec\delta$. Along the same lines we prove $\delta_-(\tau)\prec\delta$. 
Thus \eqref{upper bound for delta+ and delta-} and with it \eqref{rigidity of eigenvalues in support} are proven.

The statement about the fluctuation of the eigenvalues at the leftmost edge, \eqref{extreme edge fluctuation} follows directly from \eqref{rigidity of eigenvalues in support} and \eqref{gaps in spectrum} in Corollary \ref{crl:Convergence of cumulative eigenvalue distribution}. Indeed, for $\tau \in [\alpha_1,\alpha_1+\eps_0)$ we have $\lambda_{i(\tau)}\leq \lambda_{i(\alpha_1+\eps_0)}$ and from \eqref{rigidity of eigenvalues in support} with $\Delta(\tau)=1$, as well as $\rho(\alpha_1+\eps_0)\sim \eps_0^{1/2}$, and from the definition of $\eps_0$ we see that 
\[
\lambda_{i(\alpha_1+\eps_0)}\,\leq\, \alpha_1 +\eps_0+N^{\gamma-2/3}\,\leq\, \tau +2\1N^{\gamma-2/3}\qquad \text{a.w.o.p.}\,.
\]
On the other hand, \eqref{gaps in spectrum} shows that a.w.o.p. $\lambda_{i(\tau)}\geq \alpha_1-N^{\gamma-2/3}$. Since $\gamma$ was arbitrary,  \eqref{extreme edge fluctuation} follows. 
The rigidity at the rightmost edge is proven along the same lines.  

The claim, \eqref{edge eigenvalue rigidity}, about the remaining eigenvalues follows from a similar argument. For $\tau \in (\beta_k-\eps_k,\alpha_{k+1}+\eps_k)$, as a consequence of \eqref{gaps in spectrum}, we have 
\[
\lambda_{i(\tau)} \,\in\, \bigl[\2\lambda_{\1i(\beta_k-\eps_k)}, \beta_k+\delta_k\bigr]
\cup 
\bigl[\2\alpha_{k+1}-\delta_k\2, \lambda_{\1i(\alpha_{k+1}+\eps_k)}\bigr]
\qquad \text{a.w.o.p.}\,.
\]
From \eqref{rigidity of eigenvalues in support} and the definition of $\eps_k$ we infer $\lambda_{i(\beta_k-\eps_k)}\geq \beta_k-2\1\eps_k $ a.w.o.p., as well as $\lambda_{i(\alpha_{k+1}+\eps_k)} \leq \alpha_{k+1}+2\1\eps_k$ a.w.o.p., which finishes the proof of  \eqref{edge eigenvalue rigidity}. 
\qed

\subsection{Proof of Corollary \ref{crl:Delocalization of eigenvectors}}

The delocalization of eigenvectors is a simple consequence of the anisotropic local law Theorem~\ref{thr:Anisotropic law} using the argument from \cite{EKYY}. 
Expressing the resolvent in the eigenbasis, we have
\bels{spectral representation of b cdot Gb}{
\brm{b} \cdot \brm{G}(z)\brm{b}\,=\, \sum_{i=1}^N \frac{\abs{\2\brm{b} \cdot \brm{u}^{(i)}}^2\msp{-7}}{\!\lambda_i \,-\,z}
\,,
}
where $\brm{u}^{(i)} $ is the $ \ell^{\12}$-normalised eigenvector corresponding to the eigenvalue $\lambda_i$. 
We evaluate this at $z:=\lambda_k+\cI\1N^{\gamma-1}$ with $\gamma>0$ as in the statement of Theorem \ref{thr:Anisotropic law}. 
The anisotropic local law implies that also  $\brm{b} \cdot \brm{G}(z)\brm{b}$ is uniformly bounded. Hence we get
\[
1\,\gtrsim\, 
\Im\;\brm{b} \cdot \brm{G}(z)\brm{b}
\,\ge\,  
N^{1-\gamma}\2
\abs{\2\brm{b} \cdot \brm{u}^{(k)}}^2
\,,
\]
by keeping only a single summand $ i = k $ from \eqref{spectral representation of b cdot Gb}. As $\gamma>0$ was arbitrary we conclude 
that 
\[ 
\abs{\2\brm{b} \cdot \brm{u}^{(k)}} \,\prec\, N^{-1/2}\,, 
\]
uniformly in $ k $.
\qed

\section{Anisotropic law and universality}

\subsection{Proof of Theorem \ref{thr:Anisotropic law}}

Given the entrywise local law, Theorem \ref{thr:Local law}, the proof of the anisotropic law follows exactly as Section 7 in \cite{EKYYiso}, where the same argument was presented for generalized Wigner matrices (this argument itself mimicked the detailed proof of the isotropic law for sample covariance matrices in Section 5 of \cite{EKYYiso}).
The only difference is that in our case $G_{ii}(z)$ is close to $m_i(z)$, the $ i$-th component of the solution to the QVE, which now  genuinely depends on $i$, while in \cite{EKYYiso} we had 
$ G_{ii}(z)\approx m_{\mathrm{sc}}(z)$ for every $i$, where $m_{\mathrm{sc}}(z)$ is the solution to \eqref{scemsc}.
However, the diagonal resolvent elements played no essential role in \cite{EKYYiso}. 
We now explain the small modifications.

Recall from Section 5.2 of \cite{EKYYiso} that by polarization it is sufficient to prove \eqref{isoeq} for $\ell^{\12}$-normalized vectors $\vect{w}= \vect{v}$.  
We can then write 
\[
\sum_{i,j=1}^N \overline{v}_i\, G_{ij} v_j -\sum_{i=1}^N m_i |v_i|^2 \,=\, \sum_i (G_{ii}-m_i)|v_i|^2 + \mathcal{Z}
,\qquad
\mcl{Z}:= \sum_{i\ne j}^N \overline{v_i}\2 G_{ij} v_j
\,.
\]
The first term containing the diagonal elements $G_{ii}$ is clearly bounded by the right hand side of \eqref{isoeq} by Theorem \ref{thr:Local law}. This is the first instant  where the nontrivial $i$-dependence of $m_i$ is used. 

The main technical part of the proof in  \cite{EKYYiso} is then to control $ \mcl{Z} $, the contribution of the off diagonal terms. We can follow this proof in our case to the letter; the nontrivial $i$-dependence of $m_i$ requires a slight modification only at one point. 
To see this, we recall the main structure of the proof. For any even $ p $, the moment 
\bels{EZ}{
\EE |{\cal Z}|^p \,=\, \EE\!
\sum_{b_{11}\ne  b_{12}} \!\cdots\! \sum_{b_{p1}\ne b_{p2}} 
\biggl(\; 
\prod_{k\1=\11}^{p/2} \overline{v}_{b_{k1}} G_{b_{k1}b_{k2}}v_{b_{k2}}\biggr)
\biggl(\msp{-13} 
\prod_{\msp{20}k=p/2+1}^{p} \msp{-15}\overline{v}_{b_{k1}} G^*_{b_{k1}b_{k2}}v_{b_{k2}}
\biggr)
\,,
}
is computed. 
Let us concentrate on a fixed summand in \eqref{EZ}  and let $B=\{b_{k_1}\}\cup\{b_{k_2}\}$ be the set of $\vect{v}$-indices appearing in that term. 
Using the resolvent identity \eqref{resolvent identity} 
we successively expand the resolvents until each of them appears in a \emph{maximally expanded} form, where every resolvent entry is of the form $ G^{(B\backslash ab)}_{ab} $, for some $ a,b \in B $ (cf. Definition 5.4 of \cite{EKYYiso}).
Each time a maximally expanded off-diagonal element is produced we use \eqref{off-diagonal expansion}. 
Finally, unless we end up with an expression that contains a very large numbers of off-diagonal resolvent entries (such \emph{trivial leaves} are treated separately in Subsection 5.11 of \cite{EKYYiso})
we apply \eqref{diagschur} to expand the remaining maximally expanded diagonal resolvent entries.
This way we end up with an expression where only the resolvent entries of the type $ G^{(B)}_{ij} $, with $ i,j \notin B $, appear. In other words, the $ \brm{v} $-indices and the indices of the resolvent entries are completely decoupled;  only explicit products of entries of $ \brm{H} $ represent the connections between them.
We can now take partial expectation w.r.t. the rows and columns of these $h$-terms. In this way we guarantee that each index in $ B $ appears at least twice as a value of $ b_{k1} $ or $ b_{k2}$ in \eqref{EZ}, i.e., the entries of $ \brm{v} $  must be at least paired, and therefore the $2p$-fold  summation in \eqref{EZ} effectively becomes at most a $p$-fold summation. 
This renders the uncontrolled $\ell^1$-norm of $ \brm{v} $ to $\ell^{\1q}$-norms of $ \brm{v}$, with $ q \ge 2$, which are bounded by one by normalization.

Along this procedure it is only at the treatment of the maximally expanded diagonal resolvent elements appearing in the non-trivial leaves (cf. Subsection 5.12 of \cite{EKYYiso}) where we need to slightly adjust the proof to the setting where $ \brm{S} $ is not stochastic. 
Using the QVE \eqref{SCE} and Schur's formula, similarly as in \eqref{diagschur}, we obtain a representation, where all the dependence of the $ B$-columns and -rows of $ \brm{H} $ is explicit
\bels{Laszlo's: w-expansion for 1/G}{
\frac{1}{G^{(B\backslash b)}_{bb}\!} 
\,&=\,
\frac{1}{m_b} 
\,-\,
\sum_{i,j}^{(B)}
\Bigl(\,h_{bi}\2G^{(B)}_{ij} h_{jb}-s_{bi}\1m_i\2\delta_{ij}\Bigr)
\,+
\sum_{a \in B} s_{ba}m_a
\,+\,h_{bb} 
\,,
\quad b \in B
\,.
}
This formula replaces (5.41) from \cite{EKYYiso}.
Taking the inverse of this formula and expanding around the leading term $ m_b  $, we get a geometric series representation for $ G^{(B\backslash b)}_{bb} $ in terms of powers of the last three term in \eqref{Laszlo's: w-expansion for 1/G}.
The resulting formula is analogous to (5.42) in \cite{EKYYiso}. 
The geometric series converges 
because the last three term on the right hand side of \eqref{Laszlo's: w-expansion for 1/G} are much smaller than $ \abs{1/m_b}\sim 1 $ a.w.o.p..  
Indeed,  the last two terms in \eqref{Laszlo's: w-expansion for 1/G} are of size $ N^{-1} $ and $ N^{-1/2+c} $ a.w.o.p., respectively. 
The double sum in \eqref{Laszlo's: w-expansion for 1/G} is small by using the large deviation estimates \eqref{large deviation off diagonal}--\eqref{large deviation  diagonal}, similarly as in the proof of Lemma \ref{lmm:Bound on perturbation}. When estimating the  diagonal sum $i=j$,  
we note that $ \abs{G^{(B)}_{ii}-m_i} $ is small by first estimating $ \abs{G^{(B)}_{ii}-G_{ii}} $ similarly to \eqref{G^T-G is small}, and then we use the local law Theorem \ref{thr:Local law} to see that also $ \abs{G_{ii}-m_i}  $ is small.

The proof in \cite{EKYYiso} did not use the specific form of the subtracted term $s_{bi}m_i \delta_{ij} $ in \eqref{Laszlo's: w-expansion for 1/G}, just the fact that the subtraction made  \eqref{large deviation diagonal} applicable for the double summation in \eqref{Laszlo's: w-expansion for 1/G}. After this slight modification, the rest of the proof in \cite{EKYYiso} goes through without any further changes. 
\qed

\subsection{Proof of Theorem \ref{thr:Universality}}\label{sec:uni}

For the proof of Theorem \ref{thr:Universality} we follow the method developed in \cite{ESYY, EYY, EKYY}.
Theorem 2.1 from \cite{EK} was designed for proving universality for a random matrix with a small independent Gaussian component and densities of state that may differ from Wigner's semicircle law. 
The main theorem 
in \cite{EK} asserts that if local laws hold in a sufficiently strong sense then bulk universality holds locally for matrices with a small Gaussian component.  
We remark that a similar approach was independently developed in \cite{LY} that can also be easily used to conclude bulk universality from Theorem \ref{thr:Local law}, but here we follow \cite{EK}.
In Section 2.5 of \cite{EK} a recipe was given how to use this theorem to establish universality for a quite general class of random matrix models even without the Gaussian component, as long as uniform local laws on the optimal scale are known and the matrix satisfies  the appropriate $ q $-fullness condition (cf. Definition \ref{def:q-full})  that allows for an application of the moment matching (Lemma 6.5 in \cite{EYY}) and the Green's function comparison theorem (Theorem 2.3 in \cite{EYY}). 

Let $ \brm{H} $ be the Wigner-type matrix satisfying the hypotheses of Theorem \ref{thr:Universality}, and for which the universality is to be proven. 
Let $ \tau $ be a bulk point of $ \rho $, so that $ \rho(\tau) \ge \eps $, for some $ \eps > 0 $, and let $ I := [\tau-\delta,\tau+\delta\1] $ be some environment of size $ \delta \sim 1  $ around $ \tau $. Following the above recipe, it remains to show that the local law holds for the random matrices
\[
\vect{H}_t\,=\,  \nE^{-t/2}\1\vect{H}_0 + (1- \nE^{-t}\1)^{\msp{-1}1/2}\1\vect{U} ,
\]
uniformly in both  $t \in [0,T] $ and the spectral parameters $ z \in I +\cI\2[N^{\gamma-1},\infty) $. 
Here $ T $ is a small negative power of $N $, i.e., $T=N^{-\xi}$ for some $\xi> 0 $, such that $ \brm{H} $ and $ \brm{H}_T $ are close in the four moment comparison sense (cf. Theorem 2.3 of \cite{EYY}), and $ \brm{U}$ is a standard GUE/GOE random matrix.
The random matrix $ \brm{H}_0$ has independent entries, is independent of $\vect{U}$, and has a variance matrix 
\[
\vect{S}_0\,:=\, \nE^{T}\vect{S}\2-\2(\nE^{T}-1)\1\vect{S}_{\rm G}
\,,
\]
with $ \brm{S}$ and $\brm{S}_{\rm G}$ denoting the variance matrices of $ \brm{H} $ and the standard GUE/GOE-matrix, respectively.
It follows that the variance matrix of $\vect{H}_t$ is
\[
\vect{S}_t\,=\, \nE^{-t}\vect{S}_0+(1-\nE^{-t})\1\vect{S}_{\mathrm{G}}
\,,
\]
and hence $ \brm{S}_T =  \vect{S} $ as required by the moment matching. 

We will now show that $ \brm{H}_t $ satisfy the hypotheses of Corollary~\ref{crl:Local law in the bulk} uniformly in $ t$.
Since $ T=N^{-\xi} $ is small, the variance matrices $ \brm{S}_t $ are all small perturbations of $ \brm{S} $. 
In particular, $ \brm{S}_t $, $ t \in [0,T] $, are hence $ q/2 $-full.

Next we show that the interval $ I $ is inside the bulk of $ \brm{H}_t $. To this end, we consider the QVE associated to the variance matrix $ \brm{S}_t $,
\[
-\frac{1}{\2m_{t\1;\1i}\!} \,=\, z \2+\2 (\brm{S}\brm{m}_t)_i \1+\2 d_i\,,
\qquad
\brm{d} \,=\, (\1\brm{S}_t-\brm{S}\1)\1\brm{m}_t\,,
\]
as a perturbation of the original QVE with $ \brm{S} = \brm{S}_T $.
In order to use our stability results we show $ \norm{\brm{d}}_\infty \lesssim T $. 
Since $ \brm{H}_t $ is $ q/2$-full we have  $ s_{t;ij} \ge q/2 $ and hence using (i) of Theorem~6.1 of \cite{AEK1} we see that there is a constant $ \delta' \sim 1  $ such that $ \norm{\brm{m}_t(z)}_\infty \sim 1 $ uniformly for $ \abs{\Re\,z} \leq \delta' $.
Moreover,  the structural $ \Lp{2}$-bound from Theorem~2.1 of \cite{AEK1} implies
\[
\frac{\norm{\brm{m}_t(z)}^2_{\ell^2}}{N} = \frac{1}{N}\sum_{i\1=\11}^N \2\abs{\1m_{t\1;\1i}(z)}^2\,\leq\, \frac{4}{\abs{z}^2}\,,\qquad z \in \Cp\,,\; t \in [0,T]\,.
\]
Combining these estimates we see that $ \sup_{t,z}\norm{\brm{m}_t(z)}^2_{\ell^2} \lesssim N $, and consequently the perturbation is small in the uniform norm: $ \norm{\brm{d}}_\infty \lesssim N\sup_{i,j}\abs{s_{t;ij}-s_{ij}} \lesssim N^{-\xi} $.
Applying the stability (Theorem~\ref{thr:Stability around support} or Theorem~2.12 from \cite{AEK1}) of the QVE associated to $ \brm{S} $ we conclude that $ \norm{\1\brm{m}_t(z)-\brm{m}(z)}_\infty \lesssim N^{-\xi}\eps^{-2} $, and hence $ \rho_t(\omega) \ge \eps/2 $ for $\omega \in I $ and all $ t $, provided $ N $ is sufficiently large. 

The moment condition (D) is automatically satisfied uniformly for every $\vect{H}_t$ by construction.
Since the condition (A) is merely a matter of normalization we have now shown that $ \brm{H}_t $ satisfy the hypotheses of Corollary~\ref{crl:Local law in the bulk} uniformly in $ t $. Thus $ \brm{H}_t $ satisfy local law uniformly in $t \in [0,T] $ and $ z \in I +\cI\2[N^{\gamma-1},\infty) $. This finishes the proof of universality. 
\qed

\appendix 
\section{Appendix}

The relation $ \prec $ is transitive and 
it satisfies the following arithmetic rules:

\begin{lemma}[Basic facts about stochastic domination]
\label{lmm:Basic facts about stochastic domination}
We have:
\begin{itemize}
\titem{i} If $ \phi \prec N^\delta\psi $, for every $ \delta > 0 $, then $ \phi \prec \psi\, $;  
\titem{ii} If $ \phi \prec N^{-\delta}\phi+\psi $, for some $ \delta > 0$, then $ \phi \prec \psi\, $.
\end{itemize}
Let $ \phi_u $ and $ \psi_u $ be some non-negative random variables parametrized by elements $ u $ of some set $ \mathbb{U} $, such that $ \phi_u \prec \psi_u $, uniformly in $ u $.
If $ \mathbb{U}' \subset \mathbb{U} $, then  
%
\begin{itemize}
\titem{iii} 
$ \sum_{u \in \mathbb{U}'} \phi_u \prec\, \sum_{u \in \mathbb{U}'} \psi_u $, provided $ \abs{\1\mathbb{U'}} \leq N^C $ for some $ C < \infty\, $;
\titem{iv} 
$ \prod_{u \in \mathbb{U}'} \phi_u \prec\, \prod_{u \in \mathbb{U}'} \psi_u $, provided $ \abs{\1\mathbb{U'}} \leq C $, for some $ C < \infty\, $.
\end{itemize}
\end{lemma}
These properties follow directly from the definition (Definition~\ref{def:Stochastic domination}) of stochastic domination. For further details see~\cite{EKYY}.

\NLemma{Bound propagation}{Suppose $C_1,D_1,D_2,D_3$ and $\eps_1$ are positive constants, depending explicitly on $p$, $P$, $L$, $\ul{\mu}$, $\gamma$ and possible on additional parameters in some set $V$. Suppose further that the threshold function $N_0$ from Definition \ref{def:Overwhelming probability} depends on the same parameters.
Let $\D^{(N)} \subseteq \Cp$ be a sequence of connected subsets of the complex upper half plane with only polynomially growing diameter, $\sup\{|z_1-z_2|: z_1,z_2 \in\D^{(N)} \}\leq N^{D_1}$. Let
$\varphi = (\varphi^{(N)}(z): z \in  \D^{(N)})_{N \in \N}$ be a sequence of non-negative random functions and $\Phi^{(N)}:\D^{(N)} \to (N^{-D_3},\infty)$ a sequence of deterministic functions on these sets. Suppose they satisfy the following conditions:
\begin{itemize}
\item Uniformly for all $z_1,z_2 \in \D^{(N)}$
\bels{weak continuity condition}{
|\varphi^{(N)}(z_1)-\varphi^{(N)}(z_2)|+|\Phi^{(N)}(z_1)-\Phi^{(N)}(z_2)|\,\leq\, C_1\2N^{D_2}|z_1-z_2|^{\eps_1} \,.
}
\item Uniformly for all $z \in \D^{(N)}$ 
\bels{gap condition}{
\text{a.w.o.p.}\qquad 
\varphi^{(N)} \,\notin\, \bigl[\2\Phi^{(N)}(z)-N^{-D_3},\2\Phi^{(N)}(z)\2\bigr]
\,. 
}
\item There is a sequence $z_0^{\2(N)} \in \D^{(N)}$ such that
\bels{initial value condition}{
\text{a.w.o.p.}\qquad
\varphi^{(N)}(\2z_0^{\2(N)})\,\leq\, \Phi^{(N)}(\2z_0^{\2(N)})
\,.
}
\end{itemize}
Then the sequence $\varphi$ satisfies the bound
\bels{}{
\text{a.w.o.p.}\qquad \text{for all }\; z \in \D^{(N)}\;:\quad
\varphi^{(N)}(z)\,\leq\, \Phi^{(N)}(z)\,.
}
}
\begin{Proof}
We will not carry the upper index $N$ in this proof.
First we choose a grid $\G \subseteq \D$ with the following properties
\begin{itemize}
\item The number of points in $\G $ is polynomially large, i.e.,
$|\G|\leq C_2N^{D_4}$.
\item The grid is connected and sufficiently dense in $\D$, i.e., for any two points $z_1,z \in \G$ there is a path $(z_i)_{i=2}^K \subseteq \G$, such that $\max\{|z_K-z|,|z_{i+1}-z_i|\}\,\leq\, N^{-D_5}$ for all $i=1,\dots,K-1$.
\end{itemize}
Here, the positive exponent $D_5$ is chosen  sufficiently large such that $C_1\1N^{D_2-\eps_1\1D_5}\leq N^{-D_3}/2$. Then an upper bound on the positive constants $D_4$ and $C_2$ is determined by the choice of $D_5$ and the diameter of $\D$, i.e., by $D_1$. 

Now let $z \in \G$. Then we find a path $(z_i)_{i=1}^{K}$ in $ \G$ that connects $z_0$ with $z_{K+1}:=z$ in the sense of the second property of $\G$. We may assume the length of the path, $K$, to be bounded by $|\G|$. Inductively we show that 
for all $ i=0,\dots, K+1$
\[
\text{a.w.o.p.}\qquad \varphi(z_i)\,\leq\, \Phi(z_i)-N^{-D_3}.
\]
For $i=0$ this follows from \eqref{initial value condition} and \eqref{gap condition}. For all other $i$ it follows by induction using  the continuity condition \eqref{weak continuity condition}, which implies
$|\varphi(z_{i+1})-\varphi(z_{i})|+|\Phi(z_{i+1})-\Phi(z_{i})|\leq N^{-D_3}/2$. This shows that if $\varphi(z_{i})\leq \Phi(z_{i})-N^{-D_3}$, then $\varphi(z_{i+1})\leq \Phi(z_{i+1})$ and with \eqref{gap condition} even that $\varphi(z_{i+1})\leq \Phi(z_{i+1})-N^{-D_3}$. In particular, $\varphi(z)\leq \Phi(z)-N^{-D_3}$ a.w.o.p..

Using a union bound we infer that
\[
\text{a.w.o.p.}\qquad \text{for all } z \in \G \qquad \varphi(z)\,\leq\, \Phi(z)-N^{-D_3}.
\]
By \eqref{weak continuity condition} and since $\G$ is sufficiently dense in $\D$ this bound extends to all $z \in \D$ and the lemma is proven.
\end{Proof}

\begin{Proof}[Proof of Lemma \ref{lmm:Bounding measures by Stieltjes transforms}]
For $f, \chi$ compactly supported on $\R$ the Cauchy integral formula holds true,
\begin{equation*}
\begin{split}
&f(\tau)\,=\, \frac{1}{\pi}\int_{\R^2}\frac{\partial_{\overline z}\widetilde f(\sigma + \cI \1\eta)}{\tau-\sigma -\cI\1\eta}\dif \sigma \dif \eta 
\,=\, \frac{1}{2\pi} \int_{\R^2} \frac{\cI\1\eta  f''(\sigma)\chi(\eta) + \cI(f(\sigma) + \cI\1\eta \1 f'(\sigma))\chi'(\eta)}{\tau-\sigma - \cI\1\eta}\dif \sigma \dif \eta\,,
\\
&\widetilde f(\sigma + \cI\1\eta)\,:=\, (f(\sigma)+ \cI\1\eta\1 f'(\sigma))\chi(\eta)\,.
\end{split}
\end{equation*}
For a signed measure $\nu$ on $\R$ this implies the formula
\[
\int_\R f(\tau)\1\nu(\dif \tau)\,=\, \Re\2\int_\R f(\tau)\1\nu(\dif \tau)\,=\,-\frac{1}{2\pi}\bigl(\,
L_1(\nu)+L_2(\nu)+L_3(\nu)\,
\bigr)\,,
\]
where the three integrals $L_1$, $L_2$ and $L_3$ are given as
\begin{equation*}
\begin{split}
L_1(\nu)\,&:=\, \int_{\R^2} \eta \1 f''(\sigma)\1\chi(\eta)\1\Im\2 m_\nu(\sigma+\cI\1\eta)\2 \dif \sigma \dif \eta\,,
\\
L_2(\nu)\,&:=\, \int_{\R^2} f(\sigma)\chi'(\eta) \1\Im\2 m_\nu(\sigma+\cI\1\eta)\2\dif \sigma \dif \eta\,,
\\
L_3(\nu)\,&:=\,\int_{\R^2} \eta \1 f'(\sigma)\1\chi'(\eta)\Re\2 m_\nu(\sigma+\cI\1\eta)\2\dif \sigma \dif \eta\,,
\end{split}
\end{equation*}
and $m_\nu$ is the Stieltjes transform of $\nu$.

Now we choose $f\geq 0$, such that $f|_{[\tau_1,\tau_2]}=1$ and $f|_{\R\backslash[\tau_1-\1\eta_1, \tau_2 + \1\eta_2]}=0$. Furthermore, we assume that the derivatives of $f$ satisfy
\bea{
&\norm{f'|_{[\1\tau_1-\1\eta_1,\2\tau_1]}}_\infty\,\lesssim\,\eta_1^{-1}\,,\qquad \norm{f''|_{[\1\tau_1-\1\eta_1,\2\tau_1]}}_\infty\,\lesssim\,\eta_1^{-2}
\,,
\\
& \norm{f'|_{[\1\tau_2,\2\tau_2+\1\eta_2]}}_\infty\,\lesssim\,\eta_2^{-1}\,,\qquad \norm{f''|_{[\1\tau_2,\2\tau_2+\1\eta_2]}}_\infty\,\lesssim\,\eta_2^{-2}
\,.
}
The function $\chi\geq0$ is chosen to be symmetric and such that $\chi|_{[-\eps,\2\eps\1]}=1$, $\chi|_{\R\backslash[-2\1\eps\1,2\1\eps\1]}=0$, as well as $\norm{\chi'}_\infty\lesssim \eps^{-1}$. Here the constant $\eps$ is chosen to satisfy $\eps\geq \max\{\eta_1,\eta_2\}$.
We now derive bounds on $L_{k}(\nu_1-\nu_2)$ for $k=1,2,3$.

We split the integral, $L_1$, into the contributions,
\begin{equation*}
\begin{split}
L_1(\nu)\,&=\, 2\2\bigl(\2L_{1,1,<}(\nu)+ L_{1,1,>}(\nu) + L_{1,2,<}(\nu)+ L_{1,2,>}(\nu)\2\bigr)\,,
\\
L_{1,1,<}(\nu)\,&:=\,  \int_{\tau_1-\eta_1}^{\tau_1} \msp{-8}\dif \sigma \int_{0}^{\eta_1} \dif \eta\,\eta \1 f''(\sigma)\1\Im\2 m_\nu(\sigma+\cI\1\eta)\,,
\\
L_{1,1,>}(\nu)\,&:=\,  \int_{\tau_1-\eta_1}^{\tau_1} \msp{-8}\dif \sigma \int_{\eta_1}^{2\1\eps} \dif \eta\,\eta \1 f''(\sigma)\1\chi(\eta)\1\Im\2 m_\nu(\sigma+\cI\1\eta)\,,
\\
L_{1,2,<}(\nu)\,&:=\,  \int_{\tau_2}^{\tau_2+\eta_2} \msp{-8}\dif \sigma \int_{0}^{\eta_2} \dif \eta\,\eta \1 f''(\sigma)\1\Im\2 m_\nu(\sigma+\cI\1\eta)\,,
\\
L_{1,2,>}(\nu)\,&:=\,  \int_{\tau_2}^{\tau_2+\eta_2} \msp{-8}\dif \sigma \int_{\eta_2}^{2\1\eps} \dif \eta\,\eta \1 f''(\sigma)\1\chi(\eta)\1\Im\2 m_\nu(\sigma+\cI\1\eta)\,.
\end{split}
\end{equation*}
For a positive measure $\nu$ the function $\eta \mapsto \eta \,\Im \2m_\nu(\sigma + \cI\1\eta)$ is monotonously increasing. Thus, we estimate
\begin{equation*}
\begin{split}
|L_{1,1,<}(\nu)|\,&\leq\max_{\sigma \in [0, \eta_1]}|f''(\tau_1-\sigma)|
\int_{\tau_1-\eta_1}^{\tau_1} \msp{-8}\dif \sigma \int_{0}^{\eta_1} \dif \eta\,\eta_1 \1\Im\2 m_\nu(\sigma+\cI\1\eta_1)
\\
\,&\leq
\int_{\tau_1-\eta_1}^{\tau_1} \msp{-8}\dif \sigma \, \Im\2 m_\nu(\sigma+\cI\1\eta_1)\,,\msp{80} \nu\2\geq\20\,.
\end{split}
\end{equation*}
We conclude that
\bels{I11< bound}{
|L_{1,1,<}(\nu_1-\nu_2)|\,&\leq\, 
\int_{\tau_1-\eta_1}^{\tau_1} \msp{-8}\dif \sigma\,\Bigl(\,2\2 \Im\2 m_{\nu_1}(\sigma+\cI\1\eta_1)+\absb{\1m_{\nu_1-\nu_2}(\sigma+\cI\1\eta_1)}\Bigr)
\,.
}
In the same way we find
\bels{I12< bound}{
|L_{1,2,<}(\nu_1-\nu_2)|\,&\leq\, 
\int_{\tau_2}^{\tau_2+\eta_2} \msp{-8}\dif \sigma\,\Bigl(2\2 \Im\2 m_{\nu_1}(\sigma+\cI\1\eta_2)+\absb{\1m_{\nu_1-\nu_2}(\sigma+\cI\1\eta_2)}\1\Bigr)\,.}
For the treatment of $L_{1,1,>}$ we integrate by parts, first in $\sigma$ and then in $\eta$,
\[
L_{1,1,>}(\nu)\,=\, 
-\eta_1\!\int_{\tau_1-\eta_1}^{\tau_1} \msp{-8}\dif \sigma \,  f'(\sigma)\1\Re\2  m_\nu(\sigma+\cI\1\eta_1)
-
 \int_{\eta_1}^{2\1\eps} \dif \eta\int_{\tau_1-\eta_1}^{\tau_1} \msp{-8}\dif \sigma \,\partial_\eta(\eta \1\chi(\eta))\1 f'(\sigma)\1\Re\2  m_\nu(\sigma+\cI\1\eta)\,.
\]
We use $\max_\eta|\chi(\eta)+ \eta\1\chi'(\eta)|\lesssim 1$ and $\max_{\sigma \in [\10\1,\2\eta_1]}|f'(\tau_1-\sigma)|\lesssim \eta_1^{-1} $. 
In this way we estimate for $\nu=\nu_1-\nu_2$,
\bels{I11> bound}{
L_{1,1,>}(\nu_1-\nu_2)\,\lesssim
\int_{\tau_1-\eta_1}^{\tau_1} \msp{-8}\dif \sigma \1  | m_{\nu_1-\nu_2}(\sigma+\cI\1\eta_1)|
+
\frac{1}{\eta_1} \int_{\eta_1}^{2\1\eps} \dif \eta\int_{\tau_1-\eta_1}^{\tau_1} \msp{-8}\dif \sigma \1|  m_{\nu_1-\nu_2}(\sigma+\cI\1\eta)|\,.
}
Going through the same steps we also arrive at
\bels{I12> bound}{
L_{1,2,>}(\nu_1-\nu_2)\,\lesssim
\int_{\tau_2}^{\tau_2+\eta_2} \msp{-8}\dif \sigma \1  | m_{\nu_1-\nu_2}(\sigma+\cI\1\eta_2)|
+
\frac{1}{\eta_2} \int_{\eta_2}^{2\1\eps} \dif \eta\int_{\tau_2}^{\tau_2+\eta_2} \msp{-8}\dif \sigma \1|  m_{\nu_1-\nu_2}(\sigma+\cI\1\eta)|\,.
}

We continue by estimating $L_2$ from above. 
\bels{I2 bound}{
|L_2(\nu_1-\nu_2)|\,\lesssim\,\frac{1}{\eps}\int_{\tau_1-\eta_1}^{\tau_2+\eta_2}\msp{-8}\dif \sigma \int_{\eps}^{2\1\eps}\dif \eta\1 | m_{\nu_1-\nu_2}(\sigma+\cI\1\eta)|\,.
}
Finally we derive a bound for $L_3$. We split the integral into two components,
\begin{equation*}
\begin{split}
L_3(\nu)\,&=\,2\2\big( L_{3,1}(\nu)+L_{3,2}(\nu) \big)\,,
\\
L_{3,1}(\nu)\,&:=\,
\int_{\tau_1-\eta_1}^{\tau_1}\msp{-8}\dif \sigma \int_{\eps}^{2\1\eps}\!\dif \eta\, \eta \1 f'(\sigma)\1\chi'(\eta)\Re\2 m_\nu(\sigma+\cI\1\eta)\,,
\\
L_{3,2}(\nu)\,&:=\,
\int_{\tau_2}^{\tau_2+\eta_2}\msp{-8}\dif \sigma \int_{\eps}^{2\1\eps}\!\dif \eta\, \eta \1 f'(\sigma)\1\chi'(\eta)\Re\2 m_\nu(\sigma+\cI\1\eta)\,.
\end{split}
\end{equation*}
We arrive at the bound
\[
L_3(\nu_1-\nu_2)\,\lesssim\, 
\frac{1}{\eta_1}\int_{\tau_1-\eta_1}^{\tau_1}\msp{-8}\dif \sigma \int_{\eps}^{2\1\eps}\!\dif \eta \1 | m_{\nu_1-\nu_2}(\sigma+\cI\1\eta)|
+
\frac{1}{\eta_2}\int_{\tau_2}^{\tau_2+\eta_2}\msp{-8}\dif \sigma \int_{\eps}^{2\1\eps}\!\dif \eta \1 | m_{\nu_1-\nu_2}(\sigma+\cI\1\eta)|
\,.
\]
We combine this with the estimates from \eqref{I11< bound}, \eqref{I12< bound}, \eqref{I11> bound}, \eqref{I12> bound} and \eqref{I2 bound}. Altogether we have
\[
\absB{\int f\,
\dif(\nu_1-\nu_2)}
\,\lesssim\, J_1 + J_2 +J_3\,,
\]
where the three terms on the right hand side are given by
\begin{equation*}
\begin{split}
J_1\,&:=\, 
\int_{\tau_1-\eta_1}^{\tau_1} \msp{-8}\dif \sigma\,\Big( \Im\2 m_{\nu_1}(\sigma+\cI\1\eta_1)+|m_{\nu_1-\nu_2}(\sigma+\cI\1\eta_1)|+\frac{1}{\eta_1} \int_{\eta_1}^{2\1\eps} \dif \eta\1|  m_{\nu_1-\nu_2}(\sigma+\cI\1\eta)|
\Big)\,,
\\
J_2\,&:=\, 
\int_{\tau_2}^{\tau_2+\eta_2} \msp{-8}\dif \sigma\,\Big( \Im\2 m_{\nu_1}(\sigma+\cI\1\eta_2)+|m_{\nu_1-\nu_2}(\sigma+\cI\1\eta_2)|+\frac{1}{\eta_2} \int_{\eta_2}^{2\1\eps} \dif \eta\1|  m_{\nu_1-\nu_2}(\sigma+\cI\1\eta)|
\Big)\,,
\\
J_3\,&:=\, 
\frac{1}{\eps}\int_{\tau_1-\eta_1}^{\tau_2+\eta_2}\msp{-8}\dif \sigma \int_{\eps}^{2\1\eps}\dif \eta\1 | m_{\nu_1-\nu_2}(\sigma+\cI\1\eta)|\,.
\end{split}
\end{equation*}

Now we use this bound for the smoothed out indicator function to derive a bound on the difference of number of eigenvalues in the interval $[\tau_1,\tau_2]$ and the predicted number, given by the integral over the density of states. We use
\bels{}{
\nu_2([\tau_1,\tau_2]) \,\leq\,\int f\2\dif \nu_1 + 
\int f \dif(\nu_1-\nu_2)
\,,
}
for $f$ defined as above. Then we get
\[
\nu_2([\tau_1,\tau_2]) \,\leq\,\nu_1([\tau_1,\tau_2])+\nu_1([\tau_1-\eta_1,\tau_1]\cup[\tau_2,\tau_2+\eta_2])
\,+\, 
\absB{\int f\2
\dif(\nu_1-\nu_2)}
\,.
\]
Similarly we use 
\[
\nu_1([\tau_1,\tau_2])\,\geq\, \int f\1\dif\nu_2 \,-\;\nu_1(\1[\tau_1-\eta_1,\tau_1]\cup[\tau_2,\tau_2+\eta_2]\1)\,,
\]
to get the bound
\[
\nu_1([\tau_1,\tau_2])
\;\ge\; 
\nu_2([\tau_1,\tau_2])
\,-\,
\absB{ \int f\2\dif(\nu_2-\nu_1)\1}\,-\,\nu_1(\1[\tau_1-\eta_1,\tau_1]\cup[\tau_2,\tau_2+\eta_2]\1)\,.
\]
Together, the two bounds imply
\[
\abs{\1\nu_1([\tau_1,\tau_2])-\nu_2([\tau_1,\tau_2])\1}
\;\lesssim\; 
\nu_1(\1[\tau_1-\eta_1,\tau_1]\cup[\tau_2,\tau_2+\eta_2]\1)+ J_1 +J_2+J_3
\,.
\]
\end{Proof}

\newpage

%

\end{document}